\pdfoutput=1
\RequirePackage{ifpdf}
\ifpdf 
\documentclass[pdftex]{sigma}
\else
\documentclass{sigma}
\fi

\usepackage[all]{xy}

\usepackage{array}
\usepackage{mathtools}
\usepackage{mathdots}
\usepackage{ytableau}

\usepackage{tikz}
\usetikzlibrary{patterns}
\usetikzlibrary{decorations.pathreplacing}
\usepackage{tkz-graph}
\usetikzlibrary{arrows}
\usetikzlibrary{decorations.markings}
\usepackage{cellspace}
\newtheorem{Theorem}{Theorem}[section]
\newtheorem{Corollary}[Theorem]{Corollary}
\newtheorem{Lemma}[Theorem]{Lemma}
\newtheorem{Proposition}[Theorem]{Proposition}
\newtheorem{Conjecture}[Theorem]{Conjecture}
\newtheorem{Question}[Theorem]{Question}
 { \theoremstyle{definition}
\newtheorem{Definition}[Theorem]{Definition}
\newtheorem{Example}[Theorem]{Example}
\newtheorem{Remark}[Theorem]{Remark} }

\numberwithin{equation}{section}

\usepackage{etoolbox}
\apptocmd{\sloppy}{\hbadness 10000\relax}{}{}

\def\SetFancyGraph {
	\SetVertexMath
	\GraphInit[vstyle=Art]
	\SetUpVertex[MinSize=2pt]
	\SetVertexLabel
	\tikzset{VertexStyle/.style = {shape = circle,shading = ball,ball color = black,inner sep = 1.5pt}}
	\SetUpEdge[color=black]
	\tikzset{->-/.style={decoration={ markings, mark=at position 0.8 with {\arrow{>}}},postaction={decorate}}}
	\tikzset{->--/.style={decoration={ markings, mark=at position 0.55 with {\arrow{>}}},postaction={decorate}}}
}

\newcommand{\widetriangle}{%
\begin{tikzpicture}%
\draw (-1.5ex,0) -- (1.5ex, 0) -- (0, 2ex) -- (-1.5ex,0);
\end{tikzpicture}%
}

\newcommand{\talltriangle}{%
\begin{tikzpicture}%
\draw (0,-1.5ex) -- (0,1.5ex) -- (2ex,0) -- (0,-1.5ex);
\end{tikzpicture}%
}

\newcommand{\cornertriangle}{%
\begin{tikzpicture}%
\draw (0,0) -- (2ex, 0) -- (0, 2ex) -- (0,0);
\end{tikzpicture}%
}

\begin{document}
\allowdisplaybreaks

\newcommand{\arXivNumber}{1907.09337}

\renewcommand{\PaperNumber}{130}

\FirstPageHeading

\ShortArticleName{Cyclic Sieving for Plane Partitions and Symmetry}

\ArticleName{Cyclic Sieving for Plane Partitions and Symmetry}

\Author{Sam HOPKINS}

\AuthorNameForHeading{S.~Hopkins}

\Address{School of Mathematics, University of Minnesota, Minneapolis, MN 55455, USA}
\Email{\href{mailto:shopkins@umn.edu}{shopkins@umn.edu}}
\URLaddress{\url{http://www-users.math.umn.edu/~shopkins/}}

\ArticleDates{Received May 05, 2020, in final form December 06, 2020; Published online December 09, 2020}

\Abstract{The cyclic sieving phenomenon of Reiner, Stanton, and White says that we can often count the fixed points of elements of a cyclic group acting on a combinatorial set by plugging roots of unity into a polynomial related to this set. One of the most impressive instances of the cyclic sieving phenomenon is a theorem of Rhoades asserting that the set of plane partitions in a rectangular box under the action of promotion exhibits cyclic sieving. In Rhoades's result the sieving polynomial is the size generating function for these plane partitions, which has a well-known product formula due to MacMahon. We extend Rhoades's result by also considering symmetries of plane partitions: specifically, complementation and transposition. The relevant polynomial here is the size generating function for symmetric plane partitions, whose product formula was conjectured by MacMahon and proved by Andrews and Macdonald. Finally, we explain how these symmetry results also apply to the rowmotion operator on plane partitions, which is closely related to promotion.}


\Classification{05E18; 05E10; 17B10; 17B37}

\section{Introduction and statement of results} \label{sec:intro}

\subsection{Plane partitions}

An \emph{$a\times b$ plane partition of height $m$} is an $a\times b$ array $\pi = (\pi_{i,j})_{\substack{1\leq i \leq a,\\ 1\leq j \leq b}}$ of nonnegative integers $\pi_{i,j} \in \mathbb{N}$ which is weakly decreasing in rows and columns (i.e., $\pi_{i,j} \geq \pi_{i+1,j}$ and~$\pi_{i,j} \geq \pi_{i,j+1}$ for all~$i$,~$j$) and for which the largest entry is less than or equal to~$m$ (i.e., $\pi_{1,1} \leq m$). We denote the set of such plane partitions by $\mathrm{PP}^{m}(a\times b)$. For a plane partition~$\pi \in \mathrm{PP}^{m}(a\times b)$, we define its \emph{size} to be $|\pi| \coloneqq \sum_{\substack{1\leq i \leq a, \\ 1\leq j \leq b}} \pi_{i,j}$.

MacMahon's celebrated product formula~\cite[Section~495]{macmahon1915combinatory} for the size generating function for $a\times b$ plane partitions of height $m$ is:
\[\mathrm{Mac}(a,b,m;q) \coloneqq \sum_{\pi \in \mathrm{PP}^{m}(a\times b)} q^{|\pi|} = \prod_{\substack{1 \leq i \leq a,\\ 1 \leq j \leq b}} \frac{\big(1-q^{i+j+m-1}\big)}{\big(1-q^{i+j-1}\big)}. \]
See~\cite[Theorem 7.21.7]{stanley1999ec2} for a modern presentation of this result. Note $\mathrm{Mac}(a,b,m;q)$ is (essentially) a principal specialization of a Schur polynomial $s_{\lambda}(x_1,x_2,\dots,x_{a+b})$:
\[\mathrm{Mac}(a,b,m;q) = q^{-\kappa(m^a)} s_{m^a}\big(1,q,q^2,\dots,q^{a+b-1}\big),\]
where $m^a$ is the $a\times m$ rectangle shape, and $\kappa(\lambda) \coloneqq 0\lambda_1 + 1\lambda_2 + 2\lambda_3 + \cdots$. The Schur polynomials~$s_{\lambda}$ occur in many contexts, but of particular relevance is the fact that they are characters of general linear group representations.

\subsection{Promotion}

\emph{Promotion} is a certain invertible operator acting on the set of these plane partitions. It can be defined as a composition of piecewise-linear involutions, as follows. We define the \emph{piecewise-linear toggle} $\tau_{i,j}\colon \mathrm{PP}^{m}(a\times b)\to\mathrm{PP}^{m}(a\times b)$ for $1 \leq i \leq a$, $1\leq j\leq b$ by
\[(\tau_{i,j}\pi)_{k,l} \coloneqq \begin{cases}\pi_{k,l} &\textrm{if $(k,l)\neq (i,j)$}, \\ \min (\pi_{i,j-1},\pi_{i-1,j}) + \max (\pi_{i+1,j},\pi_{i,j+1}) - \pi_{i,j} &\textrm{if $(k,l)=(i,j)$}, \\ \end{cases}\]
with the convention that $\pi_{0,j} \coloneqq \pi_{i,0} \coloneqq m$ and $\pi_{a+1,j} \coloneqq \pi_{i,b+1} \coloneqq 0$. The $\tau_{i,j}$ are involutions. Also, toggles $\tau_{i,j}$ and $\tau_{i',j'}$ commute unless $(i,j)$ and $(i',j')$ are directly adjacent, where by ``directly adjacent'' we mean that $|i-i'|+|j-j'|=1$.

Then for~$-a+1 \leq k \leq b-1$ we define $\mathcal{F}_{k} \coloneqq \prod_{\substack{1\leq i \leq a, \\1\leq j \leq b, \\ j-i=k}} \tau_{i,j}$ to be the composition of all the toggles along the ``$k$th diagonal'' of our array (note that all these toggles commute).

Finally, we define \emph{promotion} $\mathrm{Pro}\colon \mathrm{PP}^{m}(a\times b)\to\mathrm{PP}^{m}(a\times b)$ as the composition of these diagonal toggles $\mathcal{F}_k$ from left to right:
\[\mathrm{Pro} \coloneqq \mathcal{F}_{b-1} \cdot \mathcal{F}_{b-2} \cdots \mathcal{F}_{-a+2} \cdot \mathcal{F}_{-a+1}.\]

\begin{Example}
Suppose $a \coloneqq 2$, $b\coloneqq 2$, and $m \coloneqq 4$. We can compute an application of promotion on a plane partition $\pi \in \mathrm{PP}^{4}(2\times 2)$ as follows:
\[\pi = \begin{ytableau}
2 & 2\\
*(lightgray) 1 & 0
\end{ytableau} \xrightarrow{\tau_{2,1}} \begin{ytableau}
2 & 2\\
 1 & *(lightgray) 0
\end{ytableau} \xrightarrow{\tau_{2,2}} \begin{ytableau}
*(lightgray) 2 & 2\\
 1 & 1
\end{ytableau} \xrightarrow{\tau_{1,1}} \begin{ytableau}
 4 & *(lightgray) 2\\
 1 & 1
\end{ytableau} \xrightarrow{\tau_{1,2}} \begin{ytableau}
 4 & 3\\
 1 & 1
\end{ytableau} = \mathrm{Pro}(\pi).\]
\end{Example}

This description of promotion in terms of piecewise-linear involutions goes back to Berenstein and Kirillov~\cite{kirillov1995groups} and Berenstein and Zelevinsky~\cite{berenstein1996canonical}, building on work of Bender and Knuth~\cite{bender1972enumeration} and Gansner~\cite{gansner1980equality}. More recently, interest in these piecewise-linear toggles has been rekindled in connection with another related operator called \emph{rowmotion} (see, e.g.,~\cite{einstein2013combinatorial} or the survey~\cite{roby2016dynamical}). We will discuss rowmotion later.

It is more common, following the seminal work of Sch\"{u}tzenberger~\cite{schuztenberger1963quelques, schutzenberger1972promotion, schutzenberger1977robinson} (see also Haiman~\cite{haiman1992dual}), to consider promotion as an operator on semistandard Young tableaux defined in terms of ``jeu de taquin'' sliding moves. But in fact, via a simple change of coordinates using Gelfand--Tsetlin patterns, promotion of plane partitions in $\mathrm{PP}^{m}(a\times b)$ as we have just defined it exactly corresponds to the usual promotion of semistandard tableaux of shape $a \times m$ with entries in the set $\{1,\dots,a+b\}$ (see Appendix~\ref{sec:appendix} for the details of this correspondence).

\subsection{The sieving phenomenon}

The \emph{sieving phenomenon} of Reiner, Stanton, and White~\cite{reiner2004cyclic} is, loosely speaking, the philosophy that we can often count fixed points for a nice group action on a set of combinatorial objects by plugging roots of unity into a polynomial related to this set. Initially the philosophy was considered only for cyclic group actions, and in this context it is usually called the \emph{cyclic sieving phenomenon}; but there has also been interest in broadening the philosophy to include other groups as well~\cite{barcelo2008bimahonian, rao2017dihedral}. In fact, as we explain below, the sieving phenomenon grew out of Stembridge's \emph{``$q=-1$'' phenomenon}, which is basically the case where the group has order two.

Sieving phenomena involving polynomials which have simple product formulas in terms of ratios of $q$-numbers (like the MacMahon formula) are especially valuable, because these imply that every symmetry class has a product formula.

One of the most impressive such examples of the cyclic sieving phenomenon is:

\begin{Theorem}[{\cite[Theorem 1.4]{rhoades2010cyclic}, \cite[Theorem 1.3]{shen2018cyclic}}] \label{thm:ppart_csiev}
For any $k \in \mathbb{Z}$, we have
\[\#\big\{\pi \in \mathrm{PP}^{m}(a\times b)\colon \mathrm{Pro}^k(\pi)=\pi\big\} = \mathrm{Mac}\big(a,b,m;q \mapsto \zeta^k\big),\]
where $\zeta \coloneqq {\rm e}^{2\pi {\rm i}/(a+b)}$ is a primitive $(a+b)$th root of unity.
\end{Theorem}

Theorem~\ref{thm:ppart_csiev} says that promotion acting on plane partitions has a very regular orbit structure. For instance, note that one consequence of Theorem~\ref{thm:ppart_csiev} is that $\mathrm{Pro}$ acting on $\mathrm{PP}^{m}(a\times b)$ has order dividing $a+b$. But, as mentioned, Theorem~\ref{thm:ppart_csiev} implies much more than this: it also means that every symmetry class has a product formula.

Theorem~\ref{thm:ppart_csiev} was first proved by Rhoades~\cite{rhoades2010cyclic}. To prove this theorem he used Kazhdan--Lusztig theory and the related theory of quantum groups. In particular, he employed the dual canonical basis for representations of the general linear group~\cite{kashiwara1993global, lusztig1990canonical}. More recently, Shen and Weng~\cite{shen2018cyclic} gave a different proof of Theorem~\ref{thm:ppart_csiev}. Their approach was from the perspective of cluster algebras and the ``cluster duality'' conjecture of Fock and Goncharov~\cite{fock2009cluster}. Specifically, their proof employed the Gross--Hacking--Keel--Kontsevich~\cite{gross2018canonical} canonical basis (or ``theta basis'') for cluster algebras, in the particular case of the coordinate ring of the Grassmannian.

In either case, the proof of Theorem~\ref{thm:ppart_csiev} followed the ``linear algebraic'' paradigm. This means that the desired equality is established by computing the trace of a linear operator on a vector space in two different bases. The first basis should be indexed by the combinatorial set in question, and the linear operator should permute this basis according to the cyclic action, so that its trace computes the number of fixed points of the action. The second basis is an eigenbasis where we can compute trace by considering eigenvalues. See, e.g.,~\cite{reiner2014whatis} or~\cite[Section~4]{sagan2011cyclic}.

In both the Rhoades~\cite{rhoades2010cyclic} and Shen--Weng~\cite{shen2018cyclic} proofs the vector space in question actually carries more structure: it is a ${\rm GL}(a+b)$ representation. And in both proofs the linear operator corresponding to promotion is the action of a particular lift to the general linear group of the~\emph{long cycle} (a.k.a.~standard \emph{Coxeter element}) $c\in \mathfrak{S}_{a+b}$ in the symmetric group. Geometrically, this map is the \emph{twisted cyclic shift}~$\chi \in {\rm GL}(a+b)$ acting on the Grassmannian $\mathrm{Gr}(a,a+b)$.

\subsection{Symmetries of plane partitions}

In the present paper, we extend Theorem~\ref{thm:ppart_csiev} by considering symmetries of plane partitions. The study of plane partitions with symmetry goes back to MacMahon~\cite{macmahon1899partitions}, but really took off in the 1970s and 80s: see for instance the seminal paper of Stanley~\cite{stanley1986symmetries} which identified 10 symmetry classes of plane partitions, and see~\cite{krattenthaler2016plane} for a modern update to Stanley's paper. Here we will be concerned exclusively with the involutive symmetries of plane partitions.

The first symmetry we consider is \emph{complementation} $\mathrm{Co}\colon \mathrm{PP}^{m}(a\times b)\to \mathrm{PP}^{m}(a\times b)$, which is defined by $\mathrm{Co}(\pi)_{i,j} \coloneqq m - \pi_{a+1-i,b+1-j}$. A plane partition $\pi \in \mathrm{PP}^{m}(a\times b)$ can be viewed as the $3$-dimensional stack of cubes inside of an $a\times b \times m$ rectangular box which has~$\pi_{i,j}$ cubes stacked at position $(i,j)$, in which case complementation is set-theoretic complementation inside this box.

Plane partitions $\pi \in \mathrm{PP}^{m}(a\times b)$ with $\mathrm{Co}(\pi)=\pi$ are called \emph{self-complementary}. Stanley~\cite{stanley1986symmetries} was the first to enumerate self-complementary plane partitions. The enumeration of self-complementary plane partitions is one of the prototypical examples of Stembridge's \emph{``$q=-1$'' phenomenon}~\cite{stembridge1994minuscule, stembridge1996canonical}, the precursor to the cyclic sieving phenomenon. Namely:

\begin{Theorem}[{\cite[Theorem~1.1]{stembridge1994hidden}, \cite[Theorem~4.1]{stembridge1994minuscule}, \cite[Theorem~3.1]{kuperberg1994minuscule}}] \label{thm:com_qneg}
We have
\[\#\big\{\pi \in \mathrm{PP}^{m}(a\times b)\colon \mathrm{Co}(\pi)=\pi\big\} = \mathrm{Mac}(a,b,m;q \mapsto -1).\]
\end{Theorem}

Stembridge~\cite{stembridge1994minuscule} (and, independently, Kuperberg~\cite{kuperberg1994minuscule}) gave a ``linear algebraic'' proof of Theo\-rem~\ref{thm:com_qneg}, computing the trace of a linear operator in two ways. The linear operator corresponding to complementation is the action of a particular lift $\overline{w_0}\in {\rm GL}(a+b)$ of the \emph{longest element}~$w_0 \in \mathfrak{S}_{a+b}$ in the symmetric group. Geometrically, this map is the \emph{twisted reflection} on the Grassmannian.

Complementation and promotion together generate a dihedral group: $\mathrm{Co} \cdot \mathrm{Pro} = \mathrm{Pro}^{-1} \cdot \mathrm{Co}$. (In the context of tableaux, complementation is more often referred to as the \emph{``Sch\"{u}tzenberger involution''} or \emph{``evacuation''}; see for instance Stanley's survey paper~\cite{stanley2009promotion}.) This means that $\mathrm{Co} \cdot \mathrm{Pro}^k$ is conjugate to $\mathrm{Co} \cdot \mathrm{Pro}^j$ whenever $k$ and $j$ have the same parity. So if we want to count the number of fixed points of $\mathrm{Co} \cdot \mathrm{Pro}^k$, there are only two cases to consider: $k$ even (which is addressed by Theorem~\ref{thm:com_qneg} above), and $k$ odd.

Recently there has been interest in sieving for dihedral group actions~\cite{rao2017dihedral, stier2020dihedral}, where fixed points of both the rotations and reflections are counted by plugging roots of unity into a polynomial. And in fact a ``dihedral sieving'' result for the action of $\langle \mathrm{Pro}, \mathrm{Co}\rangle$ on plane partitions has already been considered by Abuzzahab--Korson--Li--Meyer~\cite{abuzzahab2005cyclic} and obtained by Rhoades~\cite{rhoades2010cyclic}:\footnote{We stated this result in a somewhat more compact way than it is stated in \cite[Theorem~7.6]{rhoades2010cyclic}.}

\begin{Theorem}[{\cite[Theorem 1.2 and Conjecture~1.3]{abuzzahab2005cyclic},~\cite[Theorem 7.6]{rhoades2010cyclic}}] \label{thm:com_csiev}
 For any even~$k\in \mathbb{Z}$, we have $\#\big\{\pi \in \mathrm{PP}^{m}(a\times b)\colon \mathrm{Co}\cdot \mathrm{Pro}^k(\pi)=\pi\big\} = \mathrm{Mac}(a,b,m;q \mapsto -1)$.
For any odd $k \in \mathbb{Z}$, we have
\[\#\left\{ \begin{matrix}\pi \in \mathrm{PP}^{m}(a\times b) \colon \\ \mathrm{Co}\cdot \mathrm{Pro}^k(\pi)=\pi \end{matrix} \right\} = \begin{cases} \mathrm{Mac}(a,b,m;q \mapsto -1) &\textrm{if $a$ or $b$ is even}, \\ (-1)^{\kappa(m^a)} s_{m^a}(\overbrace{1,-1,1,\dots,-1,1}^{\textrm{$a+b-1$ values}},1) &\textrm{if $a$ and $b$ are odd}. \end{cases}\]
\end{Theorem}

We note that there are product formulas for the Schur function evaluation appearing in Theorem~\ref{thm:com_csiev}: see~\cite[Lemma 8.2]{abuzzahab2005cyclic}. However, Theorem~\ref{thm:com_csiev} is certainly not quite as clean a result as Theorem~\ref{thm:ppart_csiev}. By considering additional plane partition symmetries, we will actually discover fixed point enumerations which are as pleasant as Theorem~\ref{thm:ppart_csiev}.

The next symmetry we consider is \emph{transposition}, or in other words, reflection across the main diagonal. In order for transposition to act on a fixed set of plane partitions, we need $a=b$. So let $n\coloneqq a=b$. Then $\mathrm{Tp}\colon \mathrm{PP}^{m}(n\times n)\to \mathrm{PP}^{m}(n\times n)$ is defined by $\mathrm{Tp}(\pi)_{i,j} \coloneqq \pi_{j,i}$. Plane partitions $\pi \in \mathrm{PP}^{m}(n\times n)$ with $\mathrm{Tp}(\pi)=\pi$ are usually just called \emph{symmetric plane partitions}. In 1899 MacMahon conjectured~\cite{macmahon1899partitions}, and in 1978 Andrews~\cite{andrews1978plane} proved, the following product formula for the size generating function for symmetric plane partitions:
\[\mathrm{SymMac}(n,m;q) \coloneqq \sum_{\substack{\pi \in \mathrm{PP}^{m}(n\times n), \\ \mathrm{Tp}(\pi)=\pi}} q^{|\pi|} = \prod_{1\leq i < j \leq n}\frac{\big(1-q^{2(i+j+m-1)}\big)}{\big(1-q^{2(i+j-1)}\big)} \cdot \prod_{1\leq i \leq n} \frac{\big(1-q^{2i+m-1}\big)}{\big(1-q^{2i-1}\big)}. \]
Macdonald~\cite[Example~17, p.~52]{macdonald1979symmetric} also proved MacMahon's conjecture independently and essentially simultaneously. See also~\cite{gordon1983bender} and~\cite{proctor1984bruhat} for other proofs of MacMahon's conjecture.

Transposition and complementation commute, and so together generate a group isomorphic to $\mathbb{Z}/2\mathbb{Z}\times \mathbb{Z}/2\mathbb{Z}$. We will refer to their composition $\mathrm{Tp}\cdot \mathrm{Co}$ as \emph{transpose-complementation}. Plane partitions $\pi \in \mathrm{PP}^{m}(n\times n)$ with $\mathrm{Tp}\cdot \mathrm{Co}(\pi) = \pi$ are called \emph{transpose-complementary}. Proctor~\cite{proctor1990new} was the first to count transpose-complementary plane partitions. Later, Stembridge observed~\cite{stembridge1994hidden}, and Kuperberg~\cite{kuperberg1994minuscule} explained, the following $q=-1$ phenomenon for the enumeration of transpose-complementary plane partitions:

\begin{Theorem}[{\cite[Theorem 1.2]{stembridge1994hidden},~\cite[Theorem 5.1]{kuperberg1994minuscule}}] \label{thm:trcom_qneg}
We have
\begin{gather*} \#\big\{\pi \in \mathrm{PP}^{m}(n\times n)\colon \mathrm{Tp}\cdot \mathrm{Co}(\pi)=\pi\big\} = \mathrm{SymMac}(n,m;q \mapsto -1).\end{gather*}
\end{Theorem}

Kuperberg~\cite{kuperberg1994minuscule} gave a linear algebraic proof of Theorem~\ref{thm:trcom_qneg} where transposition corresponds to the outer automorphism of ${\rm GL}(2n)$ induced by the symmetry of the Dynkin diagram of type~$A_{2n-1}$. Geometrically, the outer automorphism is the \emph{symplectic orthogonal complement} on the Grassmannian.

There is another $q$-analog for the symmetric plane partitions. For $\pi \in \mathrm{PP}^{m}(n\times n)$, set $|\pi|' \coloneqq \sum_{1\leq i \leq j \leq n} \pi_{i,j}$, which roughly speaking is the `size of~$\pi/\mathrm{Tp}$'. Bender and Knuth~\cite{bender1972enumeration} conjectured that:
\[\mathrm{SymMac}'(n,m;q) \coloneqq \sum_{\substack{\pi \in \mathrm{PP}^{m}(n\times n), \\ \mathrm{Tp}(\pi)=\pi}} q^{|\pi|'} = \prod_{1\leq i \leq j \leq n}\frac{\big(1-q^{i+j+m-1}\big)}{\big(1-q^{i+j-1}\big)}. \]
In~\cite{andrews1977plane}, Andrews showed that in fact the MacMahon and Bender--Knuth conjectures are equivalent, and in doing so proved the Bender--Knuth conjecture.

This second $q$-analog of symmetric plane partitions also has an associated $q=-1$ result, this time enumerating \emph{symmetric self-complementary plane partitions}. Symmetric self-complementary plane partitions were again first enumerated by Proctor~\cite{proctor1983trapezoid}. Later, Stembridge~\cite{stembridge1994minuscule, stembridge1994hidden} interpreted this enumeration as a $q=-1$ result:

\begin{Theorem}[{\cite[Theorem~4.1, Example~4.3]{stembridge1994minuscule}, \cite[Theorem~1.1]{stembridge1994hidden}}] \label{thm:sym_com_qneg}
We have
\[\#\big\{\pi \in \mathrm{PP}^{m}(n\times n)\colon \mathrm{Tp}(\pi)=\pi, \mathrm{Co}(\pi)=\pi\big\} = \mathrm{SymMac}'(n,m;q \mapsto -1).\]
\end{Theorem}

Stembridge~\cite{stembridge1994minuscule} gave a linear algebraic proof of Theorem~\ref{thm:sym_com_qneg}, but this time using representations of the \emph{special orthogonal group} ${\rm SO}(2n+1)$ (or at least its Lie algebra) instead of the general linear group.

\subsection{New sieving results}

To summarize the above, the interaction of promotion and complementation of plane partitions, and also the interaction of transposition and complementation, are understood. The main undertaking of this paper is to understand how promotion and transposition, and promotion and transpose-complementation, interact.

Transposition and promotion together generate a dihedral group: $\mathrm{Tp} \cdot \mathrm{Pro} = \mathrm{Pro}^{-1} \cdot \mathrm{Tp}$. Our first main result is the following ``dihedral sieving''-style result concerning fixed points for elements of $\langle \mathrm{Pro}, \mathrm{Tp}\rangle$:

\begin{Theorem} \label{thm:tr_csiev}
For any $k \in \mathbb{Z}$, we have
\[\#\big\{\pi \in \mathrm{PP}^{m}(n\times n)\colon \mathrm{Tp}\cdot \mathrm{Pro}^k(\pi)=\pi\big\} = \mathrm{SymMac}\big(n,m;q \mapsto (-1)^k\big).\]
\end{Theorem}

To prove Theorem~\ref{thm:tr_csiev}, we note that the evaluation $\mathrm{SymMac}(n,m;q \mapsto -1)$ is nonzero if and only if $m=2M$ is even, and in this case is equal to
\[\mathrm{SymMac}(n,2M;q \mapsto -1) = \prod_{1\leq i \leq j \leq n-1} \frac{i+j+2M}{i+j}.\]
Using algebraic techniques, Proctor~\cite{proctor1990new} demonstrated that a couple of different combinatorial sets of plane partition flavor are enumerated by~$\prod_{1\leq i \leq j \leq n-1} \frac{i+j+2M}{i+j}$. We will establish that the set $\big\{\pi \in \mathrm{PP}^{2M}(n\times n)\colon \mathrm{Tp}\cdot \mathrm{Pro}(\pi)=\pi\big\}$ is in bijection with one of the sets Proctor showed is counted by~$\prod_{1\leq i \leq j \leq n-1} \frac{i+j+2M}{i+j}$.

Transpose-complementation and promotion commute. Hence the group they generate is a~product of two cyclic groups: $\langle \mathrm{Pro}, \mathrm{Tp}\cdot \mathrm{Co}\rangle\simeq \mathbb{Z}/2n\mathbb{Z} \times \mathbb{Z}/2\mathbb{Z}$. Enumerating fixed points for elements of a~cyclic group times $\mathbb{Z}/2\mathbb{Z}$ is a bit more interesting than for elements of a dihedral group because there are more conjugacy classes. Barcelo, Reiner, and Stanton~\cite{barcelo2008bimahonian} considered an extension of cyclic sieving to products of two cyclic groups, which they called ``bicyclic sieving''. Our second main result is the following ``bicyclic sieving''-style result concerning fixed points for elements of~$\langle \mathrm{Pro}, \mathrm{Tp}\cdot \mathrm{Co}\rangle$:

\begin{Theorem} \label{thm:trcom_csiev}
For any $k \in \mathbb{Z}$, we have
\[\#\big\{\pi \in \mathrm{PP}^{m}(n\times n)\colon (\mathrm{Tp}\cdot\mathrm{Co})\cdot \mathrm{Pro}^{n+k}(\pi)=\pi\big\} = \mathrm{SymMac}\big(n,m;q \mapsto \zeta^k\big),\]
where $\zeta \coloneqq {\rm e}^{\pi {\rm i}/n}$ is a primitive $(2n)$th root of unity.
\end{Theorem}

The proof of Theorem~\ref{thm:trcom_csiev} is more involved than the proof of Theorem~\ref{thm:tr_csiev}. We use a linear algebraic approach, extending the work of Rhoades~\cite{rhoades2010cyclic}. Basically, we show that $\overline{w_0} \in {\rm GL}(2n)$ (corresponding to complementation) and the outer automorphism of ${\rm GL}(2n)$ (corresponding to transposition) behave in the appropriate way on the dual canonical basis of the relevant ${\rm GL}(2n)$ representation. In fact, this has essentially already been done: Berenstein--Zelevinsky~\cite{berenstein1996canonical} and Stembridge~\cite{stembridge1996canonical} showed that~$w_0$ behaves as evacuation on the dual canonical basis of any irreducible general linear group representation; and Berenstein--Zelevinsky~\cite{berenstein1996canonical} also described the effect of the outer automorphism on the dual canonical basis. (General results of Lusztig~\cite{lusztig1993quantum} imply that these automorphisms permute the canonical basis in \emph{some} way.) We just have to put all these results together and compute the trace of the appropriate composition of these operators.

Theorems~\ref{thm:ppart_csiev}, \ref{thm:com_csiev}, \ref{thm:tr_csiev} and~\ref{thm:trcom_csiev} together imply that for any element $g \in \langle \mathrm{Pro}, \mathrm{Co}, \mathrm{Tp}\rangle$, the number of plane partitions in $\mathrm{PP}^{m}(n\times n)$ fixed by $g$ is given by some kind of evaluation at a root of unity of a polynomial which has a nice product formula representation as a rational expression. However, it is unclear how to package all of these results together into one theorem.

\begin{Example} \label{ex:m1}
In this example we consider the case $m=1$.

For a subset $I\subseteq \{1,2,\dots,a+b\}$ and a permutation $\sigma \in \mathfrak{S}_{a+b}$ in the symmetric group on $a+b$ letters, we write $\sigma(I) \coloneqq \{\sigma(i)\colon i \in I\}$. We also use the notation $-I \coloneqq \{1,2,\dots,a+b\} \setminus I$. For us, two significant elements of $\mathfrak{S}_{a+b}$ are the \emph{long cycle} $c \coloneqq \left(\begin{smallmatrix} 1 & 2 & \dots & a+b-1 & a+b \\ 2 & 3 & \dots & a+b & 1 \end{smallmatrix}\right)$ and the \emph{longest element} $w_0 \coloneqq \left(\begin{smallmatrix} 1 & 2 & \dots & a+b \\ a+b & a+b-1 & \dots & 1 \end{smallmatrix}\right)$.

We define a bijection $\Psi\colon \mathrm{PP}^{1}(a\times b) \xrightarrow{\sim} \{I\subseteq \{1,\dots,a+b\} \textrm{ of size $a$}\}$ in Appendix~\ref{sec:appendix}. Under this bijection we have $\Psi(\mathrm{Pro}(\pi)) = c(\Psi(\pi))$, $\Psi(\mathrm{Co}(\pi)) = w_0(\Psi(\pi))$, and (in the case~$a=b=n$) $\Psi(\mathrm{Tp}(\pi))=-w_0(\Psi(\pi))$ for all $\pi \in \mathrm{PP}^{1}(a\times b)$. Hence, when $m=1$, the results described above reduce to formulas counting the fixed points of rotation, reversal, and complementation acting on subsets.

For instance, Theorem~\ref{thm:ppart_csiev} says that
\[ \# \big\{ I\subseteq \{1,\dots,a+b\} \textrm{ of size $a$}\colon c^k(I) = I\big\} = \mathrm{Mac}\big(a,b,1;q \mapsto \zeta^k\big), \]
where $\zeta \coloneqq {\rm e}^{2\pi {\rm i}/(a+b)}$ is a primitive $(a+b)$th root of unity, and $\mathrm{Mac}(a,b,1;q)$ is the usual \emph{$q$-binomial coefficient}:
\[\mathrm{Mac}(a,b,1;q) = \genfrac{\lbrack}{\rbrack}{0pt}{}{a+b}{a}_{q} \coloneqq q^{-a(a+1)/2} \cdot \sum_{\substack{I \subseteq \{1,\dots,a+b\}, \\ \# I = a}} q^{\sum_{i\in I} i} = \prod_{1\leq i \leq a} \frac{\big(1-q^{a+b+1-i}\big)}{\big(1-q^{i}\big)}.\]
This is one of the most prototypical cyclic sieving results, going back to the original Reiner--Stanton--White paper~\cite[Theorem~1.1(b)]{reiner2004cyclic}. Theorem~\ref{thm:com_csiev} offers a dihedral extension of this prototypical cyclic sieving result. This dihedral extension, which combines rotation of subsets with reversal of subsets, is less well known, but is discussed for instance in~\cite[Proposition~4.1]{rao2017dihedral} (at least for $a+b$ odd).

Now let us assume $a=b=n$ so that we can also consider what Theorems~\ref{thm:tr_csiev} and~\ref{thm:trcom_csiev} say about the case $m=1$. Actually, Theorem~\ref{thm:tr_csiev} is not so interesting in this case because, as mentioned, $\mathrm{SymMac}(n,m;q \mapsto -1) = 0$ when $m$ is odd, and moreover it is clear that there are no~$I$ with $-w_0(I)=c(I)$ since $1 \in c(I) \iff n \in I$ while $1 \in -w_0(I) \iff n \notin I$. But Theorem~\ref{thm:trcom_csiev} is quite interesting in this case. It says that
\[ \# \big\{ I\subseteq \{1,\dots,2n\} \textrm{ of size $n$}\colon c^{n+k}(I) = -I\big\} = \mathrm{SymMac}\big(n,1;q \mapsto \zeta^k\big),\]
where $\zeta \coloneqq {\rm e}^{\pi {\rm i}/n}$ is a primitive $(2n)$th root of unity, and $ \mathrm{SymMac}(n,1;q)$ has the simple form
\[ \mathrm{SymMac}(n,1;q) = \sum_{I \subseteq \{1,3,5,\dots,2n-1\}} q^{\sum_{i\in I} i} = \prod_{1\leq i \leq n} \big(1+q^{2i-1}\big).\]
For example, taking $n=2$, we have
\[\mathrm{SymMac}(2,1;q) = 1 + q + q^3 + q^4.\]
The relevant evaluations are
\[\mathrm{SymMac}(2,1;q \mapsto 1) = 4, \qquad \mathrm{SymMac}(2,1;q \mapsto \pm i) = 2, \qquad \mathrm{SymMac}(2,1;q \mapsto -1) = 0.\]
In agreement with these evaluations: the $I \subseteq \{1,2,3,4\}$ of size $2$ with $c^2(I) = -I$ are $\{1,2\}$, $\{2,3\}$, $\{3,4\}$, and $\{1,4\}$; while the $I$ for which $c(I) = -I$ (and also the ones with $c^3(I)=-I$) are $\{1,3\}$ and $\{2,4\}$; and there are no $I$ with $I=-I$.

It is worth contrasting the previous paragraph with a known ``type~B'' version of cyclic sieving for subsets under rotation. For a subset $I\subseteq \{1,2,\dots,n\}$, we define~$\widetilde{c}(I) \coloneqq c(I) \Delta \{1\}$, where here $c \coloneqq \left(\begin{smallmatrix} 1 & 2 & \dots &n-1 & n \\ 2 & 3 & \dots & n & 1 \end{smallmatrix}\right) \in \mathfrak{S}_{n}$, and $\Delta$ denotes symmetric difference. This $\widetilde{c}$ action might be called \emph{twisted rotation}, and its order is $2n$. We have
\[ \# \big\{ I\subseteq \{1,\dots,n\} \colon \widetilde{c}^{\, k}(I) = I\big\} = \mathrm{SymMac}'\big(n,1;q \mapsto \zeta^k\big),\]
where $\zeta \coloneqq {\rm e}^{\pi {\rm i}/n}$ is a primitive $(2n)$th root of unity, and $\mathrm{SymMac}'(n,1;q)$ has the simple form:
\[ \mathrm{SymMac}'(n,1;q) = \sum_{I \subseteq \{1,2,\dots,n\}} q^{\sum_{i\in I} i} = \prod_{1\leq i \leq n} \big(1+q^{i}\big).\]
This result also appears in the original Reiner--Stanton--White paper~\cite[Corollary~8.5]{reiner2004cyclic}, and is further discussed in~\cite[Section~6.1]{alexandersson2019cyclic}. Work of Rush and Shi~\cite{rush2013orbits} implies that subsets under twisted rotation are in equivariant bijection with $\mathrm{PP}^{1}(\talltriangle_{n})$ under $\mathrm{Row}$ (these notions are defined in Section~\ref{sec:row_tri}); hence, this result can be seen as the case $m=1$ of Conjecture~\ref{conj:tall_csiev} as well. While evidently quite similar, we know of no direct connection between counting subsets fixed by powers of twisted rotation, and counting subsets whose rotations are equal to their complements.
\end{Example}

\begin{Remark}Some of the root of unity evaluations of polynomials appearing in Rhoades's paper~\cite{rhoades2010cyclic} had prior combinatorial interpretations, for instance in terms of border-strip tableaux (see~\cite{lascoux1994green}). However, we are not aware of any prior combinatorial interpretations of the evaluations appearing in~Theorems~\ref{thm:tr_csiev} and~\ref{thm:trcom_csiev}.
\end{Remark}

\subsection{Rowmotion}

In the last sections of the paper we consider another invertible operator on plane partitions called \emph{rowmotion}, $\mathrm{Row}\colon \mathrm{PP}^m(a\times b)\to \mathrm{PP}^m(a\times b)$. Rowmotion and promotion are closely related. Rowmotion is, like promotion, a composition of all of the piecewise-linear toggles acting on~$\mathrm{PP}^m(a\times b)$. However, whereas promotion is a composition of these toggles ``from left to right'', rowmotion is a composition ``from top to bottom''. Striker and Williams~\cite{striker2012promotion} explained that the actions of promotion and rowmotion are conjugate; in particular, there is some composition $D$ of toggles so that $D \cdot \mathrm{Row} \cdot D^{-1} = \mathrm{Pro}$. We show that this conjugating map $D$ behaves nicely with respect to complementation and transposition. We conclude that versions of Theorems~\ref{thm:ppart_csiev}, \ref{thm:com_csiev}, \ref{thm:tr_csiev} and~\ref{thm:trcom_csiev} hold for $\mathrm{Row}$ (but with slight differences since, e.g., $\mathrm{Row}$ commutes with $\mathrm{Tp}$, while $(\mathrm{Tp}\cdot \mathrm{Co})\cdot \mathrm{Row} = \mathrm{Row}^{-1} \cdot (\mathrm{Tp}\cdot \mathrm{Co})$, et cetera). In particular, for any $g \in \langle \mathrm{Row}, \mathrm{Co}, \mathrm{Tp}\rangle$, we can again count the number of plane partitions in $\mathrm{PP}^{m}(n\times n)$ fixed by $g$ by some kind of sieving phenomenon evaluation.

One reason to consider rowmotion instead of promotion is because rowmotion makes sense acting on any partially ordered set (not all posets have a notion of left and right, but they all have a notion of top and bottom). Our original motivation for studying the way symmetries interact with rowmotion was a series of cyclic sieving conjectures we made in~\cite{hopkins2019minuscule} concerning rowmotion acting on the $P$-partitions of other posets $P$ besides the rectangle poset. Many of the posets with conjectured cyclic sieving for rowmotion are ``triangular'' posets which can be obtained from the rectangle by enforcing certain symmetries. More precisely, in the final section we show, following Grinberg and Roby~\cite{grinberg2015birational2}, that the $P$-partitions for these various triangular posets $P$ are in rowmotion-equivariant bijection with the set of plane partitions in $\mathrm{PP}^{m}(n\times n)$ fixed by various subgroups of~$\langle \mathrm{Row}, \mathrm{Tp}\rangle$. While our results concerning plane partitions fixed by elements of~$\langle \mathrm{Row}, \mathrm{Tp}\rangle$ do not directly imply anything about rowmotion for these triangular posets, they do lend credence to the idea that there are nice sieving phenomenon formulas counting plane partitions fixed by many subgroups of $\langle \mathrm{Row}, \mathrm{Tp}\rangle$.

\section{Promotion and transposition} \label{sec:pro_tr}

In this section we prove Theorem~\ref{thm:tr_csiev}. First we recall a combinatorial interpretation of the quantity $\prod_{1\leq i \leq j \leq n-1} \frac{i+j+2M}{i+j}$ which appeared in the aforementioned paper of Proctor~\cite{proctor1990new}. Consider a triangular array $\pi = (\pi_{i,j})_{1 \leq i < j \leq n}$ of nonnegative integers~$\pi_{i,j}\in \mathbb{N}$ like so (we depict the case~$n=4$):
\[\begin{matrix}
 & \pi_{1,2} & \pi_{1,3} & \pi_{1,4} \\
 & & \pi_{2,3} & \pi_{2,4} \\
 & & & \pi_{3,4}
\end{matrix}\]
Suppose further that $\pi$ satisfies the following conditions:
\begin{itemize}\itemsep=0pt
\item $\pi$ is weakly decreasing in rows and columns (i.e., $\pi_{i,j} \geq \pi_{i+1,j}$, $\pi_{i,j} \geq \pi_{i,j+1}$ for all $i$, $j$),
\item the largest entry is less than or equal to~$m$ (i.e., $\pi_{1,2} \leq m$),
\item $m=2M$ is even, and moreover every entry along the first diagonal is even (i.e.,~$\pi_{i,i+1}$ is even for all $i$).
\end{itemize}
Denote the set of such triangular arrays $\pi$ by $\mathrm{CY}(n,2M)$. Proctor~\cite[Theorem~1, formula~`CYI']{proctor1990new} showed that $\#\mathrm{CY}(n,2M) = \prod_{1\leq i \leq j \leq n-1} \frac{i+j+2M}{i+j}$ (and see also the addendum of that paper for other references for this formula).

With this result in hand we are now ready to prove Theorem~\ref{thm:tr_csiev}.

\begin{proof}[Proof of Theorem~\ref{thm:tr_csiev}]
Since $\langle \mathrm{Pro}, \mathrm{Tp}\rangle$ is a dihedral group, we have that $\mathrm{Tp}\cdot \mathrm{Pro}^k$ is conjugate to $\mathrm{Tp}\cdot \mathrm{Pro}^j$ as long as $k$ and $j$ have the same parity, so there are only two cases of the theorem we need to address: $k=0$ and $k=1$. The case $k=0$ is clear from the definition of $ \mathrm{SymMac}(n,m;q)$. So now let us count plane partitions in $\mathrm{PP}^{m}(n\times n)$ fixed by $\mathrm{Tp}\cdot \mathrm{Pro}$.

As mentioned in Section~\ref{sec:intro}, we have
\[\mathrm{SymMac}(n,m;q \mapsto -1) = \begin{cases} \displaystyle \prod_{1\leq i \leq j \leq n-1} \frac{i+j+2M}{i+j} & \textrm{if $m=2M$ is even},\\ 0 &\textrm{otherwise}. \end{cases}\]
We want to show that this quantity is also $\#\big\{\pi \in \mathrm{PP}^{2M}(n\times n)\colon \mathrm{Tp}\cdot \mathrm{Pro}(\pi)=\pi\big\}$.

Recall the diagonal toggles $\mathcal{F}_{k} \coloneqq \prod_{\substack{1\leq i,j \leq n, \\ j-i=k}} \tau_{i,j}$ for~$-n+1 \leq k \leq n-1$ which make up promotion. Let us also use the notation $\mathcal{F}_{-} \coloneqq \mathcal{F}_{-1} \cdots \mathcal{F}_{-n+2}\cdot \mathcal{F}_{-n+1}$ and $\mathcal{F}_{+} \coloneqq \mathcal{F}_{n-1} \cdots \mathcal{F}_{1}$, so that $\mathrm{Pro} = \mathcal{F}_{+} \cdot \mathcal{F}_0 \cdot \mathcal{F}_{-}$. We claim that there is a bijection
\begin{align*}
\big\{\pi \in \mathrm{PP}^{m}(n\times n)\colon \mathrm{Tp}\cdot \mathrm{Pro}(\pi)=\pi\big\} &\xrightarrow{\sim} \big\{\pi \in \mathrm{PP}^{m}(n\times n)\colon \mathrm{Tp}(\pi)=\pi, \mathcal{F}_0(\pi) = \pi\big\}, \\
\pi &\mapsto \mathcal{F}_{-} (\pi).
\end{align*}
An example of this map in the case $n=3$, $m=4$ is
\[ \parbox{0.8in}{\begin{ytableau} 4 & 4 & 3 \\ 3 & 3 & 2 \\ 2 & 2 & 1 \end{ytableau}} \mapsto \parbox{0.8in}{\begin{ytableau} 4 & 4 & 3 \\ 4 & 3 & 2 \\ 3 & 2 & 1 \end{ytableau}} \]

The key point to showing this map is a bijection is that for all $\pi \in \mathrm{PP}^{m}(n\times n)$, and all $-n+1 \leq k \leq n-1$, $\mathrm{Tp}(\mathcal{F}_i(\pi)) = \mathcal{F}_{-i}(\mathrm{Tp}(\pi))$; so that in particular $\mathrm{Tp}(\mathcal{F}_{-}(\pi))=\mathcal{F}^{-1}_+(\mathrm{Tp}(\pi))$.

First let us show that for $\pi \in \mathrm{PP}^{m}(n\times n)$ with $\mathrm{Pro}(\pi)=\mathrm{Tp}(\pi)$, we have $\mathrm{Tp}(\mathcal{F}_{-}(\pi))=\mathcal{F}_{-}(\pi)$ and $\mathcal{F}_0(\mathcal{F}_{-}(\pi))=\mathcal{F}_{-}(\pi)$. Indeed, if $\mathcal{F}_0$ acts nontrivially on $\mathcal{F}_{-}(\pi)$ then $\mathrm{Pro}=\mathcal{F}_{+}\cdot\mathcal{F}_0\cdot\mathcal{F}_{-}$ will alter the main diagonal of~$\pi$ and so we cannot have $\mathrm{Pro}(\pi)=\mathrm{Tp}(\pi)$. And if $\mathcal{F}_0(\mathcal{F}_{-}(\pi))=\mathcal{F}_{-}(\pi)$ and $\mathrm{Pro}(\pi)=\mathrm{Tp}(\pi)$, then $\mathcal{F}_{+}(\mathcal{F}_{-}(\pi))=\mathcal{F}_{+}(\mathcal{F}_0(\mathcal{F}_{-}(\pi))=\mathrm{Pro}(\pi)=\mathrm{Tp}(\pi)$, which means that $\mathcal{F}_{-}(\pi)=\mathcal{F}^{-1}_{+}(\mathrm{Tp}(\pi))=\mathrm{Tp}(\mathcal{F}_{-}(\pi))$.

Next let's show that if $\pi \in \mathrm{PP}^{m}(n\times n)$ satisfies $\mathrm{Tp}(\mathcal{F}_{-}(\pi))=\mathcal{F}_{-}(\pi)$ and $\mathcal{F}_0(\mathcal{F}_{-}(\pi))=\mathcal{F}_{-}(\pi)$, then $\mathrm{Pro}(\pi)=\mathrm{Tp}(\pi)$. The argument is basically the same. The assumptions imply that
\begin{align*} \mathrm{Pro}(\pi)& =\mathcal{F}_{+}(\mathcal{F}_{0}(\mathcal{F}_{-}(\pi)))=\mathcal{F}_{+}(\mathcal{F}_{-}(\pi))
=\mathcal{F}_{+}(\mathrm{Tp}(\mathcal{F}_{-}(\pi)))\\
& =\mathcal{F}_{+}\big(\mathcal{F}^{-1}_{+}(\mathrm{Tp}(\pi))\big)=\mathrm{Tp}(\pi),
\end{align*}
as required.
	
So indeed the map $\pi \mapsto \mathcal{F}_{-}(\pi)$ is a bijection between the claimed sets. Thus we reduced the problem to counting $\#\{\pi \in \mathrm{PP}^{m}(n\times n)\colon \mathrm{Tp}(\pi)=\pi, \mathcal{F}_0(\pi) = \pi\}$.

Now for $\pi \in \mathrm{PP}^{m}(n\times n)$ with $\mathrm{Tp}(\pi)=\pi$, we have that $\mathcal{F}_0(\pi) = \pi$ if and only if we have~$\pi_{i,i} = \frac{1}{2}(\pi_{i,i-1} + \pi_{i,i+1})$ for all $1 \leq i \leq n$ (where we recall the conventions that~$\pi_{0,j}\coloneqq m$ and~$\pi_{i,n+1}\coloneqq 0$). Because $\pi$ has integer entries, this can only happen if $\pi_{i,i-1}$ and $\pi_{i,i+1}$ have the same parity for all $1 \leq i \leq n$; since $\pi_{n,n+1}=0$, this means we need all the $\pi_{i,i+1}$ to be even, including $\pi_{0,1}=m$. So if $m$ is odd then $\#\{\pi \in \mathrm{PP}^{m}(n\times n)\colon \mathrm{Tp}(\pi)=\pi, \mathcal{F}_0(\pi) = \pi\} = 0$. And moreover, it follows from what we have just explained that if $m=2M$ is even then there is a bijection from $\{\pi \in \mathrm{PP}^{m}(n\times n)\colon \mathrm{Tp}(\pi)=\pi, \mathcal{F}_0(\pi) = \pi\}$ to $\mathrm{CY}(n,2M)$ which takes the square array $(\pi_{i,j})_{1 \leq i,j \leq n}$ to the triangular array $(\pi_{i,j})_{1 \leq i < j \leq n}$ of entries in the upper right corner. Together with Proctor's formula for $\#\mathrm{CY}(n,2M)$, this completes the proof.
\end{proof}

\section{Promotion and transpose-complementation}

In this section we prove Theorem~\ref{thm:trcom_csiev}. We do this by extending Rhoades's~\cite{rhoades2010cyclic} approach to cyclic sieving for tableaux using the dual canonical basis of ${\rm GL}(a+b)$ representations. But actually, rather than hew closely to Rhoades's presentation, we instead follow the presentation of Lam~\cite{lam2019cyclic}. Lam explained how the relevant ${\rm GL}(a+b)$ representation is the coordinate ring of the Grassmannian $\mathrm{Gr}(a,a+b)$. We find this geometric perspective useful. Also, as hinted at in Section~\ref{sec:intro}, we owe a great debt to the papers of Stembridge~\cite{stembridge1994minuscule, stembridge1996canonical} and especially Kuperberg~\cite{kuperberg1994minuscule} for explaining how involutive symmetries of plane partitions can be realized as algebra automorphisms on these coordinate rings.

In this section we work with semistandard tableaux rather than plane partitions. We recall that the correspondence between plane partitions and tableaux of rectangular shape is explained in Appendix~\ref{sec:appendix}. For a partition $\lambda$ we use $\mathrm{SSYT}(\lambda,k)$ to denote the set of semistandard Young tableaux of shape $\lambda$ with entries less than or equal to $k$. We will freely use the bijection $\Psi\colon \mathrm{PP}^{m}(a\times b) \xrightarrow{\sim} \mathrm{SSYT}(m^a,a+b)$ defined in the appendix. Via this bijection promotion $\mathrm{Pro}$ and complementation $\mathrm{Co}$ are viewed as operators on $\mathrm{SSYT}(m^a,a+b)$ and transposition $\mathrm{Tp}$ is viewed as an operator on $\mathrm{SSYT}(m^n,2n)$. The behavior of these operators on tableaux is explained in Proposition~\ref{prop:tableau_ops}.

\subsection{Background on Grassmannian coordinate rings}

Before we can prove Theorem~\ref{thm:trcom_csiev} we have to review a bit about Grassmannians and the representation theory arising from their study. We start with the Grassmannian.

The \emph{Grassmannian} $\mathrm{Gr}(a,a+b)$ is the space of $a$-dimensional subspaces of the complex vector space $V=\mathbb{C}^{a+b}$. There is a very well-known system of coordinates on the Grassmannian called the \emph{Pl\"{u}cker coordinates}. Let $U \in \mathrm{Gr}(a,a+b)$ and choose an ordered basis $v_1,\dots,v_a \in \mathbb{C}^{a+b}$ of~$U$; let $I=\{i_1, i_2, \dots, i_a\}\subseteq \{1,2,\dots,a+b\}$ be a subset of size~$a$; then the \emph{Pl\"{u}cker coordinate} $\Delta_{I}(U)$ is equal to the maximal minor of the $(a+b) \times a$ matrix with column vectors $v_1,\dots,v_a$ given by selecting rows~$i_1,\dots,i_a$. The Grassmannian is a projective variety and the map
\[U \mapsto [\Delta_I(U)\colon I\subseteq \{1,\dots,a+b\} \textrm{ of size $a$}]\]
is an embedding of $\mathrm{Gr}(a,a+b)$ into $\mathbb{P}^{\binom{a+b}{a}-1}$ known as the \emph{Pl\"{u}cker embedding}. We use the notation $\widehat{\mathrm{Gr}}(a,a+b)\subseteq \mathbb{C}^{\binom{a+b}{a}}$ to denote the \emph{affine cone over $\mathrm{Gr}(a,a+b)$} in its Pl\"{u}cker embedding. And we use $R(a,a+b)$ to denote the coordinate ring of $\widehat{\mathrm{Gr}}(a,a+b)$. In other words, $R(a,a+b)$ is the commutative ring
\[ R(a,a+b) = \mathbb{C}[\Delta_I\colon I\subseteq \{1,\dots,a+b\} \textrm{ of size $a$}]/ \langle \textrm{Pl\"{u}cker relations}\rangle,\]
where the \emph{Pl\"{u}cker relations} are the well-known relations cutting out $\mathrm{Gr}(a,a+b)$ as a subset of $\mathbb{P}^{\binom{a+b}{a}-1}$. Equivalently we may think of $R(a,a+b)$ as the homogeneous coordinate ring of $\mathrm{Gr}(a,a+b)$.

See for instance~\cite[Chapter 1]{seshadri2014standard} for the basics concerning the coordinate ring of the Grassmannian. We use $R(a,a+b)_m$ to denote the functions in $R(a,a+b)$ of homogeneous degree~$m$. For $T \in \mathrm{SSYT}(m^a,a+b)$ we set $M(T) \coloneqq \Delta_{I_1} \Delta_{I_2} \cdots \Delta_{I_m}$, where~$I_1,I_2,\dots,I_m$ are the columns of $T$. Note $M(T)$ belongs to $R(a,a+b)_m$. The set~$\{M(T)\colon T \in \mathrm{SSYT}(m^a,a+b)\}$ is the \emph{standard monomial basis} of $R(a,a+b)_m$. It is a classical result, going back to Young, that $\{M(T)\colon T \in \mathrm{SSYT}(m^a,a+b)\}$ is indeed a linear basis of $R(a,a+b)_m$.

Now we review representations of the general linear group and canonical bases.

We will find the following notation for matrices useful: $\mathrm{diag}(x_1,x_2,\dots,x_k)$ is the diagonal $k\times k$ matrix with diagonal entries~$x_1$ to $x_k$ from upper-left to lower-right; $\mathrm{antidiag}(x_1,x_2,\dots,x_k)$ is the anti-diagonal $k\times k$ matrix with anti-diagonal entries~$x_1$ to $x_k$ from upper-right to lower-left; and of course $\mathrm{Id}_k \coloneqq \mathrm{diag}(\overbrace{1,1,\dots,1}^{k})$ is the $k\times k$ identity matrix. We also often write matrices in block form.

The \emph{general linear group} ${\rm GL}(a+b)$ is the group of invertible linear transformations acting on $V=\mathbb{C}^{a+b}$. We usually think of the elements of the general linear group as~$(a+b)\times (a+b)$ $\mathbb{C}$-matrices with nonzero determinant, having implicitly chosen an ordered basis $e_1,e_2,\dots,e_{a+b}$ of~$V$. The \emph{special linear group} ${\rm SL}(a+b) \subseteq {\rm GL}(a+b)$ is the subgroup of matrices in ${\rm GL}(a+b)$ of determinant one. The Lie algebra corresponding to the Lie group ${\rm SL}(a+b)$ is the Lie algebra~$\mathfrak{sl}_{a+b}$ of $(a+b)\times (a+b)$ $\mathbb{C}$-matrices with trace zero, with Lie bracket given by the commutator. The Lie algebra $\mathfrak{sl}_{a+b}$ is simple.

Inside of ${\rm GL}(a+b)$ is the algebraic torus~$T$ of diagonal matrices $\mathrm{diag}(x_1,\dots,x_{a+b})$ with $\prod_{i=1}^{a+b}x_i\neq 0$. The torus of~${\rm SL}(a+b)$ is the subgroup of those diagonal matrices with~$\prod_{i=1}^{a+b}x_i=1$. The symmetric group $\mathfrak{S}_{a+b}$ on $a+b$ letters is the quotient of the normalizer in ${\rm GL}(a+b)$ of $T$ by~$T$. In other words, $\mathfrak{S}_{a+b}$ is the \emph{Weyl group} of~${\rm GL}(a+b)$. The symmetric group is also the Weyl group of ${\rm SL}(a+b)$. Thus elements of the symmetric group can be lifted to the general linear group in various ways; but conjugation by elements of the symmetric group gives a well-defined action on the torus of~${\rm GL}(a+b)$, and on the torus of ${\rm SL}(a+b)$.

See for instance~\cite[Chapter 15]{fulton1991representation} for the basics concerning the representation theory of \mbox{${\rm GL}(a+b)$}. If $V$ is a ${\rm GL}(a+b)$ representation, then a vector $v \in V$ is said to be a \emph{weight vector} with weight $\alpha=(\alpha_1,\dots,\alpha_{a+b})$ if every torus element $\mathrm{diag}(x_1,x_2,\dots,x_{a+b}) \in {\rm GL}(a+b)$ sends the vector $v$ to $x_1^{\alpha_1}\cdots x_{a+b}^{\alpha_{a+b}} \cdot v$. For any partition $\lambda = (\lambda_1 \geq \lambda_2 \geq \cdots \geq \lambda_{a+b} \geq 0)$ with at most $a+b$ parts we have an irreducible, finite-dimensional representation $V(\lambda)$ of ${\rm GL}(a+b)$ of highest weight~$\lambda$. Each of these ${\rm GL}(a+b)$ representations $V(\lambda)$ is also an irreducible ${\rm SL}(a+b)$ representation, and hence an irreducible $\mathfrak{sl}_{a+b}$ representation. For $0 \leq i \leq a+b$ we define the partition $\omega_{i} \coloneqq (\overbrace{1,1,\dots,1}^{i},\overbrace{0,\dots,0}^{a+b-i})$. The $V(\omega_i)$ for $0 < i < a+b$ are the \emph{fundamental representations} of $\mathfrak{sl}_{a+b}$. Because ${\rm SL}(a+b)$ is connected and simply connected, representations of ${\rm SL}(a+b)$ and $\mathfrak{sl}_{a+b}$ are in exact correspondence.

Lusztig~\cite{lusztig1990canonical} and Kashiwara~\cite{ kashiwara1993global} constructed a \emph{canonical} (or \emph{global}) \emph{basis} of the irreducible $U_q(\mathfrak{sl}_{a+b})$-module~$V_q(\lambda)$ (and their two constructions are known to give the same basis~\cite{grojnowski1993comparison}). Here $U_q(\mathfrak{sl}_{a+b})$ is the quantized universal enveloping algebra of~$\mathfrak{sl}_{a+b}$, a deformation of the universal enveloping algebra $U(\mathfrak{sl}_{a+b})$. By setting $q \mapsto 1$ in their work, and by fixing a particular highest weight vector $v_+$ of~$V(\lambda)$, we have a \emph{canonical basis} $\{G(T)\colon T \in \mathrm{SSYT}(\lambda,a+b)\}$ of~$V(\lambda)$. For a~tableau $T \in \mathrm{SSYT}(\lambda,a+b)$, the \emph{weight} of $T$ is $\mathrm{wt}(T) \coloneqq (\alpha_1,\dots,\alpha_{a+b})$ where $\alpha_i$ is the number of $i$'s in $T$. The canonical basis vector $G(T)$ is a weight vector of~$V(\lambda)$ with weight $\mathrm{wt}(T)$. We use $\{H(T)\colon T \in \mathrm{SSYT}(\lambda,a+b)\}$ to denote the dual basis to~$\{G(T)\colon T \in \mathrm{SSYT}(\lambda,a+b)\}$; to be clear, this is a basis of the dual space $V(\lambda)^{*}$. We refer to $\{H(T)\colon T \in \mathrm{SSYT}(\lambda,a+b)\}$ as the \emph{dual canonical basis} of $V(\lambda)^{*}$.

We can view the Grassmannian as a quotient $\mathrm{Gr}(a,a+b) = \mathrm{Mat}^{\bullet}(a+b,a) / {\rm GL}(a)$, where $\mathrm{Mat}^{\bullet}(a+b,a)$ is the space of $(a+b)\times a$ $\mathbb{C}$-matrices of rank $a$, and ${\rm GL}(a)$ acts on $\mathrm{Mat}^{\bullet}(a+b,a)$ on the right in the obvious way. This space of matrices $\mathrm{Mat}^{\bullet}(a+b,a)$ caries an obvious left action of ${\rm GL}(a+b)$ which commutes with the right ${\rm GL}(a)$ action, and in this way we obtain an action of ${\rm GL}(a+b)$ on $\mathrm{Gr}(a,a+b)$. Similarly, we can view $\widehat{\mathrm{Gr}}(a,a+b) \setminus \{0\}$ as the quotient $\widehat{\mathrm{Gr}}(a,a+b) \setminus \{0\} = \mathrm{Mat}^{\bullet}(a+b,a) / {\rm SL}(a)$, and in this way we obtain an action of ${\rm GL}(a+b)$ on $\widehat{\mathrm{Gr}}(a,a+b) \setminus \{0\}$ which is compatible with the action of~${\rm GL}(a+b)$ on $\mathrm{Gr}(a,a+b)$. We extend this action to an action of~${\rm GL}(a+b)$ on all of $\widehat{\mathrm{Gr}}(a,a+b)$ by declaring $g \cdot 0 = 0$ for all $g \in {\rm GL}(a+b)$. We then get an action of~${\rm GL}(a+b)$ on $R(a,a+b)$ via algebra automorphisms by inverting and pulling back: for~$g \in {\rm GL}(a+b)$ and $f \in R(a,a+b)$ we set $(g\cdot f)(U) \coloneqq f\big(g^{-1} U\big)$ for all $U \in \widehat{\mathrm{Gr}}(a,a+b)$. It is well known, for instance via the classical Borel--Weil theorem, that as ${\rm GL}(a+b)$ representations we have~$R(a,a+b)_m \simeq V(m \omega_a)^{*}$, where $V(m \omega_a)^{*}$ is the \emph{dual} of the irreducible representation $V(m \omega_a)$. (Recall that if $\rho\colon G \to {\rm GL}(V)$ is a representation of a~group~$G$, then the \emph{dual representation} $\rho^{*}\colon G \to {\rm GL}(V^*)$ is the representation where $g\in G$ acts on the dual space~$V^*$ by~$\rho(g^{-1})^T$, with the $T$ superscript denoting transpose.) This means we can consider the dual canonical basis $\{H(T)\colon T \in \mathrm{SSYT}(m^a,a+b)\}$ as a basis of $R(a,a+b)_m$.

\subsection{Grassmannian coordinate ring automophisms}

We now define several algebra automorphisms of $R(a,a+b)$. These automorphisms are at the heart of our proof of Theorem~\ref{thm:trcom_csiev}: they will correspond to promotion, complementation, and transposition of plane partitions.

Some of these automorphisms are (the actions of) elements of ${\rm GL}(a+b)$. First we define the \emph{twisted cyclic shift} $\chi \in {\rm GL}(a+b)$:
\[\chi \coloneqq \begin{pmatrix}\def\arraystretch{1.4}\begin{array}{c|c} 0 & (-1)^{a-1} \\ \hline \parbox{2cm}{\begin{center}\vspace{0.5cm} $\mathrm{Id}_{a+b-1}$ \vspace{0.5cm} \end{center}} & 0 \end{array}\end{pmatrix} \in {\rm GL}(a+b).\]
The matrix $\chi$ is a lift of the \emph{long cycle} $c \coloneqq \left(\begin{smallmatrix} 1 & 2 & \dots & a+b-1 & a+b \\ 2 & 3 & \dots & a+b & 1 \end{smallmatrix}\right) \in \mathfrak{S}_{a+b}$. Note that $\chi$ has order $a+b$ acting on $\widehat{\mathrm{Gr}}(a,a+b)$ because $\chi^{a+b}$ multiplies each vector in $\mathbb{C}^{a+b}$ by $(-1)^{(a-1)}$ and hence each Pl\"{u}cker coordinate by $(-1)^{a(a-1)}=1$.

We next define the \emph{twisted reflection} $\overline{w_0} \in {\rm GL}(a+b)$:
\[\overline{w_0} \coloneqq i^{(a-1)} \cdot \mathrm{antidiag}(1,1,\dots,1) \in {\rm GL}(a+b).\]
We denote this element of ${\rm GL}(a+b)$ by $\overline{w_0}$ because it is a particular lift of the \emph{longest element} $w_0 \coloneqq \left(\begin{smallmatrix} 1 & 2 & \dots & a+b \\ a+b & a+b-1 & \dots & 1 \end{smallmatrix}\right) \in \mathfrak{S}_{a+b}$ in the symmetric group. Note that $\overline{w_0}$ is an involution acting on $\widehat{\mathrm{Gr}}(a,a+b)$ because $w_0^2$ multiplies each vector in~$\mathbb{C}^{a+b}$ by $(-1)^{(a-1)}$ and hence each Pl\"{u}cker coordinate by $(-1)^{a(a-1)}=1$.

For the next several paragraphs we suppose that $a=b=n$. Let $B$ be the following skew-symmetric $2n\times 2n$ matrix:
\[ B \coloneqq \mathrm{antidiag}(1,-1,1,-1,\dots,-1) \in {\rm GL}(2n).\]
Then $B$ defines a symplectic form $\langle\cdot , \cdot \rangle_B$ on $V=\mathbb{C}^{2n}$ by $\langle x , y \rangle_B \coloneqq x^T \cdot B \cdot y$ for all~$x,y\in\mathbb{C}^{2n}$. This symplectic form defines an outer automorphism $\phi_{B}$ of ${\rm GL}(2n)$ where a linear transformation is sent by $\phi_{B}$ to the inverse of its transpose with respect to the identification of $V$ and $V^{*}$ induced by $\langle\cdot , \cdot \rangle_B$. At the level of matrices we have
\begin{align*}
\phi_{B}\colon \ {\rm GL}(2n) &\to {\rm GL}(2n), \\
A &\mapsto B^{-1} \cdot \big(A^{T}\big)^{-1} \cdot B,
\end{align*}
where the superscript $T$ denotes usual matrix transposition. Clearly $\phi_B$ is an involution. And~$\phi_B$ restricts to an involutive outer automorphism $\phi_B\colon {\rm SL}(2n) \to {\rm SL}(2n)$ of the special linear group.

The \emph{symplectic group} $\mathrm{Sp}(2n)\subseteq {\rm GL}(2n)$ is the subgroup of ${\rm GL}(2n)$ consisting of those matrices $A \in {\rm GL}(2n)$ with $\phi_{B}(A)=A$. In fact $\mathrm{Sp}(2n)\subseteq \mathrm{SL}(2n)$ (this can be seen by consideration of Pfaffians). Inside of $\mathrm{Sp}(2n)$ we have the algebraic torus consisting of those diagonal matrices $D=\mathrm{diag}(x_1,x_2,\dots,x_{2n})$ with $x_i=x_{2n+1-i}^{-1}$ for all $i$. The analog of the symmetric group here is the \emph{hyperoctahedral group}. The hyperoctahedral group $(\mathbb{Z}/2\mathbb{Z}) \wr \mathfrak{S}_n$ is the subset of the symmetric group $\mathfrak{S}_{2n}$ consisting of those permutations $\sigma \in \mathfrak{S}_{2n}$ for which~$\sigma(i) = (2n+1)-\sigma((2n+1)-i)$ for all $1 \leq i \leq 2n$. The hyperoctahedral group is the Weyl group of $\mathrm{Sp}(2n)$. Thus, elements of~$(\mathbb{Z}/2\mathbb{Z}) \wr \mathfrak{S}_n$ act on the torus of $\mathrm{Sp}(2n)$ by conjugation.

The symplectic form $\langle\cdot , \cdot \rangle_B$ also gives rise to the following automorphism (of projective algebraic varieties) on the middle-dimensional Grassmannian $\mathrm{Gr}(n,2n)$:
\begin{align*}
\widetilde{\phi}_{B}\colon \ \mathrm{Gr}(n,2n) &\to \mathrm{Gr}(n,2n), \\
U &\mapsto U^{\perp},
\end{align*}
where $U^{\perp}$ is the orthogonal complement of the subspace $U$ with respect to~$\langle\cdot , \cdot \rangle_B$, that is, $U^{\perp} \coloneqq \big\{x\in \mathbb{C}^{2n}\colon \langle x,y\rangle_B = 0 \textrm{ for all $y\in U$}\big\}$. Clearly~$\widetilde{\phi}_B$ is an involution. We extend~$\widetilde{\phi}_B$ to an automorphism $\widetilde{\phi}_B\colon \widehat{\mathrm{Gr}}(n,2n) \to \widehat{\mathrm{Gr}}(n,2n)$ (of affine algebraic varieties) in a unique way by requiring that $\Delta_{\{1,2,\dots,n\}}(U) = \Delta_{\{1,2,\dots,n\}}\big(\widetilde{\phi}_B(U)\big)$ for all~$U \in \widehat{\mathrm{Gr}}(n,2n)$. By abuse of notation we also use $\widetilde{\phi}_{B}$ to denote the induced algebra automorphism on the coordinate ring~$R(n,2n)$ of $\widehat{\mathrm{Gr}}(n,2n)$ given by pulling back $\widetilde{\phi}_{B}$. All of these $\widetilde{\phi}_{B}$ remain involutions.

The essential property connecting $\widetilde{\phi}_{B}$ and $\phi_{B}$ is that for any $A \in {\rm GL}(2n)$, we have for all $U \in \mathrm{Gr}(n,2n)$ that~$\widetilde{\phi}_B(A \cdot U) = \phi_{B}(A) \cdot \widetilde{\phi}_B(U)$. This is easy to see from the fact that $\langle A x , \phi_{B}(A) y\rangle_B = \langle x, y\rangle_B$ for all $x,y\in \mathbb{C}^{2n}$. Moreover, the amount that the action of a matrix $A \in {\rm GL}(2n)$ scales $\Delta_{\{1,2,\dots,n\}}(U)$ for $U \in \widehat{\mathrm{Gr}}(n,2n)$ is given by the principal $n\times n$ minor of $A$; a simple computation with block matrices shows that for $A \in {\rm SL}(2n)$ this principal minor is the same for $A$ and for $\phi_{B}(A)$. Hence for any $A\in {\rm SL}(2n)$, we in fact have that $\widetilde{\phi}_B(A \cdot U) = \phi_{B}(A) \cdot \widetilde{\phi}_B(U)$ for all~$U \in \widehat{\mathrm{Gr}}(n,2n)$. That is to say, for any $A\in {\rm SL}(2n)$ we have the following equality of automorphisms of the coordinate ring~$R(n,2n)$: $\widetilde{\phi}_B \cdot A = \phi_B(A) \cdot \widetilde{\phi}_B$.

\subsection{Behavior of coordinate ring automorphisms on bases}

Now we study how these automorphisms behave on the various bases of $R(a,a+b)$.

As in Example~\ref{ex:m1}, for a permutation $\sigma \in \mathfrak{S}_{a+b}$ and a subset $I\subseteq \{1,\dots, a+b\}$ we write~$\sigma(I) \coloneqq \{\sigma(i)\colon i \in I\}$; and we also write $- I \coloneqq \{1,\dots,a+b\}\setminus I$.

\begin{Lemma} \label{lem:plucker_actions}
The actions of the above automorphisms on the Pl\"{u}cker coordinates generating the coordinate ring $R(a,a+b)$ are:
\begin{itemize}\itemsep=0pt
\item $\chi(\Delta_I) = \Delta_{c(I)}$,
\item $\overline{w_0}(\Delta_I) = \Delta_{w_0(I)}$,
\item (if $a=b=n$) $\widetilde{\phi}_{B}(\Delta_I) = \Delta_{-w_0(I)}$.
\end{itemize}
\end{Lemma}
\begin{proof}The first two bulleted items are stating simple facts about how matrix minors behave when rows are permuted. The factor $i^{(a-1)}$ in the definition of $\overline{w_0}$ is there because $i^{(a-1)} \, \mathrm{Id}_{a+b}$ multiplies each Pl\"{u}cker coordinate by $i^{-a(a-1)} = (-1)^{a(a-1)}$, which is exactly the right number of minus signs to cancel the number of row transpositions we need to vertically flip an $a\times a$ submatrix. Similarly, the entry of~$(-1)^{a-1}$ in the definition of $\chi$ cancels out the row transpositions needed to bring the last row of an $a\times a$ submatrix to the front.

The statement about $\widetilde{\phi}_B$ is explained, in the somewhat different but equivalent language of alternating forms and the Hodge star, in~\cite[proof of Theorem 4.1]{kuperberg1994minuscule}. It is also not hard to see directly. Let $U \in \mathrm{Gr}(a,a+b)$, and suppose $U$ lies in the dense open subset of the Grassmannian where $\Delta_{\{1,2,\dots,n\}}(U)\neq 0$. Let's represent~$U$ by a matrix in reduced column echelon form whose column span is $U$ (and note that $\Delta_{\{1,2,\dots,n\}}(U)\neq 0$ implies the upper $n\times n$ square submatrix of this matrix is~$\mathrm{Id}_n$). Then the effect of~$\widetilde{\phi}_B$ is to ``transpose'' the lower $n\times n$ square submatrix of this matrix across its main anti-diagonal, while also multiplying the entries in this square submatrix by $\pm 1$ in a checkerboard pattern, as the following diagrams depict in the cases $n=2,3,4$:
\begin{gather*} \begin{pmatrix} 1 & 0 \\ 0 & 1 \\ a & c \\ b & d \end{pmatrix} \xrightarrow{\widetilde{\phi}_B} \begin{pmatrix} 1 & 0 \\ 0 & 1 \\ -d & c \\ b & -a \end{pmatrix}, \qquad \begin{pmatrix} 1 & 0 & 0 \\ 0 & 1 & 0 \\ 0 & 0 & 1 \\ a & d & g \\ b & e & h \\ c & f & i \end{pmatrix} \xrightarrow{\widetilde{\phi}_B} \begin{pmatrix} 1 & 0 & 0 \\ 0 & 1 & 0 \\ 0 & 0 & 1 \\ i & -h & g \\ -f & e & -d \\ c & -b & a \end{pmatrix}, \\
 \begin{pmatrix} 1 & 0 & 0 & 0 \\ 0 & 1 & 0 & 0 \\ 0 & 0 & 1 & 0 \\ 0 & 0 & 0 & 1 \\ a & e & i & m \\ b & f & j & n \\ c & g & k & o \\ d & h & l & p \end{pmatrix} \xrightarrow{\widetilde{\phi}_B} \begin{pmatrix} 1 & 0 & 0 & 0 \\ 0 & 1 & 0 & 0 \\ 0 & 0 & 1 & 0 \\ 0 & 0 & 0 & 1 \\ -p & o & -n & m \\ l & -k & j & -i \\ -h & g & -f & e \\ d & -c & b & -a \end{pmatrix}.
\end{gather*}
This matrix representation makes it easy to check that $\Delta_I\big(\widetilde{\phi}_B(U)\big) = \Delta_{-w_0(I)}(U)$. This is because the maximal minors of the $2n\times n$ matrix correspond, up to sign, to \emph{all} of the minors of its lower $n\times n$ square submatrix (see~\cite[Lemma~3.9]{postnikov2006total}); and the checkerboard pattern of signs exactly addresses the sign issue. Then we observe that $\Delta_I\big(\widetilde{\phi}_B(U)\big) = \Delta_{-w_0(I)}(U)$ in fact holds for all $U\in\mathrm{Gr}(a,a+b)$ since it holds on a dense open subset. It also holds for $U \in \widehat{\mathrm{Gr}}(a,a+b)$ since we declared that $\Delta_{\{1,2,\dots,n\}}(U) = \Delta_{\{1,2,\dots,n\}}\big(\widetilde{\phi}_B(U)\big)$ and $-w_0(\{1,2,\dots,n\})=\{1,2,\dots,n\}$.
\end{proof}

\begin{Corollary} \label{cor:std_monomial_actions}
The actions of the automorphisms $\overline{w_0}$ and $\widetilde{\phi}_{B}$ on the standard monomial basis of the ring~$R(a,a+b)_m$ are:
\begin{itemize}\itemsep=0pt
\item $\overline{w_0}(M(T))=M(\mathrm{Co}(T))$,
\item $($if $a=b=n)$ $\widetilde{\phi}_{B}(M(T)) = M(\mathrm{Tp}(T))$.
\end{itemize}
\end{Corollary}
\begin{proof}
These follow immediately from Lemma~\ref{lem:plucker_actions} if we recall the effects of $\mathrm{Co}$ and $\mathrm{Tp}$ on tableaux as described in Proposition~\ref{prop:tableau_ops}. For a tableau $T \in \mathrm{SSYT}(m^a,a+b)$, the columns of the complementary tableau $\mathrm{Co}(T)$ are $w_0(I_m), w_0(I_{m-1}),\dots, w_0(I_1)$, where $I_1,I_2,\dots,I_m$ are the columns of $T$. Similarly, for a tableau $T \in \mathrm{SSYT}(m^n,2n)$, the columns of the transposed tableau $\mathrm{Tp}(T)$ are $-w_0(I_1), -w_0(I_{2}),\dots, -w_0(I_m)$, where $I_1,I_2,\dots,I_m$ are the columns of $T$.
\end{proof}

Essentially via Corollary~\ref{cor:std_monomial_actions}, Stembridge~\cite{stembridge1994minuscule} and Kuperberg~\cite{kuperberg1994minuscule} were able to deduce the $q=-1$ results discussed in the Section~\ref{sec:intro}: Theorems~\ref{thm:com_qneg} and~\ref{thm:trcom_qneg}. Note crucially, however, that we {\it do not} have $\chi(M(T)) = M(\mathrm{Pro}(T))$. Indeed, the whole point of using sophisticated bases like the dual canonical basis is that the naive bases like the standard monomial basis fail to behave well under the action of the long cycle.

This brings us to the main algebraic result we need to prove Theorem~\ref{thm:trcom_csiev}:

\begin{Theorem} \label{thm:dual_canonical_actions}
The actions of the above automorphisms on the dual canonical basis of the ring~$R(a,a+b)_m$ are:
\begin{itemize}\itemsep=0pt
\item $\chi(H(T)) = H(\mathrm{Pro}(T))$,
\item $\overline{w_0}(H(T))=H(\mathrm{Co}(T))$,
\item (if $a=b=n$) $\widetilde{\phi}_{B}(H(T)) = H(\mathrm{Tp}(T))$.
\end{itemize}
\end{Theorem}
\begin{proof}
The first bulleted item is~\cite[Theorem 2]{lam2019cyclic}. But in fact this result is essentially due to Rhoades~\cite{rhoades2010cyclic}. To obtain this result Rhoades used a particular realization of the canonical basis in terms of \emph{Kazhdan--Lusztig immanants}~\cite{rhoades2006kl1,rhoades2010kl2} which was first introduced by Du~\cite{du1992canonical} and further developed by Skandera~\cite{skandera2008canonical}. (In~\cite{du1995canonical} Du showed that his canonical basis is the same as that of Lusztig~\cite{lusztig1990canonical}.) Very recently, Rush~\cite{rush2020restriction} gave a new proof of Rhoades's result which is more conceptual/abstract.

For the second bulleted item: Berenstein--Zelevinsky~\cite[Proposition 8.8]{berenstein1996canonical}, and Stembridge~\cite[Theorem 1.2]{stembridge1996canonical} showed that for \emph{any} irreducible ${\rm GL}(a+b)$ representation $V(\lambda)$, multiplication by a~lift of $w_0$ corresponds (up to an overall $\pm1$ sign) to evacuation of tableaux in the dual canonical basis. As explained in the appendix, in the case of rectangular tableaux evacuation is the same as complementation.

For the third bulleted item: Berenstein--Zelevinsky~\cite[Section~7]{berenstein1996canonical} again explained the effect of twisting by the ${\rm GL}(2n)$ automorphism $\phi_B$ on the dual canonical basis of \emph{any} irreducible representation~$V(\lambda)$. (This automorphism is denoted by $\psi$ in~\cite[Section~7]{berenstein1996canonical}.) This effect is described in terms of the ``multisegment duality''~\cite{knight1996representations}. It should be possible to show that this involution defined in terms of the multisegment duality reduces to transposition in the case of rectangular partitions, just like evacuation also radically simplifies in the rectangular case.

However, there are some annoying technicalities we would have to deal with in order to directly apply the work of Berenstein--Zelevinsky~\cite{berenstein1996canonical}. For instance, we would have to show that their indexing of the dual canonical basis is compatible with that of Du~\cite{du1992canonical}. Let us instead explain a~different way to conclude the second and third bulleted items. The idea is still to observe that these involutions correspond to automorphisms of the quantized universal enveloping algebra. But we can exploit the fact that we are working in the particularly nice ``rectangular'' setting where the standard monomial basis also behaves well with respect to these involutions (which is not always the case). It turns out that we can piggyback off of the result for the standard monomial basis to obtain the result we want for the dual canonical basis.

In discussing automorphisms of the quantized universal enveloping algebra we follow the presentation of Berenstein--Zelevinsky~\cite[Section~7]{berenstein1996canonical}. If $\phi\colon U_q(\mathfrak{sl}_{a+b}) \to U_q(\mathfrak{sl}_{a+b})$ is an algebra automorphism, then from a $U_q(\mathfrak{sl}_{a+b})$-module $V$ we get another module $\,^{\phi}V$ by twisting by $\phi$: $\,^{\phi}V = V$ as a vector space but we have $g \in U_q(\mathfrak{sl}_{a+b})$ act on $\,^{\phi}V$ by $\phi(g)$. If $V = V_q(\lambda)$ is irreducible then so is $\,^{\phi}V$: say $\,^{\phi}V\simeq V_q(\phi(\lambda))$ for the corresponding highest weight~$\phi(\lambda)$. Thus, abusing notation, we get an isomorphism of vector spaces $\phi\colon V_q(\lambda) \to V_q(\phi(\lambda))$ also denoted $\phi$ which satisfies $\phi(g v) = \phi(g) \phi(v)$ for all $g\in U_q(\mathfrak{sl}_{a+b})$ and $v\in V_q(\lambda)$. Because we have a choice of highest weight vector, the map $\phi\colon V_q(\lambda) \to V_q(\phi(\lambda))$ is uniquely defined only up to an overall scalar.

The quantized universal enveloping algebra $U_q(\mathfrak{sl}_{a+b})$ is generated by elements~$E_i$, $F_i$, $K_i$, $K_i^{-1}$ for $1\leq i <a+b$ subject to certain relations involving the parameter $q$. There are two involutive automorphisms of $U_q(\mathfrak{sl}_{a+b})$ we want to consider. Following~\cite[Section~7]{berenstein1996canonical} we call these automorphisms $\eta$ and $\psi$. They are given as follows
\begin{gather*}
\eta(E_i) = F_{a+b+1-i},\qquad \eta(F_i) = E_{a+b+1-i},\qquad \eta(K_i) = K_{a+b+1-i}^{-1}, \\
\psi(E_i) = E_{a+b+1-i}, \qquad \psi(F_i) = F_{a+b+1-i},\qquad \psi(K_i) = K_{a+b+1-i}.
\end{gather*}
At the level of weights we have $\eta(\lambda) = \lambda$ for all $\lambda$, and $\psi(m\omega_n) = m\omega_n$ for all~$m$ in the case $a=b=n$. Thus, we get involutions $\eta \colon V_q(m\omega_a)\to V_q(m\omega_a)$, and $\psi\colon V_q(m\omega_n)\to V_q(m\omega_n)$ in the case $a=b=n$. General results of Lusztig~\cite[Proposition~21.1.2]{lusztig1993quantum} (see also Berenstein--Zelevinsky~\cite[Proposition~7.1]{berenstein1996canonical}) imply that both $\eta$ and $\psi$ permute the canonical bases of these $U_q(\mathfrak{sl}_{a+b})$-modules; indeed, $\eta$ is the so-called ``Lusztig involution'', while $\psi$ is more-or-less the involutive automorphism induced by the Dynkin diagram symmetry in type~$A_{2n-1}$. These descend to involutions of the ${\rm GL}(a+b)$ representations $\eta\colon V(m\omega_a)\to V(m\omega_a)$ and $\psi\colon V(m\omega_n)\to V(m\omega_n)$ which permute the corresponding canonical bases. In this way we also get involutions $\eta\colon V(m\omega_a)^{*}\to V(m\omega_a)^{*}$ and $\psi\colon V(m\omega_n)^{*}\to V(m\omega_n)^{*}$ on the dual spaces which permute the dual canonical bases. (We are glossing over the fact that the maps $\eta$ and~$\psi$ between irreducible modules are only defined up to an overall scalar; to really get a permutation of the canonical or dual canonical bases we have to normalize properly.)

When we take the limit $q \mapsto 1$, the automorphisms of $U_q(\mathfrak{sl}_{a+b})$ above reduce to the following automorphisms of $U(\mathfrak{sl}_{a+b})$, the ordinary universal enveloping algebra of~$\mathfrak{sl}_{a+b}$:
\begin{gather*}
\eta(E_i) = F_{a+b+1-i},\qquad \eta(F_i) = E_{a+b+1-i},\qquad \eta(H_i) = -H_{a+b+1-i}, \\
\psi(E_i) = E_{a+b+1-i},\qquad \psi(F_i) = F_{a+b+1-i},\qquad \psi(H_i) = H_{a+b+1-i},
\end{gather*}
where here the $E_i$, $F_i$, $H_i$ for $1\leq i <a+b$ are the usual generators of $U(\mathfrak{sl}_{a+b})$ (see for instance \cite[Section~2.4]{lenart2007lusztig}). Recall that we can identify $E_i$, $F_i$, and~$H_i$ with matrices in $\mathfrak{sl}_{a+b}$: $E_i$ is the matrix with a $1$ in position $(i,i+1)$ and $0$'s elsewhere; $F_i$ is the matrix with a $1$ in position $(i+1,i)$ and $0$'s elsewhere; and $H_i$ is the matrix with a $1$ in $(i,i)$, a $-1$ in $(i+1,i+1)$, and $0$'s elsewhere. Thus we see that the automorphism $\eta$ of $U(\mathfrak{sl}_{a+b})$ is induced from the automorphism of~$\mathfrak{sl}_{a+b}$ given by conjugation by $\overline{w_0}$. Similarly, the automorphism $\psi$ of $U(\mathfrak{sl}_{2n})$ is induced from the automorphism $\phi_B\colon \mathfrak{sl}_{2n} \to \mathfrak{sl}_{2n}$ given by $\phi_B(A) \coloneqq B^{-1} \cdot \big({-}A^T\big) \cdot B$ for all $A\in \mathfrak{sl}_{2n}$, which is the Lie algebra version of the Lie group automorphism $\phi_B\colon {\rm SL}(2n)\to {\rm SL}(2n)$ defined earlier. But then we can note that acting by $\overline{w_0}$ on $R(a,a+b)_m \simeq V(m\omega_a)^{*}$ is clearly compatible with conjugation by $\overline{w_0}$ in the sense that $\overline{w_0} g \overline{w_0}^{-1} \cdot (\overline{w_0} v) = \overline{w_0} g v$ for all $g \in {\rm SL}(2n)$, $v\in R(a,a+b)_m$. So $\overline{w_0}$ agrees with $\eta$ up to a scalar, and hence permutes the dual canonical basis up to a sign. Similarly, $\widetilde{\phi}_B\colon R(n,2n)_m \to R(n,2n)_m$ is compatible with~$\phi_B$ in the sense that $\phi_B(g)\widetilde{\phi}_B(v) = \widetilde{\phi}_B(g v)$ for all $g \in {\rm SL}(2n)$, $v\in V(m\omega_n)^{*}$. So~$\phi_B$ agrees with $\psi$ up to a scalar, and hence permutes the dual canonical basis up to a sign.

To summarize the preceding, we have argued that $\overline{w_0}$ and $\widetilde{\phi}_{B}$ permute the dual canonical basis in some way (at least up to overall signs). But we want to conclude that these permutations correspond to complementation and transposition of tableaux. To do this, we note that the transition matrix between the dual canonical and standard monomial basis is upper unitriangular with respect to a certain order on tableaux. This is proved in a paper of Brundan~\cite[Theorem~26]{brundan2006dual}, where he in fact gives an explicit formula for this transition matrix in terms of Kazhdan--Lusztig polynomials (and he notes~\cite[Remark~10]{brundan2006dual} that his indexing of the canonical basis is consistent with that of Du~\cite{du1992canonical}). We know from Corollary~\ref{cor:std_monomial_actions} that $\overline{w_0}$ and $\widetilde{\phi}_{B}$ permute the standard monomial basis in the appropriate way. So finally we observe that if a permutation matrix is conjugated by an upper unitriangular matrix to another permutation matrix, then the two permutation matrices have to be the same (and this remains true if one of the permutation matrices is a priori only a permutation matrix up to an overall sign). Thus, the fact that $\overline{w_0}$ and $\widetilde{\phi}_{B}$ permute the standard monomial basis in the appropriate way in fact implies that they permute the dual canonical basis in the appropriate way.
\end{proof}

\subsection{Proof of bicyclic sieving for Grassmannian coordinate ring}

We are now ready to prove Theorem~\ref{thm:trcom_csiev}.

\begin{proof}[Proof of Theorem~\ref{thm:trcom_csiev}]
Define the following two matrices in ${\rm GL}(2n)$:
\begin{gather*}
 X \coloneqq \begin{pmatrix} \def\arraystretch{1.4}\begin{array}{@{}c|c@{}} \mathrm{antidiag}\big(i^{n-1},i^{n-1},\dots,i^{n-1}\big) & 0 \\ \hline 0 & \mathrm{Id}_{n} \end{array} \end{pmatrix} \in {\rm GL}(2n), \\
 Y \coloneqq \begin{pmatrix} \def\arraystretch{1.4}\begin{array}{@{}c|c@{}} \mathrm{Id}_{n} & 0 \\ \hline 0 & (-1)^{n-1}\, \mathrm{antidiag}\big(i^{n-1},i^{n-1},\dots,i^{n-1}\big) \end{array} \end{pmatrix} \in {\rm GL}(2n).
\end{gather*}
Thus,
\[ \overline{w_0} \cdot \chi^{n} = X \cdot Y. \]
It is easy to check that $X,Y \in {\rm SL}(2n)$ and that $\phi_B(Y) = X^{-1}$. So we have the following equality of algebra automorphisms $R(n,2n)\to R(n,2n)$:
\[ \widetilde{\phi}_B \cdot \overline{w_0} \cdot \chi^{n} = X^{-1} \cdot \widetilde{\phi}_{B} \cdot X. \]

Next we observe that the matrix $XBX$ has the block form
\[ \begin{pmatrix} \def\arraystretch{1.4}\begin{array}{@{}c|c@{}} 0 & \mathrm{antidiag}\big(\overbrace{i^{-(n-1)},-i^{-(n-1)},\dots,\pm i^{-(n-1)}}^{n}\big) \\ \hline \mathrm{antidiag}\big(\overbrace{{-}i^{-(n-1)},i^{-(n-1)},\dots,\mp i^{-(n-1)}}^{n}\big) & 0 \end{array} \end{pmatrix} \]
from which it is easy to see that
\[ (XBX) \cdot \chi = - \chi \cdot (XBX).\]
Then we can compute that
\begin{align*}
\phi_{B}\big(X \cdot i^{-1} \chi \cdot X^{-1}\big) &= B^{-1} \cdot \big(\big(X \cdot i^{-1} \chi\cdot X^{-1}\big)^{T}\big)^{-1} \cdot B \\
&= B^{-1} \cdot (X^T)^{-1} \cdot \big(\big(i^{-1} \chi\big)^{T}\big)^{-1} \cdot \big(\big(X^{-1}\big)^{T}\big)^{-1} \cdot B \\
&= B^{-1}X^{-1} \cdot i \chi \cdot XB \\
&= X\cdot (XBX)^{-1} \cdot i \chi \cdot (XBX) \cdot X^{-1} \\
&= X\cdot -i \chi \cdot X^{-1} \\
&= X\cdot i^{-1} \chi \cdot X^{-1}.
\end{align*}
In other words, $X^{-1} \cdot i^{-1} \chi \cdot X$ is a symplectic matrix.

By looking at the characteristic polynomial, we see that the eigenvalues of $i^{-1} \chi$ are all distinct; in fact they are precisely $-\zeta^{1/2},-\zeta^{1/2 +1},-\zeta^{1/2 +2}, \dots, -\zeta^{1/2 +(2n-1)}$, where $\zeta \coloneqq {\rm e}^{\pi {\rm i}/n}$ is a primitive $(2n)$th root of unity. Set $D_{q} \coloneqq q^{1/2-n} \, \mathrm{diag}\big(1,q,q^2,\dots,q^{2n-1}\big)$. With this notation, $i^{-1} \chi$ is conjugate to $D_{\zeta}$. Moreover, since $X \cdot i^{-1} \chi \cdot X^{-1}$ is symplectic, and it is diagonalizable, a basic result in symplectic linear algebra says this matrix is \emph{symplectically diagonalizable}; that is, we can find a symplectic matrix $S \in {\rm Sp}(2n)$ and a diagonal matrix $D$ for which we have $X \cdot i^{-1} \chi \cdot X^{-1}= S \cdot D \cdot S^{-1}$. In fact, by conjugating by an element of the hyperoctahedral group, we can assume that this diagonal matrix is $D=D_{\zeta}$. (Note that $D_{q}$ is in fact symplectic since its $i$th and $(2n+1-i)$th entries along the diagonal are inverses.) We then also clearly have that $X \cdot i^{-k} \chi^{k} \cdot X^{-1}= S \cdot D_{\zeta^{k}} \cdot S^{-1}$.

To complete the proof, we compute $\mathrm{tr}_{R(n,2n)_m}\big(\widetilde{\phi}_B \cdot \overline{w_0} \cdot i^{-k} \chi^{(n+k)}\big)$ in two ways. On the one hand,
\begin{align*}
 \mathrm{tr}_{R(n,2n)_m}\big(\widetilde{\phi}_B \cdot \overline{w_0} \cdot i^{-k} \chi^{(n+k)}\big) &= i^{knm}\!\cdot \mathrm{tr}_{R(n,2n)_m}\big(\widetilde{\phi}_B \cdot \overline{w_0} \cdot \chi^{(n+k)}\big)\\
&= i^{knm}\! \cdot \#\big\{T \!\in\! \mathrm{SSYT}(m^n,2n) \colon (\mathrm{Tp}\cdot \mathrm{Co})\cdot \mathrm{Pro}^{n+k}(T) = T\big\} \\
&= i^{knm}\! \cdot \#\big\{\pi \!\in\! \mathrm{PP}^m(n\times n) \colon (\mathrm{Tp}\cdot \mathrm{Co})\cdot \mathrm{Pro}^{n+k}(\pi) = \pi\big\}.
\end{align*}
Here the factor of $i^{knm}$ comes from the fact that $i^{-k} \, \mathrm{Id}_{2n} \in {\rm GL}(2n)$ scales each Pl\"{u}cker coordinate by $i^{kn}$ and hence scales elements of $R(n,2n)_m$ by $i^{knm}$. Meanwhile, the interpretation of the term $\mathrm{tr}_{R(n,2n)_m}\big(\widetilde{\phi}_B \cdot \overline{w_0} \cdot \chi^{(n+k)}\big)$ in terms of tableaux fixed by $(\mathrm{Tp}\cdot \mathrm{Co})\cdot \mathrm{Pro}^{n+k}$ follows from working in the dual canonical basis and recalling Theorem~\ref{thm:dual_canonical_actions}. Finally the interpretation in terms of plane partitions fixed by $(\mathrm{Tp}\cdot \mathrm{Co})\cdot \mathrm{Pro}^{n+k}$ comes from the bijection $\Psi$ between plane partitions and tableaux described in Appendix~\ref{sec:appendix}.

On the other hand, from our work above we have
\begin{align*}
\mathrm{tr}_{R(n,2n)_m}\big(\widetilde{\phi}_B \cdot \overline{w_0} \cdot i^{-k} \chi^{(n+k)}\big) &= \mathrm{tr}_{R(n,2n)_m}\big( X^{-1} \cdot \widetilde{\phi}_{B} \cdot X \cdot i^{-k} \chi^{k}\big) \\
&= \mathrm{tr}_{R(n,2n)_m}\big( \widetilde{\phi}_{B} \cdot X \cdot i^{-k} \chi^{k} \cdot X^{-1}\big) \\
&= \mathrm{tr}_{R(n,2n)_m}\big( \widetilde{\phi}_{B} \cdot S \cdot D_{\zeta^{k}} \cdot S^{-1} \big) \\
&= \mathrm{tr}_{R(n,2n)_m}\big(S\cdot \widetilde{\phi}_{B} \cdot D_{\zeta^{k}} \cdot S^{-1} \big) \\
&= \mathrm{tr}_{R(n,2n)_m}\big( \widetilde{\phi}_{B} \cdot D_{\zeta^{k}}\big).
\end{align*}
In general, the trace $\mathrm{tr}_{V(\lambda)} ( \phi \cdot D)$, where $\phi$ is the twist by an automorphism of a simple Lie group induced from a Dynkin diagram automorphism and $D$ is a torus element, is what is called a \emph{twining character}. The twining character formula, originally due to Jantzen~\cite{jantzen1973character} and later rediscovered for instance in~\cite{fuchs1997automorphisms, fuchs1996dynkin}, expresses such a twining character as an ordinary character of the so-called ``orbit Lie group''. In our case, that orbit Lie group would be the special orthogonal group ${\rm SO}(2n+1)$. But in fact, it is easy to compute $\mathrm{tr}_{R(n,2n)_m}\big( \widetilde{\phi}_{B} \cdot D_{\zeta^{k}}\big)$ directly. Indeed, as Kuperberg explained in~\cite[Section~4]{kuperberg1994minuscule}, just by considering the action on the standard monomial basis we can see~$\mathrm{tr}_{R(n,2n)_m}\big(\widetilde{\phi}_{B} \cdot D_{\zeta^{k}}\big) = \big(\zeta^{-k}\big)^{-n^2m/2} \, \mathrm{SymMac}\big(n,m; q \mapsto \zeta^{-k}\big)$. (Here we have~$\zeta^{-k}$ instead of $\zeta^{k}$ because $R(n,2n)_m\simeq V(m\omega_n)^*$ is isomorphic to the \emph{dual} representation.) But $\big(\zeta^{-k}\big)^{-n^2m/2} = i^{knm}$, so we conclude
\[\#\big\{\pi \in \mathrm{PP}^m(n\times n) \colon (\mathrm{Tp}\cdot \mathrm{Co})\cdot \mathrm{Pro}^{n+k}(\pi) = \pi\big\} = \mathrm{SymMac}\big(n,m; q \mapsto \zeta^{-k}\big).\]
Since $\zeta^{k}$ and $\zeta^{-k}$ are Galois conjugates and $\mathrm{SymMac}(n,m; q)\in\mathbb{Z}[q]$, in fact we have
\[\#\big\{\pi \in \mathrm{PP}^m(n\times n) \colon (\mathrm{Tp}\cdot \mathrm{Co})\cdot \mathrm{Pro}^{n+k}(\pi) = \pi\big\} = \mathrm{SymMac}\big(n,m; q \mapsto \zeta^{k}\big),\]
as claimed.
\end{proof}

\begin{Remark}
By now the significance (in combinatorics, algebra, geometry, et cetera) of the twisted cyclic shift acting on the Grassmannian is well appreciated. See the recent paper of Karp~\cite{karp2018moment} for a nice survey of many places in which the cyclic shift arises. Another related paper which studied the cyclic shift as well as involutive symmetries of the Grassmannian is the recent paper of Frieden~\cite{frieden2019affine}. In that paper, Frieden constructed an affine geometric crystal on the Grassmannian and in doing so showed that (a deformation of) the twisted cyclic shift tropicalizes to promotion of rectangular semistandard Young tableaux.
\end{Remark}

\begin{Remark}
As mentioned in Section~\ref{sec:intro}, there is another proof of cyclic sieving for plane partitions under promotion (Theorem~\ref{thm:ppart_csiev}) due to Shen and Weng~\cite{shen2018cyclic}. The main difference from Rhoades's proof is that, rather than use the Lusztig/Kashiwara dual canonical basis, Shen--Weng used a basis of the coordinate ring of the Grassmannian coming from its cluster algebra structure. Let us quickly review their setting. A ``cluster ensemble'' is a pair of an ``$\mathcal{X}$-cluster variety'' and an ``$\mathcal{A}$-cluster variety'' associated to a quiver, which are ``dual'' to one another. The ``Fock--Goncharov conjecture''~\cite{fock2009cluster} predicts that the tropical points of an $\mathcal{X}$-cluster variety parameterize a canonical basis of the coordinate ring of its dual $\mathcal{A}$-cluster variety, and vice-versa. Breakthrough work of Gross--Hacking--Keel--Kontsevich~\cite{gross2018canonical} establishes that the Fock--Goncharov conjecture holds as long as certain combinatorial conditions on the quiver are met. Under these conditions we have a canonical basis for the cluster algebra, the so-called ``theta basis''. Very roughly speaking, the Grassmannian carries the structure of both an $\mathcal{X}$- and an $\mathcal{A}$-cluster variety, and it is its own dual in a cluster ensemble. Shen--Weng verified the Gross--Hacking--Keel--Kontsevich combinatorial conditions for the quiver associated to the Grassmannian and so showed that the Fock--Goncharov conjecture holds in this case. (Their work is closely related to work of Rietsch and Williams~\cite{rietsch2017newton}, which also studied cluster duality for the Grassmannian.) Moreover, Shen--Weng showed that the twisted cyclic shift corresponds to an element of the ``cluster modular group'', a certain group of automorphisms of the cluster structure. The Fock--Goncharov conjecture says that the parametrization of the canonical cluster basis of one variety by the tropical points of its dual variety should be equivariant under the action of the cluster modular group. Shen--Weng showed that the action of the twisted cyclic shift on tropical points corresponds to promotion of plane partitions. They deduced that the theta basis of the coordinate ring of the Grassmannian is permuted by the twisted cyclic shift according to promotion of plane partitions (just like the dual canonical basis).

It is reasonable to ask how the involutive symmetries $\overline{w_0}$ and $\widetilde{\phi}_B$ behave on the theta basis of the coordinate ring of the Grassmannian (this could, for instance, yield a different proof of Theorem~\ref{thm:trcom_csiev}.) Understanding the behavior of the Dynkin diagram automorphism $\widetilde{\phi}_B$ seems tractable. Indeed, in the case $a=b=n$, the quiver $\Gamma_{a,a+b}$ (see~\cite[Section~2.4]{shen2018cyclic}) defining the Grassmannian cluster ensemble has a transposition symmetry $A_{i,j}\mapsto A_{j,i}$ (ignoring arrows between frozen vertices, which are largely irrelevant). Moreover, this $A_{i,j}\mapsto A_{j,i}$ symmetry exactly corresponds to the Pl\"{u}cker coordinate map $\Delta_I \mapsto \Delta_{-w_0(I)}$ which Lemma~\ref{lem:plucker_actions} says defines the action of $\widetilde{\phi}_B$. The behavior of~$\overline{w_0}$, however, is less clear to us. As Fraser explained in~\cite[Section~5]{fraser2017braid}, the twisted reflection is, together with the twisted cyclic shift, one of the most well-known and significant cluster automorphisms of the Grassmannian. (This follows from the combinatorics of weakly separated collections as elucidated in the seminal work of Scott~\cite{scott2006grassmannians} and Postnikov~\cite{postnikov2006total}.) However, the twisted reflection is, unlike the twisted cyclic shift, an ``orientation-reversing'' cluster automorphism, which means it cannot be an element of the cluster modular group. Hence, the Fock--Goncharov conjecture says nothing about the behavior of the twisted reflection on the theta basis. At the moment we have no good way of understanding the behavior of $\overline{w_0}$ on the theta basis.
\end{Remark}

\section{Promotion and rowmotion}

In this section we introduce another piecewise-linear operator on plane partitions called \emph{rowmotion}, which is different from but closely related to promotion. Let us very briefly review the history of rowmotion. \emph{Combinatorial rowmotion} is a certain invertible operator acting on the set of order ideals of any finite\footnote{Throughout, all posets will be finite.} poset~$P$ which has been studied by many authors over a~long period~\cite{armstrong2013uniform,cameron1995orbits, brouwer1974period, fonderflaass1993orbits, panyushev2009orbits, rush2013orbits, striker2012promotion}. \emph{Piecewise-linear rowmotion} is a generalization of combinatorial rowmotion which was introduced by Einstein and Propp~\cite{einstein2013combinatorial, einstein2014piecewise} about 5 years ago. Piecewise-linear rowmotion, as well as its further generalization \emph{birational rowmotion}, have already been the subject of a good deal of research~\cite{grinberg2015birational2,grinberg2016birational1, musiker2019paths}, with interesting connections to topics ranging from integrable systems to quiver representations~\cite{galashin2019rsystems, garver2018minuscule}. In the next section we will define piecewise-linear rowmotion acting on the $P$-partitions of an arbitrary poset~$P$; in this section we focus exclusively on plane partitions (which corresponds to~$P$ being the \emph{rectangle poset}, i.e., the product of two chains). Since our focus throughout will be on piecewise-linear rowmotion (as opposed to combinatorial or birational rowmotion), from now on we will drop the ``piecewise-linear'' adjective and speak simply of ``rowmotion''.

Our goal in this section is to show that versions of Theorems~\ref{thm:ppart_csiev}, \ref{thm:com_csiev}, \ref{thm:tr_csiev} and~\ref{thm:trcom_csiev} hold for rowmotion. So now let us formally define rowmotion and explain its relationship to promotion. As in the preceding sections, we fix the parameters $a$, $b$, and $m$ defining our set of plane partitions $\mathrm{PP}^{m}(a\times b)$; and sometimes (e.g., when we want to consider the transposition symmetry) we also assume that $a=b=n$.

We define rowmotion acting on $\mathrm{PP}^{m}(a\times b)$ as a composition of the piecewise-linear toggles $\tau_{i,j}\colon \mathrm{PP}^{m}(a\times b) \to \mathrm{PP}^{m}(a\times b)$ introduced in Section~\ref{sec:intro}. Recall that~$\tau_{i,j}$ and~$\tau_{i',j'}$ commute unless $(i,j)$ and $(i',j')$ are directly adjacent. For~$1 \leq k \leq a+b-1$, we define $\mathcal{R}_{k} \coloneqq \prod_{\substack{1\leq i \leq a, \\1\leq j \leq b, \\ i+j-1=k}}\tau_{i,j}$ to be the composition of all the toggles along the ``$k$th antidiagonal'' of our array (this composition is well-defined because these toggles commute).

We then define \emph{rowmotion} $\mathrm{Row} \colon \mathrm{PP}^{m}(a\times b) \to \mathrm{PP}^{m}(a\times b)$ as
\[ \mathrm{Row} \coloneqq \mathcal{R}_{a+b-1} \cdot \mathcal{R}_{a+b-2} \cdots \mathcal{R}_{2} \cdot \mathcal{R}_{1}.\]
Promotion was defined similarly but in terms of the diagonal toggles $\mathcal{F}_{k} \coloneqq \prod_{\substack{1\leq i \leq a, \\1\leq j \leq b, \\ j-i=k}} \tau_{i,j}$. Observe that promotion is a composition of the piecewise-linear toggles ``from left to right'', whereas rowmotion is a composition of the toggles ``from top to bottom''.

\begin{figure}[t]\centering
\begin{tikzpicture}
\def \x {1.4}

\node at (-0.5*\x,0.5*\x) {\bf (1,1)};
\node at (0.5*\x,0.5*\x) {\bf (1,2)};
\node at (1.5*\x,0.5*\x) {\bf (1,3)};
\node at (-0.5*\x,-0.5*\x) {\bf (2,1)};
\node at (0.5*\x,-0.5*\x) {\bf (2,2)};
\node at (1.5*\x,-0.5*\x) {\bf (2,3)};

\draw (-1*\x,1*\x) -- (2*\x,1*\x);
\draw (-1*\x,0*\x) -- (2*\x,0*\x);
\draw (-1*\x,-1*\x) -- (2*\x,-1*\x);

\draw (2*\x,-1*\x) -- (2*\x,1*\x);
\draw (1*\x,-1*\x) -- (1*\x,1*\x);
\draw (0*\x,-1*\x) -- (0*\x,1*\x);
\draw (-1*\x,-1*\x) -- (-1*\x,1*\x);

\draw[dashed, ->] (-1.75*\x,-0.5*\x) -- (3.05*\x,-0.5*\x);
\node at (3.25*\x,-0.5*\x) {$\mathcal{P}_2$};
\draw[dashed, ->] (-1.75*\x,0.5*\x) -- (3.05*\x,0.5*\x);
\node at (3.25*\x,0.5*\x) {$\mathcal{P}_1$};

\draw[dashed, ->] (-0.5*\x,1.25*\x) -- (-0.5*\x,-1.6*\x);
\node at (-0.5*\x,-1.75*\x) {$\mathcal{N}_1$};
\draw[dashed, ->] (0.5*\x,1.25*\x) -- (0.5*\x,-1.6*\x);
\node at (0.5*\x,-1.75*\x) {$\mathcal{N}_2$};
\draw[dashed, ->] (1.5*\x,1.25*\x) -- (1.5*\x,-1.6*\x);
\node at (1.5*\x,-1.75*\x) {$\mathcal{N}_3$};

\draw[dashed] (-1.75*\x,-0.75*\x) -- (1.1*\x,2.1*\x);
\node at (1.25*\x,2.25*\x) {$\mathcal{R}_1$};
\draw[dashed] (-1.25*\x,-1.25*\x) -- (1.6*\x,1.6*\x);
\node at (1.75*\x,1.75*\x) {$\mathcal{R}_2$};
\draw[dashed] (-0.25*\x,-1.25*\x) -- (2.1*\x,1.1*\x);
\node at (2.25*\x,1.25*\x) {$\mathcal{R}_3$};
\draw[dashed] (0.75*\x,-1.25*\x) -- (2.6*\x,0.6*\x);
\node at (2.75*\x,0.75*\x) {$\mathcal{R}_4$};

\draw[dashed] (0.25*\x,-1.25*\x) -- (-1.6*\x,0.6*\x);
\node at (-1.75*\x,0.75*\x) {$\mathcal{F}_{-1}$};
\draw[dashed] (1.25*\x,-1.25*\x) -- (-1.1*\x,1.1*\x);
\node at (-1.25*\x,1.25*\x) {$\mathcal{F}_{0}$};
\draw[dashed] (2.25*\x,-1.25*\x) -- (-0.6*\x,1.6*\x);
\node at (-0.75*\x,1.75*\x) {$\mathcal{F}_{1}$};
\draw[dashed] (2.25*\x,-0.25*\x) -- (-0.1*\x,2.1*\x);
\node at (-0.25*\x,2.25*\x) {$\mathcal{F}_{2}$};
\end{tikzpicture}
\caption{The various toggle compositions $\mathcal{F}_k$, $\mathcal{R}_k$, $\mathcal{P}_k$ and $\mathcal{N}_k$.} \label{fig:toggle_comps}
\end{figure}
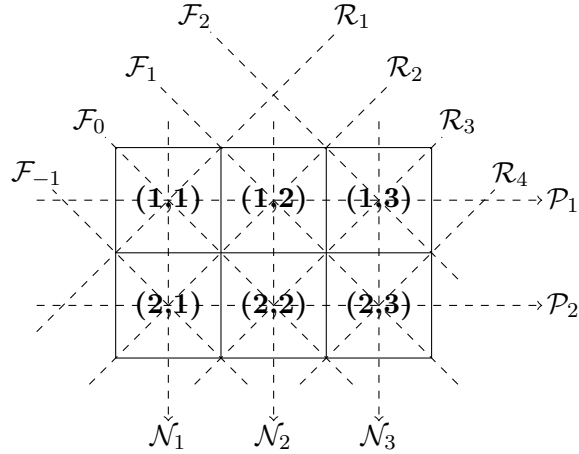

In a moment we will explain the precise relationship between promotion and rowmotion. First let us review a few other ways to express both promotion and rowmotion. In order to do this we must introduce some other compositions of toggles in addition to the~$\mathcal{F}_k$ and $\mathcal{R}_k$. For $1\leq k \leq a$ we define $\mathcal{P}_k \coloneqq \tau_{k,b} \cdot \tau_{k,b-1} \cdots \tau_{k,2} \cdot \tau_{k,1}$ and for $1\leq k \leq b$ we define $\mathcal{N}_k \coloneqq \tau_{a,k} \cdot \tau_{a-1,k} \cdots \tau_{2,k} \cdot \tau_{1,k}$. Note that while the toggles constituting $\mathcal{F}_k$ or~$\mathcal{R}_k$ all commute, this is not true for $\mathcal{P}_k$ or $\mathcal{N}_k$, so we have to be careful to specify the order of composition like we have just done. Relatedly, $\mathcal{F}_k=\mathcal{F}_k^{-1}$ and $\mathcal{R}_k=\mathcal{R}_k^{-1}$, but $\mathcal{P}_k \neq \mathcal{P}_k^{-1}$ and $\mathcal{N}_k \neq \mathcal{N}_k^{-1}$. Fig.~\ref{fig:toggle_comps} is a diagram depicting all of these various compositions of toggles. Let us explain our notation for these compositions, which perhaps appears quite strange at first. This notation derives from terminology due to Einstein and Propp~\cite[Section~2]{einstein2013combinatorial}: the~$\mathcal{F}_k$ are toggles along a \emph{file} of the rectangle; the~$\mathcal{R}_k$ are toggles along a \emph{rank} of the rectangle; the~$\mathcal{P}_k$ are toggles along a \emph{positive fiber} of the rectangle; and the~$\mathcal{N}_k$ are toggles along a~\emph{negative fiber} of the rectangle. (When we view the rectangle as a~poset, like we will do in the next section, the terms ``row'' and ``column'' do not cohere with the usual way of depicting a poset via its Hasse diagram, which is why Einstein and Propp avoided those terms.)

Both promotion and rowmotion can also be written in terms of~$\mathcal{P}_k$ (or~$\mathcal{N}_k$):

\begin{Proposition} \label{prop:pro_row_other_defs}
We have
\begin{align*}
\mathrm{Pro} =\mathcal{P}_1 \cdot \mathcal{P}_2 \cdots \mathcal{P}_{a-1} \cdot \mathcal{P}_a
 =\mathcal{N}^{-1}_{b} \cdot \mathcal{N}^{-1}_{b-1} \cdots \mathcal{N}^{-1}_{2} \cdot \mathcal{N}^{-1}_1,
\end{align*}
and
\begin{align*}
\mathrm{Row} =\mathcal{P}_{a} \cdot \mathcal{P}_{a-1} \cdots \mathcal{P}_{2} \cdot \mathcal{P}_{1}
 = \mathcal{N}_{b} \cdot \mathcal{N}_{b-1} \cdots \mathcal{N}_{2} \cdot \mathcal{N}_{1}.
\end{align*}
\end{Proposition}
\begin{proof}
This is explained, for instance, in~\cite[Section~2]{einstein2013combinatorial}. It follows easily from the commutativity properties of the toggles. As we will see in the next section, we can in fact define rowmotion to be the composition of all the toggles in the order of any linear extension of the underlying rectangle poset (and something similar is true for promotion).
\end{proof}

We now explain the relationship between promotion and rowmotion. This relationship was first discovered and explored by Striker and Williams~\cite{striker2012promotion}. It turns out that promotion and rowmotion are conjugate to one another. Moreover, there is a simple, explicit composition of toggles which conjugates $\mathrm{Row}$ to $\mathrm{Pro}$. Namely, define
\[D \coloneqq \big(\mathcal{P}^{-1}_a\big) \cdot \big(\mathcal{P}^{-1}_{a-1} \cdot \mathcal{P}^{-1}_a\big) \cdots \big(\mathcal{P}^{-1}_2 \cdot \mathcal{P}^{-1}_3 \cdots \mathcal{P}^{-1}_{a}\big) \cdot \big(\mathcal{P}^{-1}_1 \cdot \mathcal{P}^{-1}_2 \cdots \mathcal{P}^{-1}_a\big).\]
Then we have:

\begin{Lemma} \label{lem:pro_row_conjugacy}
We have $D \cdot \mathrm{Row} \cdot D^{-1} = \mathrm{Pro}$.
\end{Lemma}
\begin{proof}
This was essentially proved by Striker--Williams~\cite[Theorem 5.4]{striker2012promotion}. They were working at the level of \emph{combinatorial} rowmotion (i.e., the case $m=1$) but they only used the facts that the toggles are involutions and that non-adjacent toggles commute (in other words, they were really working in the corresponding ``right angled Coxeter group''). At any rate, this lemma follows from the description of promotion and rowmotion in terms of the $\mathcal{P}_k$ given in Proposition~\ref{prop:pro_row_other_defs}, together with the observation that $\mathcal{P}_i$ and $\mathcal{P}_j$ commute unless $|i-j|=1$.
\end{proof}

\begin{Corollary} \label{cor:ppart_csiev}
For any $k \in \mathbb{Z}$, we have
\[\#\big\{\pi \in \mathrm{PP}^{m}(a\times b)\colon \mathrm{Row}^k(\pi)=\pi\big\} = \mathrm{Mac}\big(a,b,m;q \mapsto \zeta^k\big),\]
where $\zeta \coloneqq {\rm e}^{2\pi {\rm i}/(a+b)}$ is a primitive $(a+b)$th root of unity.
\end{Corollary}
\begin{proof}
This follows from combining Theorem~\ref{thm:ppart_csiev} and Lemma~\ref{lem:pro_row_conjugacy}.
\end{proof}

Observe that Corollary~\ref{cor:ppart_csiev} implies that $\mathrm{Row}$ has order dividing $a+b$.

Next we want to understand how rowmotion interacts with the involutive symmetries of complementation and transposition. First of all, observe that $\mathrm{Co} \cdot \mathrm{Row} = \mathrm{Row}^{-1}\cdot \mathrm{Co}$, while $\mathrm{Tp} \cdot \mathrm{Row} = \mathrm{Row} \cdot \mathrm{Tp}$ (which is slightly different than for promotion).

In order to count fixed points of elements of the group $\langle \mathrm{Row}, \mathrm{Co}, \mathrm{Tp}\rangle$ acting on the plane partitions in~$\mathrm{PP}^{m}(a\times b)$, we need to understand how the conjugating map~$D$ interacts with complementation and transposition. First, we give a simple proposition about how $\mathrm{Co}$ and $\mathrm{Tp}$ interact with the individual toggles and with the~$\mathcal{P}_k$:

\begin{Proposition} \label{prop:co_tr_toggle_rels}
For $1 \leq i \leq a$ and $1 \leq j \leq b$, we have
\[ \mathrm{Co} \cdot \tau_{i,j} = \tau_{a+1-i,b+1-j} \cdot \mathrm{Co}\]
and hence for $1\leq k \leq a$ we have
\[ \mathrm{Co} \cdot \mathcal{P}_{k} = \mathcal{P}_{a+1-k}^{-1} \cdot \mathrm{Co}.\]
Similarly, if $a=b=n$, then for $1 \leq i,j \leq n$ we have
\[ \mathrm{Tp} \cdot \tau_{i,j} = \tau_{j,i} \cdot \mathrm{Tp}\]
and hence for $1\leq k \leq n$ we have
\[ \mathrm{Tp} \cdot \mathcal{P}_{k} = \mathcal{N}_{k} \cdot \mathrm{Tp}.\]
\end{Proposition}
\begin{proof}
The statements about how to commute the $\tau_{i,j}$ past $\mathrm{Co}$ or $\mathrm{Tp}$ are immediate from the definitions of complementation, transposition, and the toggles. The statements about the $\mathrm{P}_k$ then follow easily.
\end{proof}

Now we can explain how $D$ and $\mathrm{Co}$ interact.

\begin{Lemma} \label{lem:comp_d_conjugacy}
We have $D \cdot \mathrm{Co} \cdot \mathrm{Row}^{-(a+1)} = \mathrm{Co} \cdot D$.
\end{Lemma}
\begin{proof}
From Proposition~\ref{prop:co_tr_toggle_rels} we get that
\[D \cdot \mathrm{Co} = \mathrm{Co} \cdot (\mathcal{P}_1) \cdot (\mathcal{P}_{2} \cdot \mathcal{P}_1) \cdots (\mathcal{P}_{a-1} \cdot \mathcal{P}_{a-2} \cdots \mathcal{P}_{1}) \cdot (\mathcal{P}_a \cdot \mathcal{P}_{a-1} \cdots \mathcal{P}_1).\]
Hence, using the description of $\mathrm{Row}$ in terms of $\mathcal{P}_k$ from Proposition~\ref{prop:pro_row_other_defs}, as well as the fact that~$\mathcal{P}_i$ and~$\mathcal{P}_j$ commute unless $|i-j|=1$, denoting $Z \coloneqq D \cdot \mathrm{Co} \cdot \mathrm{Row}^{-(a+1)}$, we get
\begin{align*}
Z & = \mathrm{Co} \cdot (\mathcal{P}_1) \cdot (\mathcal{P}_{2} \cdot \mathcal{P}_1) \cdots (\mathcal{P}_{a-1} \cdot \mathcal{P}_{a-2} \cdots \mathcal{P}_{1}) \cdot (\mathcal{P}_a \cdot \mathcal{P}_{a-1} \cdots \mathcal{P}_1) \cdot \mathrm{Row}^{-(a+1)} \\
&= \mathrm{Co} \cdot (\mathcal{P}_1) \cdot (\mathcal{P}_{2} \cdot \mathcal{P}_1) \cdots (\mathcal{P}_{a-1} \cdot \mathcal{P}_{a-2} \cdots \mathcal{P}_{1}) \cdot \mathrm{Row}^{-a} \\
&= \mathrm{Co} \cdot (\mathcal{P}_1) \cdot (\mathcal{P}_{2} \cdot \mathcal{P}_1) \cdots (\mathcal{P}_{a-2} \cdot \mathcal{P}_{a-3} \cdots \mathcal{P}_{1}) \cdot \big(\mathcal{P}_{a}^{-1}\big) \cdot \mathrm{Row}^{-(a-1)} \\
&= \mathrm{Co} \cdot \big(\mathcal{P}_{a}^{-1}\big) \cdot (\mathcal{P}_1) \cdot (\mathcal{P}_{2} \cdot \mathcal{P}_1) \cdots (\mathcal{P}_{a-2} \cdot \mathcal{P}_{a-3} \cdots \mathcal{P}_{1}) \cdot \mathrm{Row}^{-(a-1)} \\
&= \mathrm{Co} \cdot \big(\mathcal{P}_{a}^{-1}\big) \cdot (\mathcal{P}_1) \cdot (\mathcal{P}_{2} \cdot \mathcal{P}_1) \cdots (\mathcal{P}_{a-3} \cdot \mathcal{P}_{a-4} \cdots \mathcal{P}_{1}) \cdot \big(\mathcal{P}_{a-1}^{-1} \cdot \mathcal{P}_{a}^{-1}\big) \cdot \mathrm{Row}^{-(a-2)} \\
&= \mathrm{Co} \cdot \big(\mathcal{P}_{a}^{-1}\big) \cdot \big(\mathcal{P}_{a-1}^{-1} \cdot \mathcal{P}_{a}^{-1}\big) \cdot (\mathcal{P}_1) \cdot (\mathcal{P}_{2} \cdot \mathcal{P}_1) \cdots (\mathcal{P}_{a-3} \cdot \mathcal{P}_{a-4} \cdots \mathcal{P}_{1}) \cdot \mathrm{Row}^{-(a-2)}= \cdots \\
& = \mathrm{Co} \cdot \big(\mathcal{P}^{-1}_a\big) \cdot \big(\mathcal{P}^{-1}_{a-1} \cdot \mathcal{P}^{-1}_a\big) \cdots \big(\mathcal{P}^{-1}_2 \cdot \mathcal{P}^{-1}_3 \cdots \mathcal{P}^{-1}_{a}\big) \cdot \big(\mathcal{P}^{-1}_1 \cdot \mathcal{P}^{-1}_2 \cdots \mathcal{P}^{-1}_a\big) \\
&= \mathrm{Co} \cdot D,
\end{align*}
as claimed.
\end{proof}

\begin{Corollary} \label{cor:com_csiev}
For any $k \in \mathbb{Z}$, we have
\[\#\big\{\pi \in \mathrm{PP}^{m}(a\times b)\colon \mathrm{Co}\cdot \mathrm{Row}^k(\pi)=\pi\big\} = \#\big\{\pi \in \mathrm{PP}^{m}(a\times b)\colon \mathrm{Co}\cdot \mathrm{Pro}^k(\pi)=\pi\big\},\]
with an explicit formula for this number given by Theorem~{\rm \ref{thm:com_csiev}}.
\end{Corollary}
\begin{proof}
From Lemmas~\ref{lem:pro_row_conjugacy} and~\ref{lem:comp_d_conjugacy} we have
\[ D \cdot \mathrm{Co} \cdot \mathrm{Row}^{k} \cdot D^{-1} = \mathrm{Co}\cdot D\cdot \mathrm{Row}^{k+(a+1)}\cdot D^{-1} = \mathrm{Co}\cdot\mathrm{Pro}^{k+(a+1)}.\]
In other words, $\mathrm{Co} \cdot \mathrm{Row}^{k}$ and $\mathrm{Co} \cdot \mathrm{Pro}^{k+(a+1)}$ are conjugate, and hence in particular have the same number of fixed points. But by looking at the explicit formula in Theorem~\ref{thm:com_csiev}, we see that whether $a$ is even or odd, the number of fixed points of $\mathrm{Co} \cdot \mathrm{Pro}^{k+(a+1)}$ and of $\mathrm{Co} \cdot \mathrm{Pro}^{k}$ are the same, and hence the corollary follows.
\end{proof}

In order to explain how $D$ and $\mathrm{Tp}$ interact, we a need a few more preparatory results. First we need to explain how $\mathrm{Co}$ actually can be written as a composition of toggles.

\begin{Lemma} \label{lem:co_as_toggles}
We have
\[ \mathrm{Co} = (\mathcal{F}_{-a+1}) \cdot (\mathcal{F}_{-a+2} \cdot \mathcal{F}_{-a+1}) \cdots (\mathcal{F}_{b-2} \cdots \mathcal{F}_{-a+2} \cdot \mathcal{F}_{-a+1}) \cdot (\mathcal{F}_{b-1} \cdots \mathcal{F}_{-a+2} \cdot \mathcal{F}_{-a+1}) .\]
\end{Lemma}
\begin{proof}
This follows from various results in Appendix~\ref{sec:appendix}. Namely, in Definition~\ref{def:pro_evac}, dual eva\-cuation $\varepsilon^{*}$ is defined as the corresponding composition of Bender--Knuth involutions; in Proposition~\ref{prop:bk_toggles} it is shown that these Bender--Knuth involutions correspond to the diagonal toggles~$\mathcal{F}_k$; and in Proposition~\ref{prop:tableau_ops} it is shown that complementation of plane partitions corresponds to dual evacuation of semistandard tableaux under the bijection $\Psi\colon \mathrm{PP}^{m}(a\times b) \to \mathrm{SSYT}(m^a,a+b)$ studied in the appendix.
\end{proof}

\begin{Remark}
In addition to the description of $\mathrm{Co}$ in Lemma~\ref{lem:co_as_toggles}, we also have that
\[ \mathrm{Co} = (\mathcal{R}_{1}) \cdot (\mathcal{R}_{2} \cdot \mathcal{R}_{1}) \cdots (\mathcal{R}_{a+b-2} \cdots \mathcal{R}_{2} \cdot \mathcal{R}_{1}) \cdot (\mathcal{R}_{a+b-1} \cdots \mathcal{R}_{2} \cdot \mathcal{R}_{1}).\]
This follows from a ``reciprocity'' property of rowmotion that was established by Grinberg and Roby~\cite[Theorem 32]{grinberg2015birational2}. (We discuss this reciprocity property more in the next section; see Theorem~\ref{thm:reciprocity}.) Furthermore, there are other similar ways of writing $\mathrm{Co}$ as a composition of toggles which one can obtain by considering evacuation $\varepsilon$ or by considering the reciprocity property for $\mathrm{Row}^{-1}$. The fact that complementation is an involution gives even a few more ways of writing it as a composition of toggles.
\end{Remark}

\begin{Proposition} \label{prop:co_tech_res}
We have
\[\mathrm{Co} \cdot D =(\mathcal{N}_{1} \cdot \mathcal{N}_{2} \cdots \mathcal{N}_{b-2} \cdot \mathcal{N}_{b-1}) \cdots (\mathcal{N}_1 \cdot \mathcal{N}_2) \cdot (\mathcal{N}_1).\]
\end{Proposition}
\begin{proof}
We claim that
\begin{equation} \label{eq:fs}
(\mathcal{F}_{-a+1}) \cdot (\mathcal{F}_{-a+2} \cdot \mathcal{F}_{-a+1}) \cdots (\mathcal{F}_{b-2} \cdots \mathcal{F}_{-a+2} \cdot \mathcal{F}_{-a+1}) \cdot (\mathcal{F}_{b-1} \cdots \mathcal{F}_{-a+2} \cdot \mathcal{F}_{-a+1})
\end{equation}
is equal to
\begin{equation} \label{eq:ns_ps}
(\mathcal{N}_{1} \cdots \mathcal{N}_{b-2} \cdot \mathcal{N}_{b-1}) \cdots (\mathcal{N}_{1} \cdot \mathcal{N}_{2}) \cdot (\mathcal{N}_{1}) \cdot (\mathcal{P}_a \cdots \mathcal{P}_2 \cdot \mathcal{P}_1) \cdots (\mathcal{P}_{a} \cdot \mathcal{P}_{a-1}) \cdot (\mathcal{P}_a).
\end{equation}
By Lemma~\ref{lem:co_as_toggles}, establishing this claim proves the proposition. To show that~\eqref{eq:fs} and~\eqref{eq:ns_ps} are equal, first we observe that when we expand~\eqref{eq:ns_ps} as a product of toggles, a toggle $\tau_{i,j}$ appears $i-j+b$ times; this is the same number of times as when we expand~\eqref{eq:fs} as a product of toggles. Moreover, if we read the toggles in~\eqref{eq:ns_ps} from right to left, then whenever we see a toggle $\tau_{i,j}$ for the $k$th time, we have already seen the toggles $\tau_{i+1,j}$ (assuming $i\neq a$) and $\tau_{i,j-1}$ (assuming $j\neq 1$) exactly $k$ times. This implies that we can indeed commute the toggles that make up~\eqref{eq:ns_ps} to fit the form of~\eqref{eq:fs}, as claimed.
\end{proof}

Now we can explain how $D$ and $\mathrm{Tp}$ interact.

\begin{Lemma} \label{lem:tr_d_conjugacy}
If $a=b=n$, then we have
\[ D \cdot \mathrm{Tp} \cdot \mathrm{Row}^{n} = \mathrm{Tp} \cdot \mathrm{Co} \cdot D. \]
\end{Lemma}
\begin{proof}
From Proposition~\ref{prop:co_tr_toggle_rels} we get that
\[D \cdot \mathrm{Tp} = \mathrm{Tp} \cdot \big(\mathcal{N}^{-1}_n\big) \cdot \big(\mathcal{N}^{-1}_{n-1} \cdot \mathcal{N}^{-1}_n\big) \cdots \big(\mathcal{N}^{-1}_{2} \cdot \mathcal{N}^{-1}_{3} \cdots \mathcal{N}^{-1}_{n}\big) \cdot \big(\mathcal{N}^{-1}_{1} \cdot \mathcal{N}^{-1}_{2} \cdots \mathcal{N}^{-1}_{n}\big).\]
Using the description of $\mathrm{Row}$ in terms of $\mathcal{N}_k$ from Proposition~\ref{prop:pro_row_other_defs}, as well as the fact that $\mathcal{N}_i$ and $\mathcal{N}_j$ commute unless $|i-j|=1$, denoting $Z \coloneqq D \cdot \mathrm{Tp} \cdot \mathrm{Row}^{n}$, we get
\begin{align*}
Z &= \mathrm{Tp} \cdot \big(\mathcal{N}^{-1}_n\big) \cdot \big(\mathcal{N}^{-1}_{n-1} \cdot \mathcal{N}^{-1}_n\big) \cdots \big(\mathcal{N}^{-1}_{2} \cdot \mathcal{N}^{-1}_{3} \cdots \mathcal{N}^{-1}_{n}\big) \cdot \big(\mathcal{N}^{-1}_{1} \cdot \mathcal{N}^{-1}_{2} \cdots \mathcal{N}^{-1}_{n}\big) \cdot \mathrm{Row}^{n} \\
&= \mathrm{Tp} \cdot \big(\mathcal{N}^{-1}_n\big) \cdot \big(\mathcal{N}^{-1}_{n-1} \cdot \mathcal{N}^{-1}_n\big) \cdots \big(\mathcal{N}^{-1}_{2} \cdot \mathcal{N}^{-1}_{3} \cdots \mathcal{N}^{-1}_{n}\big) \cdot \mathrm{Row}^{n-1} \\
&= \mathrm{Tp} \cdot \big(\mathcal{N}^{-1}_n\big) \cdot \big(\mathcal{N}^{-1}_{n-1} \cdot \mathcal{N}^{-1}_n\big) \cdots \big(\mathcal{N}^{-1}_{3} \cdot \mathcal{N}^{-1}_{4} \cdots \mathcal{N}^{-1}_{n}\big) \cdot \mathcal{N}_{1} \cdot \mathrm{Row}^{n-2} \\
&= \mathrm{Tp} \cdot (\mathcal{N}_{1}) \cdot \big(\mathcal{N}^{-1}_n\big) \cdot \big(\mathcal{N}^{-1}_{n-1} \cdot \mathcal{N}^{-1}_n\big) \cdots \big(\mathcal{N}^{-1}_{3} \cdot \mathcal{N}^{-1}_{4} \cdots \mathcal{N}^{-1}_{n}\big) \cdot \mathrm{Row}^{n-2} \\
&= \mathrm{Tp} \cdot (\mathcal{N}_{1}) \cdot \big(\mathcal{N}^{-1}_n\big) \cdot \big(\mathcal{N}^{-1}_{n-1} \cdot \mathcal{N}^{-1}_n\big) \cdots \big(\mathcal{N}^{-1}_{4} \cdot \mathcal{N}^{-1}_{5} \cdots \mathcal{N}^{-1}_{n}\big) \cdot (\mathcal{N}_{2} \cdot \mathcal{N}_1) \cdot \mathrm{Row}^{n-3} \\
&= \mathrm{Tp} \cdot (\mathcal{N}_{1}) \cdot (\mathcal{N}_{2} \cdot \mathcal{N}_1) \cdot \big(\mathcal{N}^{-1}_n\big) \cdot \big(\mathcal{N}^{-1}_{n-1} \cdot \mathcal{N}^{-1}_n\big) \cdots \big(\mathcal{N}^{-1}_{4} \cdot \mathcal{N}^{-1}_{5} \cdots \mathcal{N}^{-1}_{n}\big) \cdot \mathrm{Row}^{n-3} \\
&= \cdots = \mathrm{Tp} \cdot (\mathcal{N}_{1}) \cdot (\mathcal{N}_{2} \cdot \mathcal{N}_1) \cdots (\mathcal{N}_{n-1} \cdot \mathcal{N}_{n-2} \cdots \mathcal{N}_{2} \cdot \mathcal{N}_1) \\
&= \mathrm{Tp} \cdot (\mathcal{N}_{1} \cdot \mathcal{N}_{2} \cdots \mathcal{N}_{n-2} \cdot \mathcal{N}_{n-1}) \cdots (\mathcal{N}_1 \cdot \mathcal{N}_2) \cdot (\mathcal{N}_1)\\
&= \mathrm{Tp}\cdot \mathrm{Co} \cdot D,
\end{align*}
where in the last line we applied Proposition~\ref{prop:co_tech_res}.
\end{proof}

\begin{Corollary} \label{cor:tr_csiev}
For any $k \in \mathbb{Z}$, we have
\[\#\big\{\pi \in \mathrm{PP}^{m}(n\times n)\colon \mathrm{Tp} \cdot \mathrm{Row}^{k}(\pi)=\pi\big\} = \mathrm{SymMac}\big(n,m;q \mapsto \zeta^k\big),\]
where $\zeta \coloneqq {\rm e}^{\pi {\rm i}/n}$ is a primitive $(2n)$th root of unity.
\end{Corollary}
\begin{proof}
From Lemmas~\ref{lem:pro_row_conjugacy} and~\ref{lem:comp_d_conjugacy} we have
\[ D \cdot \mathrm{Tp} \cdot \mathrm{Row}^{k} \cdot D^{-1} = \mathrm{Tp}\cdot \mathrm{Co}\cdot D\cdot \mathrm{Row}^{k-n}\cdot D^{-1} = \mathrm{Tp}\cdot \mathrm{Co}\cdot\mathrm{Pro}^{k-n}.\]
In other words, $\mathrm{Tp} \cdot \mathrm{Row}^{k}$ and $ \mathrm{Tp}\cdot \mathrm{Co}\cdot\mathrm{Pro}^{k-n}$ are conjugate, and hence in particular have the same number of fixed points. By Theorem~\ref{thm:trcom_csiev}, this number is the claimed evaluation of $\mathrm{SymMac}(n,m,q)$ (where we note that $\zeta^{k-2n}=\zeta^{k}$ since $\zeta^{2n}=1$).
\end{proof}

Finally, we can explain how $D$ and $\mathrm{Tp} \cdot \mathrm{Co}$ interact.

\begin{Lemma} \label{lem:trco_d_conjugacy}
If $a=b=n$, then we have
\[ D \cdot \mathrm{Tp} \cdot \mathrm{Co}\cdot \mathrm{Row}^{-1} = \mathrm{Tp} \cdot D. \]
\end{Lemma}
\begin{proof}
By Lemmas~\ref{lem:comp_d_conjugacy} and~\ref{lem:tr_d_conjugacy} we have
\begin{align*}
D \cdot \mathrm{Tp} \cdot \mathrm{Co}\cdot \mathrm{Row}^{-1} &= \mathrm{Tp} \cdot \mathrm{Co} \cdot D \cdot \mathrm{Row}^{-n} \cdot \mathrm{Co}\cdot \mathrm{Row}^{-1} \\
&= \mathrm{Tp} \cdot \mathrm{Co} \cdot D \cdot \mathrm{Co} \cdot \mathrm{Row}^{n-1} \\
&= \mathrm{Tp} \cdot \mathrm{Co} \cdot \mathrm{Co} \cdot D \cdot \mathrm{Row}^{n+1} \cdot \mathrm{Row}^{n-1} \\
&= \mathrm{Tp} \cdot D,
\end{align*}
where we used that $\mathrm{Co}$ is an involution and $\mathrm{Row}$ has order dividing $2n$.
\end{proof}

\begin{Corollary} \label{cor:trcom_csiev}
For any $k \in \mathbb{Z}$, we have
\[\#\big\{\pi \in \mathrm{PP}^{m}(n\times n)\colon \mathrm{Tp} \cdot \mathrm{Co} \cdot \mathrm{Row}^k(\pi)=\pi\big\} = \mathrm{SymMac}\big(n,m;q \mapsto (-1)^{k+1}\big).\]
\end{Corollary}
\begin{proof}
From Lemmas~\ref{lem:pro_row_conjugacy} and~\ref{lem:trco_d_conjugacy} we have
\[ D \cdot \mathrm{Tp} \cdot \mathrm{Co} \cdot \mathrm{Row}^{k} \cdot D^{-1} = \mathrm{Tp} \cdot D \cdot \mathrm{Row}^{k+1}\cdot D^{-1} = \mathrm{Tp}\cdot \mathrm{Pro}^{k+1}.\]
In other words, $\mathrm{Tp} \cdot \mathrm{Co} \cdot \mathrm{Row}^{k}$ and $ \mathrm{Tp}\cdot \mathrm{Pro}^{k+1}$ are conjugate, and hence in particular have the same number of fixed points. By Theorem~\ref{thm:tr_csiev}, this number is the claimed evaluation of $\mathrm{SymMac}(n,m,q)$.
\end{proof}

In direct analogy with what we showed for promotion, Corollaries~\ref{cor:ppart_csiev}, \ref{cor:com_csiev}, \ref{cor:tr_csiev} and~\ref{cor:trcom_csiev} together imply that for any element $g \in \langle \mathrm{Row}, \mathrm{Co}, \mathrm{Tp}\rangle$, the number of plane partitions in $\mathrm{PP}^{m}(n\times n)$ fixed by $g$ is given by some kind of sieving phenomenon evaluation of a nice polynomial at a root of unity.

\section{Rowmotion for triangular posets} \label{sec:row_tri}

Since orbit structures are our main interest in this paper, and since, as we explained in the last section, rowmotion is conjugate to promotion, it might not be clear why we care about rowmotion at all. The reason we do care about rowmotion is that rowmotion can be defined as an action on the set of $P$-partitions of any poset~$P$. Rowmotion acting on plane partitions corresponds to taking $P$ to be the rectangle poset. Moreover, rowmotion still behaves remarkably well on the $P$-partitions of other posets $P$ besides the rectangle poset, especially certain nice posets coming from Lie theory (namely, \emph{minuscule posets} and \emph{root posets of coincidental type}). In~\cite{hopkins2019minuscule} we made a number of cyclic sieving conjectures regarding the action of rowmotion on the $P$-partitions of these other nice posets $P$. The major examples of the nice posets $P$, beyond the rectangle itself, are certain ``triangular'' posets. As we will explain in this section (following Grinberg--Roby~\cite{grinberg2015birational2}), these triangular posets can be obtained from the rectangle by enforcing certain symmetries. Understanding the behavior of rowmotion on these triangular posets was our original motivation for studying how rowmotion interacts with the symmetries of complementation and transposition. As we will see, while the results we obtained above concerning the interaction of rowmotion with these symmetries do not directly imply anything about cyclic sieving for the triangular posets, morally they are very closely related to our conjectures from~\cite{hopkins2019minuscule}.

So now we define rowmotion for arbitrary posets. We assume the reader is familiar with the basics of posets as laid out for instance in~\cite[Chapter 3]{stanley2012ec1}. Let $P$ be a finite partially ordered set. We use~$\leq$ for the partial order of $P$ and~$\lessdot$ for its cover relation. A \emph{$P$-partition of height~$m$} is a weakly order-preserving map $\pi \colon P \to \{0,1,\dots,m\}$, i.e., one for which $p\leq q\in P$ implies that~$\pi(p) \leq \pi(q)$. We use~$\mathrm{PP}^{m}(P)$ to denote the set of $P$-partitions of height $m$. For any element~$p\in P$, the \emph{piecewise-linear toggle at $p$} is the involution $\tau_p \colon \mathrm{PP}^{m}(P)\to \mathrm{PP}^{m}(P)$ defined by
\[(\tau_{p}\pi) (q) \coloneqq \begin{cases}\pi(q) &\textrm{if $p\neq q$}, \\ \min (\{\pi(r)\colon p \lessdot r\}) + \max (\{\pi(r)\colon r \lessdot p\}) - \pi(p) &\textrm{if $p=q$}, \end{cases}\]
with the conventions that $\min (\varnothing) = m$ and $\max (\varnothing)=0$. Note that toggles $\tau_p$ and $\tau_q$ commute unless there is a cover relation between $p$ and $q$.

We then define \emph{rowmotion} $\mathrm{Row} \colon \mathrm{PP}^{m}(P)\to\mathrm{PP}^{m}(P)$ by
\[ \mathrm{Row} \coloneqq \tau_{p_1} \cdot \tau_{p_2} \cdots \tau_{p_{n-1}} \cdot \tau_{p_n},\]
where $p_1,p_2,\dots,p_n$ is any linear extension of~$P$. The commutativity properties of the toggles imply that this definition does not depend on the choice of linear extension.

\begin{figure}[t]\centering
 \begin{tikzpicture}[scale=0.8]
	\SetFancyGraph
	\Vertex[NoLabel,x=0,y=0]{1}
	\Vertex[NoLabel,x=1,y=1]{2}
	\Vertex[NoLabel,x=2,y=2]{3}
	\Vertex[NoLabel,x=3,y=3]{4}
	\Vertex[NoLabel,x=-1,y=1]{5}
	\Vertex[NoLabel,x=0,y=2]{6}
	\Vertex[NoLabel,x=1,y=3]{7}
	\Vertex[NoLabel,x=2,y=4]{8}
	\Vertex[NoLabel,x=-2,y=2]{9}
	\Vertex[NoLabel,x=-1,y=3]{10}
	\Vertex[NoLabel,x=0,y=4]{11}
	\Vertex[NoLabel,x=1,y=5]{12}
	\Edges[style={thick}](1,2)
	\Edges[style={thick}](2,3)
	\Edges[style={thick}](3,4)
	\Edges[style={thick}](1,5)
	\Edges[style={thick}](2,6)
	\Edges[style={thick}](3,7)
	\Edges[style={thick}](4,8)
	\Edges[style={thick}](5,6)
	\Edges[style={thick}](6,7)
	\Edges[style={thick}](7,8)
	\Edges[style={thick}](5,9)
	\Edges[style={thick}](6,10)
	\Edges[style={thick}](7,11)
	\Edges[style={thick}](8,12)
	\Edges[style={thick}](9,10)
	\Edges[style={thick}](10,11)
	\Edges[style={thick}](11,12)
	\node[anchor=north] at (1) {$(a,b)$};
	\node[anchor=south] at (12) {$(1,1)$};
	\node[anchor=west] at (4) {$(1,b)$};
	\node[anchor=east] at (9) {$(a,1)$};
\end{tikzpicture}
\caption{The rectangle poset $a\times b$.} \label{fig:rect_poset}
\end{figure}
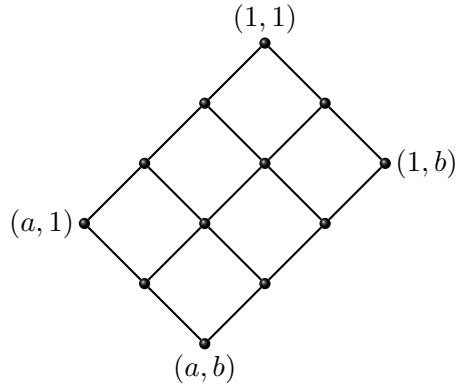

The \emph{rectangle poset}, denoted $a \times b$, is the Cartesian product of an $a$-element chain and a~$b$-element chain. Rowmotion on the rectangle poset is the same as rowmotion of~$a\times b$ plane partitions. However, note that we are working with order-\emph{preserving} maps in order to match the conventions of~\cite{einstein2013combinatorial} and~\cite{hopkins2019minuscule}; and on the other hand in order to match the traditional conventions for plane partitions we put the maximal entry of a plane partition in its upper-left corner. Thus to satisfy all of our conventions we need to view the rectangle poset $P=a\times b$ as the set~$\{(i,j)\colon 1\leq i \leq a, 1\leq j \leq b\}$ with the ``backwards'' partial order $(i,j) \leq (i',j')$ if and only if $i' \leq i$ and $j' \leq j$. Fig.~\ref{fig:rect_poset} depicts the Hasse diagram of the rectangle poset: observe how we have simply rotated the usual matrix coordinates $45^{\circ}$ clockwise. With this convention for the naming of elements of the rectangle poset, our notation $\mathrm{PP}^{m}(a \times b)$ is consistent whether we think of this as a set of plane partitions $(\pi_{i,j})$ or of $P$-partitions $\pi(i,j)$. And of course, as mentioned, the two definitions of rowmotion acting on $\mathrm{PP}^{m}(a \times b)$ agree as well.

\begin{figure}[t]\centering
\begin{tikzpicture}[scale=0.75]
	\SetFancyGraph
	\Vertex[NoLabel,x=0,y=0]{1}
	\Vertex[NoLabel,x=1,y=1]{2}
	\Vertex[NoLabel,x=2,y=2]{3}
	\Vertex[NoLabel,x=3,y=3]{4}
	\Vertex[NoLabel,x=0,y=2]{6}
	\Vertex[NoLabel,x=1,y=3]{7}
	\Vertex[NoLabel,x=2,y=4]{8}
	\Vertex[NoLabel,x=0,y=4]{11}
	\Vertex[NoLabel,x=1,y=5]{12}
	\Vertex[NoLabel,x=0,y=6]{13}
	\Edges[style={thick}](1,2)
	\Edges[style={thick}](2,3)
	\Edges[style={thick}](3,4)
	\Edges[style={thick}](2,6)
	\Edges[style={thick}](3,7)
	\Edges[style={thick}](4,8)
	\Edges[style={thick}](6,7)
	\Edges[style={thick}](7,8)
	\Edges[style={thick}](7,11)
	\Edges[style={thick}](8,12)
	\Edges[style={thick}](11,12)
	\Edges[style={thick}](12,13)
	\rotatebox{45}{\draw [decorate,decoration={brace,mirror,amplitude=10pt},yshift=-7pt] (-0.2,0) -- (4.3,0) node [black,midway,yshift=-0.6cm] {\footnotesize $n$};}
	\node at (2,-1) {$\talltriangle_n$};
\end{tikzpicture} \vline \quad \begin{tikzpicture}[scale=1.2]
	\SetFancyGraph
	\Vertex[NoLabel,x=0,y=0]{1}
	\Vertex[NoLabel,x=1,y=0]{2}
	\Vertex[NoLabel,x=2,y=0]{3}
	\Vertex[NoLabel,x=3,y=0]{4}
	\Vertex[NoLabel,x=0.5,y=0.5]{5}
	\Vertex[NoLabel,x=1.5,y=0.5]{6}
	\Vertex[NoLabel,x=2.5,y=0.5]{7}
	\Vertex[NoLabel,x=1,y=1]{8}
	\Vertex[NoLabel,x=2,y=1]{9}
	\Vertex[NoLabel,x=1.5,y=1.5]{10}
	\Edges[style={thick}](1,5)
	\Edges[style={thick}](2,5)
	\Edges[style={thick}](2,6)
	\Edges[style={thick}](3,6)
	\Edges[style={thick}](3,7)
	\Edges[style={thick}](4,7)
	\Edges[style={thick}](5,8)
	\Edges[style={thick}](6,8)
	\Edges[style={thick}](6,9)
	\Edges[style={thick}](7,9)
	\Edges[style={thick}](8,10)
	\Edges[style={thick}](9,10)
	\draw [decorate,decoration={brace,amplitude=10pt,mirror},yshift=-2pt] (-0.1,0) -- (3.1,0) node [black,midway,yshift=-0.6cm] {\footnotesize $n$};
	\node at (1.5,-1) {$\widetriangle_n$};
\end{tikzpicture} \quad \vline \quad \begin{tikzpicture}[scale=1.2]
	\SetFancyGraph
	\Vertex[NoLabel,x=0,y=0]{1}
	\Vertex[NoLabel,x=1,y=0]{2}
	\Vertex[NoLabel,x=2,y=0]{3}
	\Vertex[NoLabel,x=3,y=0]{4}
	\Vertex[NoLabel,x=0.5,y=0.5]{5}
	\Vertex[NoLabel,x=1.5,y=0.5]{6}
	\Vertex[NoLabel,x=2.5,y=0.5]{7}
	\Vertex[NoLabel,x=1,y=1]{8}
	\Vertex[NoLabel,x=2,y=1]{9}
	\Vertex[NoLabel,x=1.5,y=1.5]{10}
	\Vertex[NoLabel,x=0,y=1]{11}
	\Vertex[NoLabel,x=0.5,y=1.5]{12}
	\Vertex[NoLabel,x=0,y=2]{13}
	\Vertex[NoLabel,x=1,y=2]{14}
	\Vertex[NoLabel,x=0.5,y=2.5]{15}
	\Vertex[NoLabel,x=0,y=3]{16}
	\Edges[style={thick}](1,5)
	\Edges[style={thick}](2,5)
	\Edges[style={thick}](2,6)
	\Edges[style={thick}](3,6)
	\Edges[style={thick}](3,7)
	\Edges[style={thick}](4,7)
	\Edges[style={thick}](5,8)
	\Edges[style={thick}](6,8)
	\Edges[style={thick}](6,9)
	\Edges[style={thick}](7,9)
	\Edges[style={thick}](8,10)
	\Edges[style={thick}](9,10)
	\Edges[style={thick}](5,11)
	\Edges[style={thick}](11,12)
	\Edges[style={thick}](8,12)
	\Edges[style={thick}](12,13)
	\Edges[style={thick}](12,14)
	\Edges[style={thick}](10,14)
	\Edges[style={thick}](13,15)
	\Edges[style={thick}](14,15)
	\Edges[style={thick}](15,16)
	\draw [decorate,decoration={brace,amplitude=10pt,mirror},yshift=-2pt] (-0.1,0) -- (3.1,0) node [black,midway,yshift=-0.6cm] {\footnotesize $n$};
	\node at (1.5,-1) {$\cornertriangle_n$};
\end{tikzpicture}
\caption{The triangle posets.} \label{fig:tri_posets}
\end{figure}
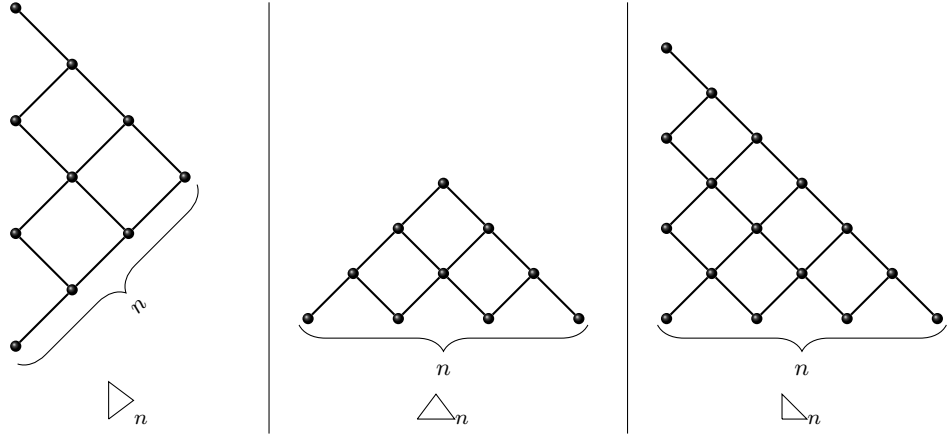

The other posets we care about in this section are three families of triangular posets which we denote $\talltriangle_n$, $\widetriangle_n$, and $\cornertriangle_n$. The Hasse diagrams of these triangular posets are depicted in Fig.~\ref{fig:tri_posets}. These three families of triangular posets are, in addition to the rectangle poset, the major examples of the nice posets for which we conjectured cyclic sieving under rowmotion in~\cite{hopkins2019minuscule}. Namely:

\begin{Conjecture}[{special case of~\cite[Conjecture~4.23]{hopkins2019minuscule}}] \label{conj:tall_csiev}
Fix $n,m \geq 1$ and let
\[F(q) \coloneqq \prod_{1\leq i \leq j \leq n}\frac{\big(1-q^{i+j+m-1}\big)}{\big(1-q^{i+j-1}\big)}.\]
Then for all $k \in \mathbb{Z}$ we have
\[ \#\big\{\pi \in \mathrm{PP}^{m}(\talltriangle_n)\colon \mathrm{Row}^{k}(\pi) = \pi\big\} = F\big(q \mapsto \zeta^{k}\big),\]
where $\zeta \coloneqq {\rm e}^{\pi {\rm i} / n}$ is a primitive $(2n)$th root of unity.
\end{Conjecture}

\begin{Conjecture}[{special case of~\cite[Conjecture 4.28]{hopkins2019minuscule}}] \label{conj:wide_csiev}
Fix $n,m \geq 1$ and let
\[F(q) \coloneqq \prod_{1\leq i \leq j \leq n}\frac{\big(1-q^{i+j+2m}\big)}{\big(1-q^{i+j}\big)}.\]
Then for all $k \in \mathbb{Z}$ we have
\[ \#\big\{\pi \in \mathrm{PP}^{m}(\widetriangle_n)\colon \mathrm{Row}^{k}(\pi) = \pi\big\} = F\big(q \mapsto \zeta^{k}\big),\]
where $\zeta \coloneqq {\rm e}^{\pi {\rm i} /(n+1)}$ is a primitive $(2(n+1))$th root of unity.
\end{Conjecture}

\begin{Conjecture}[{special case of~\cite[Conjecture 4.28]{hopkins2019minuscule}}] \label{conj:corner_csiev}
Fix $n,m \geq 1$ and let
\[F(q) \coloneqq \prod_{1 \leq i,j \leq n} \frac{\big(1-q^{2(i+j+m-1)}\big)}{\big(1-q^{2(i+j-1)}\big)}.\]
Then for all $k \in \mathbb{Z}$ we have
\[ \#\big\{\pi \in \mathrm{PP}^{m}(\cornertriangle_n)\colon \mathrm{Row}^{k}(\pi) = \pi\big\} = F\big(q \mapsto \zeta^{k}\big),\]
where $\zeta \coloneqq {\rm e}^{\pi {\rm i} /2n}$ is a primitive $(4n)$th root of unity.
\end{Conjecture}

\begin{Remark}
The $F(q)$ appearing in Conjectures~\ref{conj:tall_csiev}, \ref{conj:wide_csiev} and~\ref{conj:corner_csiev} are all actually polynomials with nonnegative integer coefficients~$F(q)\in \mathbb{N}[q]$. For example, the $F(q)$ in Conjecture~\ref{conj:tall_csiev} is the size generating function for $P$-partitions in $ \mathrm{PP}^{m}(\talltriangle_n)$: that is, $F(q) = \sum_{\pi \in \mathrm{PP}^{m}(P)} q^{|\pi|}$, where $|\pi| \coloneqq \sum_{p \in P} \pi(p)$. In fact, this is the same as $\mathrm{SymMac}'(n,m;q)$ from Section~\ref{sec:intro}. Meanwhile, the~$F(q)$ in Conjecture~\ref{conj:wide_csiev} is $q^{m \cdot \binom{n+1}{2}}$ times the quantity denoted ``(CGI)'' by Proctor in~\cite{proctor1990new}. Proctor explained how this expression is the generating function for $P$-partitions in $\mathrm{PP}^{m}(\widetriangle_n)$ with respect to a certain statistic (the statistic in question is slightly more complicated than size- it involves an alternating sum of entries). Finally, the $F(q)$ in Conjecture~\ref{conj:corner_csiev} is the result of applying the substitution $q\mapsto q^2$ to $\mathrm{Mac}(n,n,m;q)$. For all these conjectures, the case $k=0$ is known; that is, $F(1)$ is known to be equal to the total number of height~$m$ plane partitions of the poset. Also, for all these conjectures, the case $m=1$ is known~\cite{armstrong2013uniform, rush2013orbits}.
\end{Remark}

\begin{Remark}
Conjecrure~\ref{conj:tall_csiev} is part of a more general conjecture concerning all the \emph{minuscule posets}. The rectangle poset is the most prominent example of a minuscule poset, and the poset~$\talltriangle_n$ is the next most prominent example (in this context, it is commonly referred to as the ``\emph{shifted staircase}''). Besides the rectangle and the shifted staircase, the only other minuscule posets are one very simple infinite family (the ``\emph{propeller poset}'') and two exceptional posets (corresponding to $E_6$ and $E_7$). Hence, together with Theorem~\ref{thm:ppart_csiev}, establishing Conjecture~\ref{conj:tall_csiev} would resolve ``most'' of~\cite[Conjecture~4.23]{hopkins2019minuscule}.

Similarly, Conjectures~\ref{conj:wide_csiev} and~\ref{conj:corner_csiev} are part of a more general conjecture concerning the \emph{root posets of coincidental type}. The posets $\widetriangle_n$ (which is~$\Phi^{+}(A_n)$) and $\cornertriangle_n$ (which is~$\Phi^{+}(B_n)$) are the most prominent examples of root posets of coincidental type. Beyond these, there is only one very simple infinite family ($\Phi^{+}(I_2(\ell))$), and one exceptional poset ($\Phi^{+}(H_3)$). Hence, establishing Conjectures~\ref{conj:wide_csiev} and~\ref{conj:corner_csiev} would resolve ``most'' of~\cite[Conjecture~4.28]{hopkins2019minuscule}.
\end{Remark}

\begin{Remark}
In~\cite{fontaine2014cyclic}, Fontaine and Kamnitzer use some ideas from geometric representation theory to obtain a kind of refinement of Rhoades's~Theorem~\ref{thm:ppart_csiev} in which, rather than considering the action of promotion on the whole set $\mathrm{SSYT}(m^a,a+b)$, they consider only those tableaux with a fixed (cyclically symmetric) content. The relevant sieving polynomial turns out to be the corresponding \emph{Kostka--Foulkes polynomial} (in fact, in~\cite{rhoades2010cyclic} Rhoades obtained a less precise version of this result in which only the absolute value of the root of unity evaluation was considered). It is possible that Conjectures~\ref{conj:tall_csiev},~\ref{conj:wide_csiev} and~\ref{conj:corner_csiev} have similar ``content'' refinements as well. Actually, this possibility is discussed in~\cite[Remark~4.25]{hopkins2019minuscule} where it is suggested that the appropriate \emph{Lusztig's $q$-analog of weight multiplicity} could be the sieving polynomial. We will not discuss ``content'' refinements further here.
\end{Remark}

Now we explain, following Grinberg--Roby~\cite{grinberg2015birational2}, how the $P$-partitions in $\mathrm{PP}^{m}(P)$ for the triangular posets $P$ are in $\mathrm{Row}$-equivariant bijection with the plane partitions in~$\mathrm{PP}^{m}(n\times n)$ which satisfy certain symmetry properties. This allows us to reformulate Conjectures~\ref{conj:tall_csiev}, \ref{conj:wide_csiev} and~\ref{conj:corner_csiev} as assertions about the number of plane partitions in~$\mathrm{PP}^{m}(n\times n)$ fixed by various subgroups of~$\langle \mathrm{Row},\mathrm{Tp} \rangle$.

First let us explain how to relate $\talltriangle_n$ to the rectangle, which is very easy.

\begin{Lemma} \label{lem:tall_tri_embed}
There is a $\mathrm{Row}$-equivariant bijection between $\mathrm{PP}^{m}(\talltriangle_n)$ and the subset of those $\pi \in \mathrm{PP}^{m}(n\times n)$ for which $\mathrm{Tp}(\pi) = \pi$.
\end{Lemma}
\begin{proof}
This is basically proved by Grinberg--Roby in~\cite[Section~9]{grinberg2015birational2}. They were working at the birational level; the result we want, at the piecewise-linear level, is even simpler than what they did.

Let us view a $P$-partition $\pi \in \mathrm{PP}^{m}(\talltriangle_n)$ as triangular array $\pi = (\pi_{i,j})_{1\leq i \leq j \leq n}$ of nonnegative integers $\pi_{i,j}\in \mathbb{N}$ such that:
\begin{itemize}\itemsep=0pt
\item $\pi$ is weakly decreasing in rows and columns (i.e., $\pi_{i,j} \geq \pi_{i+1,j}$, $\pi_{i,j} \geq \pi_{i,j+1}$ for all $i$,~$j$),
\item the maximal entry of $\pi$ is less than or equal to~$m$ (i.e., $\pi_{1,1} \leq m$).
\end{itemize}
From such a $\pi$ we obtain a $\pi' \in \mathrm{PP}^{m}(n \times n)$ by setting
\[ \pi'_{i,j} = \begin{cases} \pi_{i,j} &\textrm{if $i\leq j$}, \\ \pi_{j,i} &\textrm{if $i > j$}.\end{cases}\]
This map $\pi \mapsto \pi'$ gives the desired bijection. In particular, it is easily seen to be rowmotion equivariant.
\end{proof}

In order to relate rowmotion of $\widetriangle_n$ to rowmotion of the rectangle, we have to review a remarkable ``reciprocity'' property of rowmotion acting on the rectangle that was established by Grinberg--Roby~\cite{grinberg2015birational2}.

\begin{Theorem}[{\cite[Theorem 32]{grinberg2015birational2}}] \label{thm:reciprocity}
For any $\pi \in \mathrm{PP}^{m}(a+b)$ we have
\[ \pi_{a+1-i,b+1-j} + \mathrm{Row}^{i+j-1}(\pi)_{i,j} = m \]
for all $1\leq i \leq a$ and $1 \leq j \leq b$.
\end{Theorem}

Actually, Grinberg--Roby proved the birational version of Theorem~\ref{thm:reciprocity}; but the birational version implies the piecewise-linear version we have stated via tropicalization. Theorem~\ref{thm:reciprocity} allows us to relate $\widetriangle_n$ to the rectangle, as follows.

\begin{Lemma} \label{lem:wide_tri_embed}
There is a $\mathrm{Row}$-equivariant bijection between $\mathrm{PP}^{M}(\widetriangle_{n-1})$ and the subset of those $\pi \in \mathrm{PP}^{2M}(n\times n)$ for which $\mathrm{Tp} \cdot \mathrm{Row}^{n}(\pi) = \pi$.
\end{Lemma}
\begin{proof}
This is basically proved by Grinberg--Roby in~\cite[Section~10]{grinberg2015birational2} (they worked at the birational level, but via tropicalization their results imply the corresponding piecewise-linear statements). They explained how to embed $\mathrm{PP}^{M}(\widetriangle_{n-1})$ into $\mathrm{PP}^{2M}(n\times n)$ in a $\mathrm{Row}$-equivariant way. We now review their embedding.

Let us view a $P$-partition $\pi \in \mathrm{PP}^{M}(\widetriangle_{n-1})$ as triangular array $\pi = (\pi_{i,j})_{\substack{1\leq i, j \leq n-1\\ i+j\leq n}}$ of nonnegative integers $\pi_{i,j}\in \mathbb{N}$ such that:
\begin{itemize}\itemsep=0pt
\item $\pi$ is weakly decreasing in rows and columns (i.e., $\pi_{i,j} \geq \pi_{i+1,j}$, $\pi_{i,j} \geq \pi_{i,j+1}$ for all $i$,~$j$),
\item the maximal entry of $\pi$ is less than or equal to~$m$ (i.e., $\pi_{1,1} \leq m$).
\end{itemize}
From such a $\pi$ we obtain a $\pi' \in \mathrm{PP}^{2M}(n \times n)$ by setting
\[ \pi'_{i,j} = \begin{cases} \pi_{i,j}+M &\textrm{if $i+j \leq n$}, \\ M &\textrm{if $i+j=n+1$}, \\ M - \mathrm{Row}^{k}(\pi)_{i-k,j-k}, \textrm{ with $k=i+j-n-1$} &\textrm{if $i+j>n+1$}.\end{cases}\]
An example of the map $\pi \mapsto \pi'$ in the case $n=3$, $M=2$ is
\[ \begin{ytableau} 2 & 1 \\ 0 \end{ytableau} \mapsto \begin{ytableau} 4 & 3 & 2 \\ 2 & 2 & 2 \\ 2 & 1 & 0 \end{ytableau} \]
In~\cite[Lemma~67]{grinberg2015birational2} it is shown that the map $\pi \mapsto \pi'$ is indeed an embedding of~$\mathrm{PP}^{M}(\widetriangle_{n-1})$ into $\mathrm{PP}^{2M}(n\times n)$ which is equivariant under rowmotion.

Moreover,~\cite[Theorem 65]{grinberg2015birational2} implies that if $\pi \in \mathrm{PP}^{2M}(n \times n)$ is in the image of this embedding of $\mathrm{PP}^{M}(\widetriangle_{n-1})$ into $\mathrm{PP}^{2M}(n\times n)$ then we will have $ \mathrm{Row}^{n}(\pi) = \mathrm{Tp}(\pi)$. Indeed, that essentially follows from the reciprocity result, Theorem~\ref{thm:reciprocity}. It also can be shown using Theorem~\ref{thm:reciprocity} that the \emph{only} $\pi \in \mathrm{PP}^{2M}(n \times n)$ with $\mathrm{Row}^{n}(\pi) = \mathrm{Tp}(\pi)$ are in the image of this embedding. But in fact, we can also establish that there are no other such $\pi$ via a counting argument. Namely, Corollary~\ref{cor:tr_csiev} implies that
\[\#\big\{\pi \in \mathrm{PP}^{2M}(n\times n) \colon \mathrm{Tp} \cdot \mathrm{Row}^{n}(\pi) = \pi\big\} = \mathrm{SymMac}(n,2M;q \mapsto -1).\]
Recall from Section~\ref{sec:intro} that
\[\mathrm{SymMac}(n,2M;q \mapsto -1) = \prod_{1\leq i \leq j \leq n-1} \frac{i+j+2M}{i+j}.\]
And it is known that $\#\mathrm{PP}^{M}(\widetriangle_{n-1}) = \prod_{1\leq i \leq j \leq n-1} \frac{i+j+2M}{i+j}$; see, for instance,~\cite[Case~``(CG)'' of Theorem~1]{proctor1990new}. This completes the proof of the lemma.
\end{proof}

\begin{Remark}
If $m$ is odd, there are no $\pi \in \mathrm{PP}^{m}(n\times n)$ for which $\mathrm{Tp} \cdot \mathrm{Row}^{n}(\pi) = \pi$. This can be seen either from the Grinberg--Roby reciprocity result, Theorem~\ref{thm:reciprocity}, or from our Corollary~\ref{cor:tr_csiev}.
\end{Remark}

\begin{Remark}
Recall the set $\mathrm{CY}(n,2M)$ defined in Section~\ref{sec:pro_tr}: this can be thought of as the subset of $\mathrm{PP}^{2M}(\talltriangle_{n-1})$ with even entries along the leftmost ``file''. In Section~\ref{sec:pro_tr} we explained how $\mathrm{CY}(n,2M)$ is naturally in bijection with the set of plane partitions~$\pi \in \mathrm{PP}^{2M}(n\times n)$ for which $\mathrm{Tp} \cdot \mathrm{Pro}(\pi) = \pi$. Meanwhile, we just explained in the proof of Lemma~\ref{lem:wide_tri_embed} that $\mathrm{PP}^{M}(\widetriangle_{n-1})$ is naturally in bijection with the set of plane partitions~$\pi \in \mathrm{PP}^{2M}(n\times n)$ for which $\mathrm{Tp} \cdot \mathrm{Row}^{n}(\pi) = \pi$. Then observe that
\[\#\mathrm{CY}(n,2M) = \prod_{1\leq i \leq j \leq n-1} \frac{i+j+2M}{i+j} = \#\mathrm{PP}^{M}(\widetriangle_{n-1}),\]
 as was proved for instance in the paper of Proctor~\cite{proctor1990new}. However, it is not at all a priori clear that $\mathrm{CY}(n,2M)$ and $\mathrm{PP}^{M}(\widetriangle_{n-1})$ have the same size, and constructing a bijection between these two sets is rather difficult (a bijection appears for instance in~\cite{sheats1999symplectic}). Thus, by studying the way these two operators interact with transposition, we uncovered another interesting ``duality'' between promotion and rowmotion.
\end{Remark}

We can relate $\cornertriangle_{n}$ to the rectangle by combining the previous two ideas.

\begin{Lemma} \label{lem:corner_tri_embed}
There is a $\mathrm{Row}$-equivariant bijection between $\mathrm{PP}^{M}(\cornertriangle_{n})$ and the subset of those $\pi \in \mathrm{PP}^{2M}(2n\times 2n)$ for which $\mathrm{Tp}(\pi) = \pi$ and $\mathrm{Row}^{2n}(\pi) = \pi$.
\end{Lemma}
\begin{proof}
This follows from combining the ideas in the proofs of Lemmas~\ref{lem:tall_tri_embed} and~\ref{lem:wide_tri_embed}. There is an obvious ``transposition'' symmetry of $\widetriangle_{2n-1}$ which reflects the poset across the vertical line of symmetry. Let $\mathrm{Tp}\colon \mathrm{PP}^{M}(\widetriangle_{2n-1}) \to \mathrm{PP}^{M}(\widetriangle_{2n-1})$ denote the induced symmetry of $P$-partitions. Then, by the same argument as in the proof of Lemma~\ref{lem:tall_tri_embed}, $\mathrm{PP}^{M}(\cornertriangle_{n})$ is in $\mathrm{Row}$-equivariant bijection with those $\pi \in \mathrm{PP}^{M}(\widetriangle_{2n-1})$ with $\mathrm{Tp}(\pi)=\pi$. Via the proof of Lemma~\ref{lem:wide_tri_embed} we can further embed $\mathrm{PP}^{M}(\cornertriangle_{n})$ into $\mathrm{PP}^{2M}(2n\times 2n)$ in desired way.
\end{proof}

The preceding lemmas about how to embed the triangular posets into the rectangle allow us to reformulate Conjectures~\ref{conj:tall_csiev},~\ref{conj:wide_csiev} and~\ref{conj:corner_csiev}, as follows:

\begin{Conjecture}[reformulation of Conjecture~\ref{conj:tall_csiev} in light of Lemma~\ref{lem:tall_tri_embed}] \label{conj:tall_csiev_redux}
For any $k \in \mathbb{Z}$ we have that
\[ \#\big\{\pi \in \mathrm{PP}^{m}(n \times n)\colon \mathrm{Tp}(\pi) = \pi, \mathrm{Row}^{k}(\pi) = \pi\big\} = F\big(q \mapsto \zeta^{k}\big),\]
where $\zeta \coloneqq {\rm e}^{\pi {\rm i} /n}$ is a primitive $(2n)$th root of unity and
\[F(q) \coloneqq \prod_{1\leq i \leq j \leq n}\frac{\big(1-q^{i+j+m-1}\big)}{\big(1-q^{i+j-1}\big)}.\]
\end{Conjecture}

\begin{Conjecture}[reformulation of Conjecture~\ref{conj:wide_csiev} in light of Lemma~\ref{lem:wide_tri_embed}] \label{conj:wide_csiev_redux}
For any $k \in \mathbb{Z}$ we have that
\[ \#\big\{\pi \in \mathrm{PP}^{2M}(n \times n)\colon \mathrm{Tp}\cdot \mathrm{Row}^{n}(\pi) = \pi, \mathrm{Row}^{k}(\pi) = \pi\big\} = F\big(q \mapsto \zeta^{k}\big),\]
where $\zeta \coloneqq {\rm e}^{\pi {\rm i} /n}$ is a primitive $(2n)$th root of unity and
\[F(q) \coloneqq \prod_{1\leq i \leq j \leq n-1}\frac{\big(1-q^{i+j+2M}\big)}{\big(1-q^{i+j}\big)}.\]
\end{Conjecture}

\begin{Conjecture}[reformulation of Conjecture~\ref{conj:corner_csiev} in light of Lemma~\ref{lem:corner_tri_embed}] \label{conj:corner_csiev_redux}
For any $k \in \mathbb{Z}$ we have that
\[ \#\big\{\pi \in \mathrm{PP}^{2M}(2n \times 2n)\colon \mathrm{Tp}(\pi) = \pi, \mathrm{Row}^{2n}(\pi) = \pi, \mathrm{Row}^{k}(\pi) = \pi\big\} = F\big(q \mapsto \zeta^{k}\big),\]
where $\zeta \coloneqq {\rm e}^{\pi {\rm i} /2n}$ is a primitive $(4n)$th root of unity and
\[F(q) \coloneqq \prod_{1 \leq i,j \leq n} \frac{\big(1-q^{2(i+j+M-1)}\big)}{\big(1-q^{2(i+j-1)}\big)}.\]
\end{Conjecture}

\begin{Remark}Conjecture~\ref{conj:corner_csiev_redux} follows easily from Conjecture~\ref{conj:wide_csiev_redux} or Conjecture~\ref{conj:tall_csiev_redux} via the same kind of argument as in~\cite[proof of Proposition~2.4]{armstrong2013uniform}.
\end{Remark}

Corollaries~\ref{cor:ppart_csiev} and~\ref{cor:tr_csiev} from the previous section say that for any $g \in \langle \mathrm{Row}, \mathrm{Tp} \rangle$, the number of plane partitions in $\mathrm{PP}^{m}(n\times n)$ fixed by $g$ is given by some kind of sieving phenomenon evaluation of a polynomial with a simple product formula as a rational expression. In other words, for any cyclic subgroup $H \subseteq \langle \mathrm{Row}, \mathrm{Tp} \rangle$, the number of plane partitions in $\mathrm{PP}^{m}(n\times n)$ fixed by $H$ is given by such an evaluation. Meanwhile, Conjectures~\ref{conj:tall_csiev_redux},~\ref{conj:wide_csiev_redux} and~\ref{conj:corner_csiev_redux} assert that for various \emph{noncyclic} subgroups $H \subseteq \langle \mathrm{Row}, \mathrm{Tp} \rangle$, the number of plane partitions in $\mathrm{PP}^{m}(n\times n)$ fixed by $H$ is given by such an evaluation. This leads us to wonder the following:

\begin{Question} \label{question:main}
Is the number of plane partitions in $\mathrm{PP}^{m}(n\times n)$ fixed by $H$, where~$H$ is any subgroup of $\langle \mathrm{Row}, \mathrm{Tp} \rangle$, given by a sieving phenomenon-type evaluation at a~root of unity of a~polynomial with a simple product formula as a rational expression?
\end{Question}

If Question~\ref{question:main} had a positive answer, it would be very similar to what happens with the ``classical'' symmetries of plane partitions, where the 10 symmetry classes all have product formulas for their enumeration (again see~\cite{krattenthaler2016plane, stanley1986symmetries}).

\begin{Remark}The polynomial $F(q)$ appearing in Conjecture~\ref{conj:tall_csiev} is the principal specialization of the character of the irreducible ${\rm SO}(2n+1)$ representation $V(m \omega_n)$, where~$\omega_n$ is the minuscule weight of type~$B_n$. (Technically one might have to work with the simply connected double cover $\widetilde{{\rm SO}}(2n+1)$, i.e., the so-called ``\emph{spin group}''.) Geometrically, this representation is the dual of the $m$th homogeneous component of the coordinate ring of the maximal orthogonal Grassmannian~$\mathrm{OG}(n,2n+1)$. Similarly, the polynomial $F(q)$ appearing in Conjecture~\ref{conj:wide_csiev} is (essentially) the principal specialization of the character of the irreducible ${\rm Sp}(2n)$ representation $V(m \omega_n)$, where~$\omega_n$ is the cominuscule weight of type~$C_n$. Geometrically, this representation is the dual of the $m$th homogeneous component of the coordinate ring of the Lagrangian Grassmannian~$\mathrm{LG}(n,2n)$. See Proctor~\cite{proctor1984bruhat, proctor1990new} or Stembridge~\cite{stembridge1994minuscule} for more information about these polynomials. At any rate, the fact that these polynomials are more-or-less Lie group characters naturally suggests an approach for resolving Conjectures~\ref{conj:tall_csiev} and~\ref{conj:wide_csiev}: find a basis of the corresponding representation indexed by the set of $P$-partitions in question and such that an appropriate group element (e.g., the lift of a Coxeter element) permutes the basis according to rowmotion (or, more likely, according to a conjugate ``promotion''-like action). In other words, extend Rhoades's approach~\cite{rhoades2010cyclic} to other types. The problem with this approach is that the naive bases like the standard monomial basis fail to behave in the appropriate way, while the sophisticated bases like the dual canonical basis or the theta basis are extremely hard to concretely get one's hands on, and doubly so outside of type~A. (We do not mean to suggest that it is \emph{totally} hopeless to work in other types. For instance, the theory of plabic graphs~\cite{postnikov2006total, scott2006grassmannians} is ultimately the combinatorial underpinning of the Shen--Weng~\cite{shen2018cyclic} proof of cyclic sieving, and the work of Karpman~\cite{karpman2018total, karpman2019purity, karpman2018symmetric} extends much of the theory of plabic graphs to the Lagrangian Grassmannian.) Our results in this paper point the way to an alternative but ultimately complementary approach to Conjectures~\ref{conj:tall_csiev} and~\ref{conj:wide_csiev}: stay in the type~A world but impose symmetries.
\end{Remark}

\appendix

\section{Plane partitions and semistandard tableaux} \label{sec:appendix}

In this appendix we explain the correspondence between plane partitions and semistandard tableaux of rectangular shape, and how the operators we are interested in (e.g., promotion) behave under this correspondence.

A \emph{partition} $\lambda = (\lambda_ 1 \geq \lambda_2 \geq \cdots)$ is an infinite nonincreasing sequence of integers for which $\lambda_i=0$ for all~$i \gg 0$. The nonzero $\lambda_i$ are called the \emph{parts} of $\lambda$. If we write a partition as~$\lambda=(\lambda_1,\dots,\lambda_k)$ that means that $\lambda_i=0$ for $i>k$. A particularly important family of partitions for us will be the \emph{rectangle} partitions $m^{a} \coloneqq (\overbrace{m,m,\dots,m}^{a})$. We represent a partition~$\lambda$ via its \emph{Young diagram}, which is the collection of boxes in rows with $\lambda_i$ boxes left-justified in row~$i$. For example, the Young diagram of $(4,3,1,1)$ is
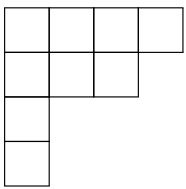
\[\ydiagram{4,3,1,1}\]

A \emph{semistandard Young tableau} of shape $\lambda$ is a filling of the boxes of the Young diagram of $\lambda$ with positive integers that is weakly increasing in rows and strictly increasing in columns. For example, the following is a semistandard Young tableau of shape $(4,3,1,1)$:
\[\begin{ytableau}1 & 3 & 3 & 6 \\ 2 & 5 & 5 \\ 4 \\ 5 \end{ytableau}\]
We use $\mathrm{SSYT}(\lambda,k)$ to denote the set of semistandard Young tableaux of shape $\lambda$ whose entries belong to $\{1,2,\dots,k\}$. Note that this set is empty if $k$ is less than the number of parts of $\lambda$.

Now we define the \emph{promotion} operator $\rho$ acting on~$\mathrm{SSYT}(\lambda,k)$.\footnote{There is no consensus in the literature about whether the operator $\rho$ or $\rho^{-1}$ is the one which should be called ``promotion''. We follow the convention of Rhoades~\cite{rhoades2010cyclic} and Lam~\cite{lam2019cyclic}; but for instance Bloom, Pechenik, and Saracino~\cite{bloom2016homomesy} take the opposite convention. Of course, ultimately these differences in convention are inconsequential.} Roughly speaking, promotion behaves as follows: first we delete all entries of $k$ in our tableau, leaving holes in their places; then we slide the remaining entries into the holes so that the holes occupy the upper-left of the Young diagram; then we increment by one all entries in the tableau; and finally we fill the holes in the upper-left with $1$'s. For example, with $k=6$, an application of promotion might look like the following:
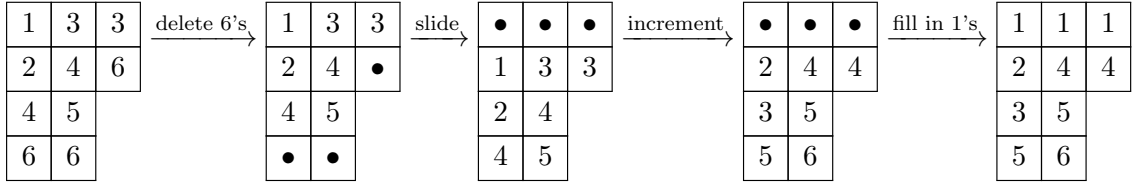
\[ \begin{ytableau}1 & 3 & 3 \\ 2 & 4 & 6 \\ 4 & 5 \\ 6 & 6 \end{ytableau} \xrightarrow{\textrm{delete $6$'s}} \begin{ytableau}1 & 3 & 3 \\ 2 & 4 & \bullet \\ 4 & 5 \\ \bullet & \bullet \end{ytableau} \xrightarrow{\textrm{slide}} \begin{ytableau}\bullet & \bullet & \bullet \\ 1 & 3 & 3 \\ 2 & 4 \\ 4 & 5 \end{ytableau} \xrightarrow{\textrm{increment}} \begin{ytableau}\bullet & \bullet & \bullet \\ 2 & 4 & 4 \\ 3 & 5 \\ 5 & 6 \end{ytableau} \xrightarrow{\textrm{fill in $1$'s}} \begin{ytableau} 1 & 1 & 1 \\ 2 & 4 & 4 \\ 3 & 5 \\ 5 & 6 \end{ytableau} \]
To formalize this definition would require more explanation of the sliding procedure. A precise description is given in~\cite[Section~2]{rhoades2010cyclic}. In fact, these slides are the ``jeu de taquin'' moves of Sch\"{u}tzenberger~\cite{schuztenberger1963quelques, schutzenberger1972promotion, schutzenberger1977robinson} (see also the presentation of Haiman~\cite{haiman1992dual}).

Another closely related operator acting on~$\mathrm{SSYT}(\lambda,k)$ is \emph{evacuation} $\varepsilon$. It can also be defined in terms of jeu de taquin slides. Evacuation roughly behaves as follows: first we rotate the tableau $180^{\circ}$; then we replace every entry by $k+1$ minus that entry; finally, we slide the entries into the upper-left so that we get back to a Young diagram shape. For example, with $k=6$, an application of evacuation might look like the following:
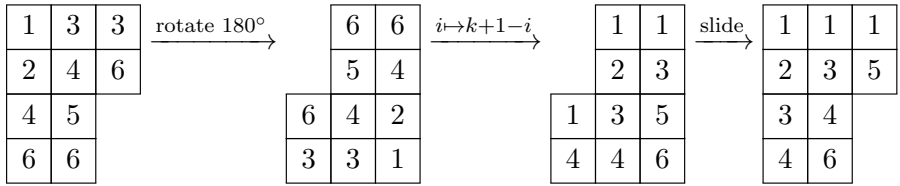
\[ \begin{ytableau}1 & 3 & 3 \\ 2 & 4 & 6 \\ 4 & 5 \\ 6 & 6 \end{ytableau} \xrightarrow{\textrm{rotate $180^\circ$}} \begin{ytableau} \none & 6 & 6 \\ \none & 5 & 4 \\ 6 & 4 & 2 \\ 3 & 3 & 1 \end{ytableau} \xrightarrow{i \mapsto k+1-i} \begin{ytableau} \none & 1 & 1 \\ \none & 2 & 3 \\ 1 & 3 & 5 \\ 4 & 4 & 6 \end{ytableau} \xrightarrow{\textrm{slide}} \begin{ytableau}1 & 1 & 1 \\ 2 & 3 & 5 \\ 3 & 4 \\ 4 & 6 \end{ytableau} \]
Again, to formalize this definition of evacuation we would need to explain the jeu de taquin slides in more detail; a precise description is given in~\cite[Section~2]{rhoades2010cyclic}.

However, rather than use jeu de taquin, we will instead work with different but equivalent definitions of promotion and evacuation in terms of the so-called ``Bender--Knuth involutions''~\cite{bender1972enumeration}.

\begin{Definition}
The $i$th \emph{Bender--Knuth involution}, denoted $\mathcal{BK}_i\colon \mathrm{SSYT}(\lambda,k)\to \mathrm{SSYT}(\lambda,k)$, for $1 \leq i < k$, is the operator which acts on a tableau $T \in \mathrm{SSYT}(\lambda,k)$ as follows: first we ``freeze'' in place all $i$'s directly above $(i+1)$'s, and all $(i+1)$'s directly below $i$'s; and then, in each row, we change unfrozen $i$'s into $(i+1)$'s and unfrozen $(i+1)$'s into $i$'s in the unique way which preserves the semistandardness condition and so that the number of unfrozen $i$'s in that row in the resulting tableau is the number of unfrozen $(i+1)$'s in that row in the original tableau, and vice-versa. In other words, considering just the unfrozen $i$'s and $(i+1)$'s in a row, we perform the transformation $i^x(i+1)^y \mapsto i^y(i+1)^x$ on these entries.
\end{Definition}

\begin{Example}
Let $T\in \mathrm{SSYT}((11,10,8,4,2), 6)$ be the following tableau:
\[ T = \begin{ytableau} 1 & 1 & 1 & 1 & 2 & 2 & 2 & 3 & 3 & 3 & 4 \\ 2 & 2 & 3 & 3 & 3 & 3 & 4 & 4 & 4 & 5 \\ 3 & 4 & 4 & 4 & 5 & 5 & 5 & 5 \\ 5 & 5 & 5 & 6 \\ 6 & 6 \end{ytableau} \]
Suppose we want to compute $\mathcal{BK}_4(T)$. First of all, we can restrict our attention to only the $4$'s and $5$'s in the tableau, resulting in a picture that looks like this:
\[ \begin{ytableau} \, & \, & \, & \, & \, & \, & \, & \, & \, & \, & 4 \\ \, & \, & \, & \, & \, & \, & 4 & 4 & 4 & 5 \\ \, & 4 & 4 & 4 & 5 & 5 & 5 & 5 \\ 5 & 5 & 5 & \, \\ \, & \, \end{ytableau} \]
Then we ``freeze'' in place $4$'s directly above $5$'s and $5$'s directly below $4$'s. In the picture below these frozen boxes have been shaded:
\[ \begin{ytableau} \, & \, & \, & \, & \, & \, & \, & \, & \, & \, & 4 \\ \, & \, & \, & \, & \, & \, & *(lightgray) 4 & *(lightgray) 4 & 4 & 5 \\ \, & *(lightgray) 4 & *(lightgray) 4 & 4 & 5 & 5 & *(lightgray) 5 & *(lightgray) 5 \\ 5 & *(lightgray) 5 & *(lightgray) 5 & \, \\ \, & \, \end{ytableau} \]
Then within each row we swap the number of unfrozen $4$'s and $5$'s, while preserving the weakly decreasing requirement, resulting in the following:
\[ \begin{ytableau} \, & \, & \, & \, & \, & \, & \, & \, & \, & \, & 5 \\ \, & \, & \, & \, & \, & \, & *(lightgray) 4 & *(lightgray) 4 & 4 & 5 \\ \, & *(lightgray) 4 & *(lightgray) 4 & 4 & 4 & 5 & *(lightgray) 5 & *(lightgray) 5 \\ 4 & *(lightgray) 5 & *(lightgray) 5 & \, \\ \, & \, \end{ytableau} \]
Finally, to obtain $\mathcal{BK}_4(T)$, we put the entries which are not $4$'s or $5$'s back in as they were in~$T$:
\[ \mathcal{BK}_4(T)= \begin{ytableau} 1 & 1 & 1 & 1 & 2 & 2 & 2 & 3 & 3 & 3 & 5 \\ 2 & 2 & 3 & 3 & 3 & 3 & 4 & 4 & 4 & 5 \\ 3 & 4 & 4 & 4 & 4 & 5 & 5 & 5 \\ 4 & 5 & 5 & 6 \\ 6 & 6 \end{ytableau} \]
\end{Example}

It is clear that the $\mathcal{BK}_i$ are indeed involutions; however, unlike the reflection operators $s_i$ which act on $\mathrm{SSYT}(\lambda,k)$ thanks to its crystal structure, note that the $\mathcal{BK}_i$ \emph{do not} satisfy the braid relations, and hence do not give an action of the symmetric group on these tableaux.

We can define promotion and evacuation as a composition of the Bender--Knuth involutions.

\begin{Definition} \label{def:pro_evac}
\emph{Promotion} $\rho\colon \mathrm{SSYT}(\lambda,k)\to \mathrm{SSYT}(\lambda,k)$ is the following composition of the Bender--Knuth involutions
\[\rho \coloneqq \mathcal{BK}_{1} \cdot \mathcal{BK}_2 \cdots \mathcal{BK}_{k-2} \cdot \mathcal{BK}_{k-1}.\]
\emph{Evacuation} $\varepsilon\colon \mathrm{SSYT}(\lambda,k)\to \mathrm{SSYT}(\lambda,k)$ is
\[\varepsilon \coloneqq (\mathcal{BK}_1) \cdot (\mathcal{BK}_2 \cdot \mathcal{BK}_1) \cdots (\mathcal{BK}_{k-2} \cdots \mathcal{BK}_2 \cdot \mathcal{BK}_1) \cdot (\mathcal{BK}_{k-1} \cdots \mathcal{BK}_2 \cdot \mathcal{BK}_{1}). \]
\emph{Dual evacuation} $\varepsilon^{*}\colon \mathrm{SSYT}(\lambda,k)\to \mathrm{SSYT}(\lambda,k)$ is
\[\varepsilon^{*} \coloneqq (\mathcal{BK}_{k-1}) \cdot (\mathcal{BK}_{k-2} \cdot \mathcal{BK}_{k-1}) \cdots (\mathcal{BK}_{2} \cdots \mathcal{BK}_{k-2} \cdot \mathcal{BK}_{k-1}) \cdot (\mathcal{BK}_{1} \cdots \mathcal{BK}_{k-2} \cdot \mathcal{BK}_{k-1}). \]
\end{Definition}

We did not discuss dual evacuation earlier but it turns out to be useful in understanding the behavior of promotion and evacuation. It is a theorem of Gansner~\cite{gansner1980equality} (see also~\cite{bloom2016homomesy, chmutov2020berenstein}) that the definitions of promotion and evacuation in terms of the Bender--Knuth involutions are the same as those in terms of Sch\"{u}tzenberger's jeu de taquin moves.

As we will see in the next proposition, evacuation has a very simple behavior on tableaux of rectangular shape. In order to record that behavior we need a little notation. So for a~tableau $T\in \mathrm{SSYT}(m^a,k)$ we define $T^{+} \in \mathrm{SSYT}(m^a,k)$ to be the tableau obtained by rotating~$T$ $180^{\circ}$ and replacing every entry by $k+1$ minus that entry.

The following proposition records some basic properties of promotion and evacuation which are well known but are ``folklore''. The best reference we have for these results is the paper of Bloom--Pechenik--Saracino~\cite{bloom2016homomesy}, who adapt the arguments presented by Stanley~\cite{stanley2009promotion} in the case of standard Young tableaux, and also use the fundamental connection of promotion and evacuation to the \emph{Robinson--Schensted--Knuth correspondence}.

\begin{Proposition}[{see~\cite[Theorem 2.9]{bloom2016homomesy}, building off of~\cite[Theorem 2.1]{stanley2009promotion}}] \label{prop:ssyt_promo_basics}
For any shape $\lambda$ and any $k$ we have the following relationship among the operators $\rho,\varepsilon, \varepsilon^{*} \colon \mathrm{SSYT}(\lambda,k)\to \mathrm{SSYT}(\lambda,k)$:
\begin{itemize}\itemsep=0pt
\item $\varepsilon^2 = ( \varepsilon^{*})^{2} = \mathrm{id}$ $($the identity operator$)$,
\item $\rho \varepsilon = \varepsilon \rho^{-1}$,
\item $\rho^k = \varepsilon \varepsilon^{*}$.
\end{itemize}
Furthermore, if $\lambda = m^a$ is a rectangle then $\varepsilon(T) = \varepsilon^{*}(T) = T^+$ for all tableaux $T \in \mathrm{SSYT}(m^a,k)$. Consequently, $\rho^k = \mathrm{id}$ if $\lambda=m^a$ is a rectangle.
\end{Proposition}

Now we explain an equivalent, but very useful, way to think about semistandard Young tableaux: namely, as \emph{Gelfand--Tsetlin patterns}. In particular, Gelfand--Tsetlin patterns will serve as the bridge between semistandard tableaux and plane partitions. The usefulness of Gelfand--Tsetlin patterns for understanding operations on semistandard tableaux in terms of piecewise-linear expressions was especially emphasized in the papers of Berenstein and Kirillov~\cite{kirillov1995groups} and Berenstein and Zelevinsky~\cite{berenstein1996canonical}.

\begin{Definition}
Let $\lambda$ be a partition and $k$ an integer greater than or equal to the number of parts of $\lambda$. A \emph{Gelfand--Tsetlin pattern} of shape $\lambda$ and length $k$ is a triangular array $\pi = (\pi_{i,j})_{1\leq i \leq j \leq k}$ of nonnegative integers $\pi_{i,j}\in \mathbb{N}$ such that:
\begin{itemize}\itemsep=0pt
\item $\pi$ is weakly decreasing in rows and columns (i.e., $\pi_{i,j} \geq \pi_{i+1,j}$, $\pi_{i,j} \geq \pi_{i,j+1}$ for all $i$,~$j$),
\item the main diagonal $(\pi_{1,1},\pi_{2,2},\dots,\pi_{k,k})$ of $\pi$ is equal to the partition~$\lambda$.
\end{itemize}
We denote the set of such Gelfand--Tsetlin patterns by $\mathrm{GT}(\lambda,k)$.
\end{Definition}

There is a well-known bijection $\Phi\colon \mathrm{GT}(\lambda,k) \to \mathrm{SSYT}(\lambda,k)$: for $\pi \in \mathrm{GT}(\lambda,k)$, the tableau $T = \Phi(\pi)$ is the unique semistandard tableau such that for all $1\leq i \leq k$ the diagonal $(\pi_{i,1}, \pi_{i+1,2},\dots, \pi_{k,k+1-i})$ of $\pi$ is the shape of the restriction of $T$ to the entries $\{1,2,\dots,k+1-i\}$. To see that this is really a bijection, observe that $\pi_{k,1}$ is the number of $1$'s in $T$, and these must all go in the first row; similarly $\pi_{k-1,1}-\pi_{k,1}$ is the number of $2$'s in the first row and $\pi_{k-1,2}$ is the number of $2$'s in the second row; and so on. In this way we can clearly reconstruct a~unique tableau $T$ from $\pi$, and the inequalities imposed on the $\pi_{i,j}$ exactly correspond to the semistandardness condition.

\begin{Example}
Let $\lambda = (3,2,1,1,0)$ and let $\pi \in \mathrm{GT}(\lambda,5)$ be the following Gelfand--Tsetlin pattern:
\[\pi = \begin{matrix}
 3 & 3 & 3 & 1 & 1 \\
 & 2 & 1 & 1 & 0 \\
 & & 1 & 1 & 0 \\
 & & & 1 & 0 \\
 & & & & 0
\end{matrix}\]
Then $T=\Phi(\pi)\in\mathrm{SSYT}(\lambda,5)$ is the semistandard tableau
\[T=\begin{ytableau} 1 & 3 & 3 \\ 3 & 5 \\ 4 \\ 5 \end{ytableau} \]
\end{Example}

Naturally we want to understand how promotion behaves in terms of Gelfand--Tsetlin patterns. This is where the piecewise-linear toggles come in. We define the \emph{piecewise-linear toggle} $\tau_{i,j}\colon \mathrm{GT}(\lambda,k)\to \mathrm{GT}(\lambda,k)$ for $1 \leq i < j \leq k$ by
\[(\tau_{i,j}\pi)_{p,q} \coloneqq \begin{cases}\pi_{p,q} &\textrm{if $(p,q)\neq (i,j)$}, \\ \min (\pi_{i,j-1},\pi_{i-1,j}) + \max (\pi_{i+1,j},\pi_{i,j+1}) - \pi_{i,j} &\textrm{if $(p,q)=(i,j)$}, \\ \end{cases}\]
where we ignore $\pi_{i,j}$ with $i$,~$j$ outside of the bounds of the triangle (at least one term in each $\max $ and $\min $ will exist). Observe that these are exactly the same as the~$\tau_{i,j}$ defined in Section~\ref{sec:intro}. We again define $\mathcal{F}_{l} \coloneqq \prod_{\substack{1\leq i\leq j \leq k, \\ j-i=l}} \tau_{i,j}$ for~$1\leq l \leq k-1$ to be the composition of all the toggles along the ``$l$th diagonal'' of our array. As Berenstein--Kirillov~\cite{kirillov1995groups} explained, these diagonal toggles are the same as the Bender--Knuth involutions:\footnote{Hence one may view $\mathcal{BK}$ as honoring either Bender and Knuth or Berenstein and Kirillov.}

\begin{Proposition}[{\cite[Proposition~2.2]{kirillov1995groups}}] \label{prop:bk_toggles}
Viewing the Bender--Knuth involutions as operators on $\mathrm{GT}(\lambda,k)$ via the bijection $\Phi$, we have $\mathcal{BK}_{i} = \mathcal{F}_{k-i}$ for $1 \leq i \leq k-1$.
\end{Proposition}

So $\rho$ and $\varepsilon$ can be described in terms of piecewise-linear dynamics on $\mathrm{GT}(\lambda,k)$.

Finally, let us concentrate on the rectangular case and the correspondence with plane partitions. So suppose that $\lambda = m^a$, and let us take $k = a+b$ to match our indexing of plane partitions. Then consider what a Gelfand--Tsetlin pattern $\pi \in \mathrm{GT}(m^a,a+b)$ looks like. In the upper-left, $\pi$ has a length $a$ triangle of entries which must all be $m$'s; in the lower-right, $\pi$ has a length $b$ triangle of entries which must all be $0$'s; and the other entries in $\pi$, whose values are not forced, form an $a \times b$ rectangle. For example, with $a=3$, $b=2$, and $m=5$, we have
\[\pi = \begin{matrix}
 5 & 5 & 5 & * & * \\
 & 5 & 5 & * & * \\
 & & 5 & * & * \\
 & & & 0 & 0 \\
 & & & & 0
\end{matrix}\]
where the asterisks denote the entries whose values are not forced. What condition is placed on these asterisk entries? Well, they certainly must be weakly decreasing in rows and columns, and they must all be integers between $0$ and $m$. In other words, they exactly form a plane partition in $\mathrm{PP}^{m}(a\times b)$. And clearly any such plane partition can be placed in the asterisk entries.

In this way we obtain a bijection $\Psi\colon \mathrm{PP}^{m}(a\times b) \to \mathrm{SSYT}(m^a,a+b)$: we extend a plane partition $\pi \in \mathrm{PP}^{m}(a\times b)$ to a Gelfand--Tsetlin pattern in $\mathrm{GT}(m^a,a+b)$ by appending a length~$a$ triangle of $m$'s to its left and a length~$b$ triangle of $0$'s below it; and then we map that Gelfand--Tsetlin pattern to a semistandard tableau in $\mathrm{SSYT}(m^a,a+b)$ via~$\Phi$. This is depicted in the following example for a plane partition $\pi \in \mathrm{PP}^{4}(2\times 2)$:
\[\pi = \begin{ytableau} 2 & 2 \\ 1 & 0 \end{ytableau} \;\; \mapsto \;\; \begin{matrix}
 4 & 4 & 2 & 2 \\
 & 4 & 1 & 0 \\
 & & 0 & 0 \\
 & & & 0 \\
\end{matrix} \;\; \mapsto \;\; \begin{ytableau} 1 & 1& 3 & 3 \\ 3 & 4 & 4 & 4 \end{ytableau} = \Psi(\pi) \]
This construction is discussed, very briefly, in~\cite[pp.~516--517]{einstein2014piecewise}.

\begin{Remark}
There is a very naive way to obtain a semistandard Young tableau of rectangular shape from a plane partition $\pi \in \mathrm{PP}^{m}(a\times b)$: rotate $\pi$ $180^{\circ}$, and then add~$i$ to all entries in the~$i$th row. The bijection $\Psi$ is {\it not} this naive procedure. Indeed, this naive procedure produces a~tableau in $\mathrm{SSYT}(b^a, m+a)$, which is not the same set of tableaux that $\Psi$ maps into.
\end{Remark}

Let us describe another way to view the bijection $\Psi$, which is also useful.

We start with the case $m=1$. Note that a single column tableau $T \in \mathrm{SSYT}(1^a,a+b)$ is exactly the same as a subset $I\subseteq \{1,2,\dots,a+b\}$ of size $a$. So in this case $\Psi$ is some bijection $\Psi\colon \mathrm{PP}^{1}(a\times b) \xrightarrow{\sim} \{I\subseteq \{1,\dots,a+b\} \textrm{ of size $a$}\}$. In fact, this bijection is a correspondence between Young diagrams that fit in an $a\times b$ rectangle and size~$a$ subsets of $\{1,2,\dots,a+b\}$ which is ubiquitous in algebraic combinatorics, as we now explain. Let $\pi \in \mathrm{PP}^{1}(a\times b)$. The boundary between the entries of $1$ and $0$ in $\pi$ determines a lattice path of down and left steps from the upper-right corner of the $a\times b$ grid to the lower-left corner. Writing this lattice path as a word in the alphabet $\{D,L\}$ with $a$ $D$'s and $b$ $L$'s (where $D$'s correspond to down steps and $L$'s to left steps), the subset $\Psi(\pi)\subseteq \{1,2,\dots,a+b\}$ is the set of positions of $D$'s in this word. This is depicted in the following example with $a=4$ and $b=5$:
\[\begin{tikzpicture}
\node at (-2,0) {$\pi = $};
\node at (0,0) {\begin{ytableau} 1 & 1 & 1 & 1 & 0 \\ 1 & 1 & 0 & 0 & 0 \\ 1 & 1 & 0 & 0 & 0 \\ 0 & 0 & 0 & 0 & 0 \end{ytableau}};
\def\x{0.6}
\draw[line width=1mm] (-2.5*\x,-2*\x) -- (-2.5*\x,-1*\x) -- (-0.5*\x,-1*\x) -- (-0.5*\x,1*\x) -- (1.5*\x,1*\x) -- (1.5*\x,2*\x) -- (2.5*\x,2*\x);
\node at (4.5,0) {$\mapsto \textrm{L\underline{D}LL\underline{DD}LL\underline{D}} \mapsto \{2,5,6,9\} =$};
\node at (7.8,0) {\begin{ytableau} 2 \\ 5 \\ 6 \\ 9 \end{ytableau}};
\node at (8.9,0) {$= \Psi(\pi)$.};
\end{tikzpicture}\]

Now we extend the construction from the previous paragraph to greater values of~$m$. For $\pi, \pi' \in \mathrm{PP}^{m}(a\times b)$ we write $\pi \geq \pi'$ to mean that $\pi$ is entrywise greater than or equal to $\pi'$; and we define $\pi + \pi' \in \mathrm{PP}^{m+m'}(a\times b)$ for $\pi \in \mathrm{PP}^{m}(a\times b)$, $\pi' \in \mathrm{PP}^{m'}(a\times b)$ to be their entrywise sum. Let $\pi \in \mathrm{PP}^{m}(a\times b)$. There are unique plane partitions $\pi^{1}, \pi^{2},\dots, \pi^{m} \in \mathrm{PP}^{1}(a\times b)$ for which $\pi=\pi^{1} + \pi^{2} + \cdots + \pi^{m}$ and~$\pi^{1} \geq \pi^{2} \geq \cdots \geq \pi^{m}$; explicitly, we have
\[ (\pi^{k})_{i,j} =\begin{cases} 1 &\textrm{if $\pi_{i,j}\geq k$}, \\ 0 &\textrm{otherwise},\end{cases}\]
for all $1\leq k \leq m$. Then $\Psi(\pi) \in \mathrm{SSYT}(m^a,a+b)$ is the tableau whose columns are the subsets $\Psi(\pi^1), \Psi(\pi^2), \dots, \Psi(\pi^m)$ in order. This is easily proven inductively: the condition $\pi^{m-1} \geq \pi^{m}$ means that placing the column $\Psi(\pi^{m})$ to the right of the column $\Psi(\pi^{m-1})$ will preserve the tableau's semistandardness; conversely, appending the column $\Psi(\pi^{m})$ changes the entries in the Gelfand--Tsetlin pattern of the tableau in exactly the way that corresponds to adding~$\pi^{m}$. This description of $\Psi$ is depicted in the following example for a plane partition $\pi \in \mathrm{PP}^{4}(2\times 2)$:
\[\begin{tikzpicture}
\node at (-1.25,0) {$\pi =$};
\node at (0,0) {\begin{ytableau} 2 & 2 \\ 1 & 0 \end{ytableau}};
\node at (1,0) {$=$};
\node at (2,0) {\begin{tikzpicture}
\node at (0,0) {\begin{ytableau} 1 & 1 \\ 1 & 0 \end{ytableau}};
\def\x{0.6}
\draw[line width=1mm] (-1*\x,-1*\x) -- (0*\x,-1*\x) -- (0*\x,0*\x) -- (0*\x,0*\x) -- (1*\x,0*\x) -- (1*\x,1*\x);
\end{tikzpicture}};
\node at (3,0) {$+$};
\node at (4,0) {\begin{tikzpicture}
\node at (0,0) {\begin{ytableau} 1 & 1 \\ 0 & 0 \end{ytableau}};
\def\x{0.6}
\draw[line width=1mm] (-1*\x,-1*\x) -- (-1*\x,0*\x) -- (1*\x,0*\x) -- (1*\x,1*\x);
\end{tikzpicture}};
\node at (5,0) {$+$};
\node at (6,0) {\begin{tikzpicture}
\node at (0,0) {\begin{ytableau} 0 & 0 \\ 0 & 0 \end{ytableau}};
\def\x{0.6}
\draw[line width=1mm] (-1*\x,-1*\x) -- (-1*\x,1*\x) -- (1*\x,1*\x);
\end{tikzpicture}};
\node at (7,0) {$+$};
\node at (8,0) {\begin{tikzpicture}
\node at (0,0) {\begin{ytableau} 0 & 0 \\ 0 & 0 \end{ytableau}};
\def\x{0.6}
\draw[line width=1mm] (-1*\x,-1*\x) -- (-1*\x,1*\x) -- (1*\x,1*\x);
\end{tikzpicture}};
\end{tikzpicture}\]
\[\begin{tikzpicture}
\node at (-1.125,0) {$\Psi(\pi) =$};
\node at (0,0) {\begin{ytableau} 1 \\ 3 \end{ytableau}};
\node at (0.5,0) {$+$};
\node at (1,0) {\begin{ytableau} 1 \\ 4 \end{ytableau}};
\node at (1.5,0) {$+$};
\node at (2,0) {\begin{ytableau} 3 \\ 4 \end{ytableau}};
\node at (2.5,0) {$+$};
\node at (3,0) {\begin{ytableau} 3 \\ 4 \end{ytableau}};
\node at (3.675,0) {$=$};
\node at (5.25,0) {\begin{ytableau} 1 & 1& 3 & 3 \\ 3 & 4 & 4 & 4 \end{ytableau}};
\end{tikzpicture}\]

We end this section by describing how the operators $\mathrm{Pro}$, $\mathrm{Co}$, and $\mathrm{Tp}$ on plane partitions defined in Section~\ref{sec:intro} behave when viewed as operators on semistandard tableaux via the bijection~$\Psi$. Unsurprisingly, $\mathrm{Pro}$ behaves as~$\rho$ (thus justifying the name promotion), while $\mathrm{Co}$ behaves as evacuation. In order to record the behavior of transposition we need a little notation. So for $T \in \mathrm{SSYT}(m^n,2n)$ we define $T^\dagger\in \mathrm{SSYT}(m^n,2n)$ to be the tableau obtained from $T$ by first replacing each entry by~$2n+1$ minus that entry, and then replacing each column $I$ by its set-theoretic complement $\{1,2,\dots,2n\}\setminus I$.

\begin{Proposition} \label{prop:tableau_ops}
Via the bijection $\Psi\colon \mathrm{PP}^{m}(a\times b) \to \mathrm{SSYT}(m^a,a+b)$, view $\mathrm{Pro}$ and~$\mathrm{Co}$, and, if $a=b=n$, also $\mathrm{Tp}$, as operators on $\mathrm{SSYT}(m^a,a+b)$. Then:
\begin{itemize}\itemsep=0pt
\item $\mathrm{Pro}(T)=\rho(T)$ for all $T \in \mathrm{SSYT}(m^a,a+b)$,
\item $\mathrm{Co}(T)=\varepsilon(T)=\varepsilon^{*}(T)=T^+$ for all $T \in \mathrm{SSYT}(m^a,a+b)$,
\item $($if $a=b=n)$ $\mathrm{Tp}(T)=T^\dagger$ for all $T \in \mathrm{SSYT}(m^n,2n)$.
\end{itemize}
\end{Proposition}
\begin{proof}
The first bulleted item is immediate from the original description of $\Psi$ in terms of Gelfand--Tsetlin patterns, together with the description of the Bender--Knuth involutions as compositions of toggles which appears in Proposition~\ref{prop:bk_toggles} above.

The second and third bulleted items are easier to see from the alternate description of~$\Psi$. (Of course, with the second bulleted item we are implicitly applying the folklore Proposition~\ref{prop:ssyt_promo_basics}.) It is easily checked that the behaviors of $\mathrm{Co}$ and $\mathrm{Tp}$ are as claimed for single column tableaux. Furthermore, for $T \in \mathrm{SSYT}(m^a,a+b)$, if the columns of $T$ are $I_1,\dots, I_m$ then the columns of $T^+$ will be $I_m^{+},\dots, I_1^{+}$; and for $T \in \mathrm{SSYT}(m^n,2n)$, if the columns of $T$ are $I_1,\dots, I_m$ then the columns of $T^{\dagger}$ will be $I_1^{\dagger},\dots, I_m^{\dagger}$. Finally, we have $\pi \geq \pi' \iff \mathrm{Co}(\pi') \geq \mathrm{Co}(\pi)$ for any pair~$\pi, \pi' \in \mathrm{PP}^{1}(a\times b)$; and similarly we have $\pi \geq \pi' \iff \mathrm{Tp}(\pi) \geq \mathrm{Tp}(\pi')$ for any pair~$\pi, \pi' \in \mathrm{PP}^{1}(n \times n)$. These observations, together with the alternate description of $\Psi$, imply the second and third bulleted items.
\end{proof}

\subsection*{Acknowledgements} I thank Chris Fraser, Gabe Frieden, Vic Reiner, Brendon Rhoades, and Jessica Striker for useful discussions related to this work. I was supported by NSF grant \#1802920. I benefited from the use of Sage mathematics software~\cite{Sage-Combinat,sagemath} in the course of this research. Finally, I thank all of the anonymous referees for their careful attention to the document and useful comments.

\pdfbookmark[1]{References}{ref}
\LastPageEnding

\end{document}